\documentclass[12pt]{amsart}
\pdfoutput=1

\usepackage{geometry}

\usepackage{amsmath}
\usepackage{amsfonts}
\usepackage{amsthm}
\usepackage{amssymb}
\usepackage{graphicx}
\usepackage{tabularx}
\usepackage{overpic}
\usepackage{tikz}
\usetikzlibrary{arrows}
\usepackage{rotating}
\usepackage{hyperref}
\usepackage{framed}
\usepackage{xcolor}
\usepackage{placeins}

\newcommand{\x}{ {\bf x} }
\newcommand{\y}{ {\bf y} }
\newcommand{\z}{ {\bf z} }
\newcommand{\boldalpha}{ {\boldsymbol \alpha} }
\newcommand{\boldbeta}{ {\boldsymbol \beta} }
\newcommand{\moduli}{ \mathcal{M}^B(\x, \y, \vecrho) }
\newcommand{\pants}{\mathcal{P}}
\newcommand{\Y}{\mathcal{Y}_\pants}
\newcommand{\Ybar}{\bar{\mathcal{Y}}_\pants}

\newcommand{\Hselfgluer}{\mathcal{H}_{SG}}
\newcommand{\Yselfgluer}{\mathcal{Y}_{SG}}

\newcommand{\RR}{ \mathbb{R} }
\newcommand{\CC}{ \mathbb{C} }
\newcommand{\HH}{ \mathbb{H} }
\newcommand{\FF}{ \mathbb{F} }
\newcommand{\Z}{\mathbb{Z} }

\newcommand{\Alg}{ \mathcal{A} }
\newcommand{\Idem}{ \mathcal{I} }

\newcommand{\spinc}{ \mathfrak{s} }
\newcommand{\Ainfty}{ \mathcal{A}_\infty }

\newcommand{\vecrho}[1][ ]{ \vec\rho^{\>{#1}} }

\newcommand{\CFDD}{ \widehat{\mathit{CFDD}} }
\newcommand{\CFD}[1]{ \widehat{\mathit{CFD^{#1}}} }
\newcommand{\CFA}[1]{ \widehat{\mathit{CFA^{#1}}} }

\newcommand{\CFAA}{ \widehat{\mathit{CFAA}} }
\newcommand{\CFAAid}{ \CFAA(\mathbb{I}) }

\newcommand{\HFhat}{ \widehat{\mathit{HF}} }
\newcommand{\CFhat}{ \widehat{\mathit{CF}} }

\newcommand{\HeegDiag}{ \mathcal{H} }
\newcommand{\torusalgebra}{ \mathcal{A} }

\newcommand{\Gen}{\mathfrak{S}}

\DeclareMathOperator{\ind}{ind}
\DeclareMathOperator{\id}{id}

\newtheorem{thm}{Theorem}
\newtheorem{prop}{Proposition}[section]
\newtheorem{lem}[prop]{Lemma}

\theoremstyle{definition}
\newtheorem{definition}[prop]{Definition}

\theoremstyle{remark}
\newtheorem{remark}[prop]{Remark}

\title{Bordered Heegaard Floer homology and graph manifolds}
\author{Jonathan Hanselman  }
\thanks{The author was partially supported by NSF grant number DMS-0739392.}

\begin{document}

\maketitle

\begin{abstract}
We perform two explicit computations of bordered Heegaard Floer invariants. The first is the type $D$ trimodule associated to the trivial $S^1$-bundle over the pair of pants $\pants$. The second is a bimodule that is necessary for self-gluing, when two torus boundary components of a bordered manifold are glued to each other. Using the results of these two computations, we describe an algorithm for computing $\HFhat$ of any graph manifold.
\end{abstract}

\section{Introduction}

Heegaard Floer homology is a collection of invariants for closed 3-manifolds introduced by Ozsv{\'a}th and Szab{\'o} \cite{OzSz:first}. The package also contains invariants for 4-dimensional cobordisms and for knots and links \cite{OzSz:smooth4mfds, OzSz:knots, Rasmussen}. It has proved to be a sensitive invariant, but in general it is difficult to compute. The definition involves a chain complex whose generators are combinatorial but whose differential requires counting pseudo-holomorphic curves.

There are a few existing algorithms for computing Heegaard Floer homology. Sarkar and Wang developed a method using nice diagrams, for which computing the differential becomes combinatorial \cite{SarkarWang}; this method has since been refined and extended \cite{modifiedSarkarWang, OzSSz:nice_diagrams_nonhat, OzSSz:nice_diagrams_invariance}. Another approach uses grid diagrams and surgery formulas \cite{MOS:combinatorialHFK, MO:surgery, MOT:grid_diagrams}. A third algorithm is based on computing the bordered Heegaard Floer invariant for the surface diffeomorphism associated with a Heegaard splitting \cite{LOT:factoring}. Of these algorithms, only the third is practical enough to have been implemented on a computer, and in its current form it can only compute for Heegaard splittings of genus at most 2. This allows computations for many interesting manifolds, but ultimately the class of 3-manifolds which admit genus 2 Heegaard splittings is small. Efficiently computing Heegaard Floer homology for general 3-manifolds remains a difficult problem.

If we restrict to particular classes of 3-manifolds, computing Heegaard Floer homology becomes easier. For example, much is known about the Heegaard Floer homology of manifolds which are obtained by plumbing circle bundles according to a negative definite tree $\Gamma$. Ozsv{\'a}th and Szab{\'o} gave a combinatorial description of $HF^+$ of  these manifolds when the tree $\Gamma$ has at most one ``bad" vertex \cite{OzSz:plumbed}. This class of manifolds includes all Seifert fibered rational homology spheres. Their algorithm for computing $HF^+$ has been useful, for instance, in determining the existence of tight contact structures on Seifert fibered spaces \cite{LiscaStipsicz}. Nemethi introduced an invariant for negative definite plumbings, \emph{lattice homology}, which is combinatorially computable and conjecturally equivalent to $HF^+$ \cite{Nemethi}. Recent work has explored this conjectured equivalence; there a spectral sequence from lattice cohomology to $HF^+$, and they are known to be isomorphic for plumbings with at most two bad vertices \cite{OzSSz:knot1, OzSSz:knot2, OzSSz:spectral}.

There has been significant interest in understanding $L$-spaces, manifolds with minimal Heegaard Floer homology,  and the relationship between this condition and the existence of taut foliations and left orderability of the fundamental group. The conjectured relationship between these conditions is known to hold for particular classes of 3-manifolds, including Seifert fibered manifolds \cite{Peters, BoyerGordonWatson}. These geometric conditions can help us determine the $L$-space condition of a manifold even if we can not compute $\HFhat$ directly. Mauricio used lattice homology and the existence of taut foliations to give sufficient conditions on the weights of a negative definite tree $\Gamma$ under which a plumbing is or is not an $L$-space \cite{Mauricio}.

The plumbings of negative definite trees mentioned above are a special case of graph manifolds. A \emph{graph manifold} is a 3-manifold whose JSJ decomposition contains only Seifert fibered pieces. The non-Seifert fibered pieces in a JSJ decomposition are hyperbolic, so with respect to geometrization a graph manifold is a manifold with no hyperbolic pieces in its geometric decomposition. Thus graph manifolds represent an important subclass of 3-manifolds. For a brief overview of graph manifolds and their place in 3-manifold topology, see \cite{Neumann:overview}. In this paper we present a method for computing $\HFhat$ for any graph manifold, which is based on computing bordered Heegaard Floer invariants for certain fundamental building blocks from which graph manifolds can be constructed. 

This method finds a middle ground between the approaches mentioned above. It is more general than results restricted to negative definite plumbing trees, since it works for arbitrary graph manifolds. At the same time, it is more computationally practical than current algorithms for general 3-manfiolds. There is a computer implementation of this algorithm that is capable of handling quite complicated manifolds; it can been used, for instance, to see that the rank of $\HFhat$ of the graph manifold represented by the weighted tree in Figure \ref{huge_graph} is $213,\!312$.

\bigbreak

Let us recall some important facts and terminology concerning graph manifolds; we will follow the notation found in \cite{Neumann}. A graph manifold can be encoded by a decorated graph:

\begin{definition}
A \emph{connected closed plumbing graph} is a finite connected graph $\Gamma$ decorated as follows:
\begin{itemize}
\item
each vertex $i$ carries two integer weights $g_i$ and $e_i$;

\item
each edge carries a sign, $+$ or $-$.

\end{itemize}
We allow $\Gamma$ to have multiple edges connecting two vertices or edges connecting a vertex to itself.
\end{definition}

A connected closed plumbing graph $\Gamma$ specifies a (prime) graph manifold $M(\Gamma)$ as follows: For each vertex $i$ of $\Gamma$, let $d_i$ be the degree of the vertex. Let $F_i$ be the compact surface of genus $g_i$ with $d_i$ boundary components, where if $g_i < 0$ we mean that $F_i$ is nonorientable of genus $| g_i |$. Let $E_i$ be the circle bundle with orientable total space over $F_i$ with a chosen trivialization on the boundary and euler number $e_i$ (the euler number is well defined once the trivialization on the boundary is chosen). The edges of $\Gamma$ dictate how the $E_i$ are glued together. For each edge connecting vertices $i$ and $j$, a component $S^1 \times S^1$ of $\partial E_i$ is glued to a component $S^1 \times S^1$ of $\partial E_j$. The gluing always exchanges base and fiber directions; for $(+)$-edges the gluing map is  $\scriptsize \left( \begin{array}{cc}  0 & 1 \\  1 & 0  \end{array} \right) $, and for $(-)$-edges the gluing map is $\scriptsize \left( \begin{array}{cc}  0 & -1 \\  -1 & 0  \end{array} \right) $. In either case the gluing map is orientation reversing, and so $M(\Gamma)$ inherits consistent orientations from all of the $E_i$. For each edge connecting a vertex $i$ to itself, two components of $\partial E_i$ are glued with the appropriate gluing map.

Every prime graph manifold can be represented by a connected closed plumbing graph (for non-prime graph manifolds we allow disconnected graphs). The representation is not unique, but \cite[Section 4]{Neumann} gives a well developed calculus for manipulating plumbing graphs. In particular, plumbing graphs can be reduced to a normal form, and graphs of this form represent manifolds uniquely. A few additional facts are worth mentioning here:

\begin{itemize}
\item
Changing the sign of an edge often does not change the resulting 3-manifold. In fact, all that matters is the total sign around each loop in $\Gamma$. In particular, for acyclic graphs we may ignore the decoration on the edges.
\item
It is possible to represent any graph manifold with a plumbing graph such that no vertex is assigned a negative genus.
\item
We can describe graph manifolds with boundary by adding an additional weight $b_i$ to each vertex $i$. In the construction, $F_i$ is the genus $g_i$ surface with $b_i$ + $d_i$ boundary components. $E_i$ is the appropriate circle bundle over $F_i$, and $b_i$ components of $\partial E_i$ are not glued to anything.
\end{itemize}

Bordered Heegaard Floer homology is an extension of Heegaard Floer homology to manifolds with boundary \cite{LOT:Bordered}. Because graph manifolds decompose so nicely, bordered Heegaard Floer homology provides a natural approach for computing their $\HFhat$. The key ingredient is to compute the bordered invariants for arbitrary $S^1$-bundles over surfaces, the building blocks of graph manifolds. Changing the euler number of one of these bundles is equivalent to changing the parametrization of the boundary, which can be accomplished by tensoring with a well understood bimodule \cite[Section 10.2]{LOT:Bimodules}. As a result, we only need to compute invariants for trivial bundles over surfaces. As noted above, it is sufficient to consider bundles over orientable surfaces. Furthermore, any orientable surface has a pants decomposition---it can be obtained by gluing together copies of the pair of pants $\pants = S^3 \backslash \{ \text{three open disks} \}$. The trivial $S^1$-bundle over the surface can be obtained by gluing copies of the trivial $S^1$-bundle over $\pants$. Thus we see that the trivial bundle $\Y = \pants \times S^1$ is the fundamental building block for graph manifolds.

In Section \ref{sec:bordered} we will review the relevant background from bordered Heegaard Floer homology. 
The trimodule $\CFD{3}(\Y)$ will be explicitly computed in Section \ref{sec:trimodule}, proving the following:

\begin{thm}
\label{theorem1}
The summand of the type $D$ trimodule $\CFD{3}(\Y)$ in the middle spin$^c$-structure has five generators as a projective module: $v, w, x, y$, and $z$. Up to quasi-isomorphism, the differential is given by the following diagram:

\begin{center}
\begin{tikzpicture}[>=latex, scale=.8]

\node (y)  at (0,0) {$y$};
\node (v)  at (-5,0) {$v$};
\node (w) at (5,0) {$w$};
\node (x)  at (0,-3) {$x$};
\node (z)  at (0,3) {$z$};

\draw[->, bend left=10] (v) to node[pos = .7, above] {\scriptsize $\rho_1 \tau_{123} \sigma_3$} (y);
\draw[->, bend right=10] (v) to node[pos = .7, below] {\scriptsize $\rho_{123} \tau_{123} \sigma_{123}$}  (y);

\draw[->, bend right=8] (v) to node[below left = -2pt] {\scriptsize $\rho_3$} (x);
\draw[->, bend right=8] (x) to node[above right = -2pt] {\scriptsize $\rho_2 \sigma_{12}$}  (v);

\draw[->, bend left=10] (w) to node[pos = .6, below] {\scriptsize $\rho_{123} \tau_1 \sigma_{123}$} (y);
\draw[->, bend right=10] (w) to node[pos = .6, above] {\scriptsize $\rho_1 \tau_1 \sigma_3$}  (y);

\draw[->, bend right=8] (w) to node[above left = -3pt] {\scriptsize $\rho_3 \sigma_{12}$} (x);
\draw[->, bend right=8] (x) to node[below right = -2pt] {\scriptsize $\rho_2$}  (w);

\draw[->, bend right=10] (y) to node[left = -2pt] {\scriptsize $\tau_2 \sigma_2$} (x);
\draw[->, bend right=10] (x) to node[right = -2pt] {\scriptsize $\tau_3 \sigma_1$}  (y);

\draw[->, bend right=10] (z) to node[left = -2pt] {\scriptsize $\rho_3 \sigma_1$} (y);
\draw[->, bend right=10] (y) to node[right = -2pt] {\scriptsize $\rho_2 \sigma_2$}  (z);

\draw[->, bend left=8] (v) to node[above left = -2pt] {\scriptsize $\tau_3$} (z);
\draw[->, bend left=8] (z) to node[above right = -2pt] {\scriptsize $\tau_2$}  (w);
\end{tikzpicture}
\end{center}

\end{thm}

\noindent $\CFD{3}(\Y)$ in the other spin$^c$-structures will also be computed. For acyclic plumbing graphs with only genus 0 vertices, $\HFhat$ of the corresponding graph manifold can be obtained from the trimodule in Theorem \ref{theorem1} and bimodules for mapping classes of the torus. If the graph has a cycle or some vertex has genus $g_i > 0$, then an additional bimodule is needed. A bordered Heegaard diagram for this bimodule was described in \cite[Section 4.4]{LipshitzTreumann}, but the bimodule was not computed. In Section \ref{sec:self_gluer} we explicitly compute this bimodule, using the bordered Heegaard diagram $\Hselfgluer$ in Figure \ref{self_gluer}.

\begin{thm}
The bimodule $\CFDD(\Hselfgluer)$ in the middle spin$^c$-structure is given by Figure \ref{selfgluer_answer}. In the extremal spin$^c$-structures, it is quasi-isomorphic to zero.
\end{thm}

\noindent Finally, given an arbitrary plumbing graph, Section \ref{sec:HFhat} will describe the procedure for piecing together the relevant bordered invariants to obtain $\HFhat$ of the corresponding graph manifold.

\subsection{Acknowledgements}
I would like to thank Robert Lipshitz for suggesting this problem and for many helpful conversations about these computations. I am also grateful to Adam Levine, Peter Ozsv{\'a}th, Dylan Thurston, and Liam Watson for helpful conversations. I especially thank Adam Levine for suggesting a finger move isotopy which simplified the computation in Section \ref{sec:trimodule}.

\section{Background on Heegaard Floer homology}
\label{sec:bordered}

We begin by recalling the essential definitions and properties concerning the bordered Heegaard Floer invariants developed by Lipshitz, Ozsv{\'a}th, and Thurston. For a full treatment of these invariants, see \cite{LOT:Bordered}, \cite{LOT:Bimodules}. We discuss only the details that will be needed in the rest of the paper. In particular, we restrict to the case of manifolds with toroidal boundary components, which simplifies many of the definitions.

\subsection{Algebraic definitions}

Let $(\Alg, d)$ be a unital differential algebra over $\FF_2$, with a subring of idempotents $\Idem$, and let $\{ \iota_i \}$ be an orthogonal basis for $\Idem$, with ${\bf 1} = \sum \iota_i$.

A \emph{(left) type $D$ structure over} $\Alg$ is a vector space $N$ over $\FF_2$ with a left action of $\Idem$ and a map
\[ \delta_1 : N \rightarrow \Alg \otimes_\Idem N \]
satisfying the relation
\begin{equation}
\label{delta1relation}
 (\mu \otimes id_N) \circ (id_\Alg \otimes \delta_1) \circ \delta_1 + (d \otimes id_N) \circ \delta_1 = 0,
 \end{equation}
where $\mu : \Alg \otimes \Alg \rightarrow \Alg$ is multiplication on $\Alg$. The tensor product $\Alg \otimes_\Idem N$ is a left differential $\Alg$ module, with module structure $a \cdot (b \otimes x) = ab \otimes x$ and differential given by $\partial(a \otimes x) = a \cdot \delta_1(x) + d(a) \otimes x$. The relation (\ref{delta1relation}) ensures that $\partial^2 = 0$. Given the map $\delta_1$, define
\[ \delta_k : N \rightarrow \underbrace{\Alg \otimes_\Idem \cdots \otimes_\Idem \Alg}_\text{k times} \otimes_\Idem N \]
inductively by $\delta_0 = \id_N$ and $\delta_k = (\id_{\Alg^{\otimes k-1}} \otimes \delta_1) \circ \delta_{k-1}$ for $k > 0$. We say that the type $D$ structure $N$ is \emph{bounded} if $\delta_k = 0$ for all $k$ sufficiently large.

We will need to work with modules with multiple left actions. Let $\Alg_1, \cdots, \Alg_k$ be differential algebras, with rings of idempotents $\Idem_1, \ldots, \Idem_k$. A \emph{$k$-fold type $D$ structure over $\Alg_1, \ldots, \Alg_k$} is a type $D$ structure over $\Alg_1 \otimes \cdots \otimes \Alg_k$. We will call the module $\left( \Alg_1 \otimes \cdots \otimes \Alg_k \right) \otimes_{ (\Idem_1 \otimes \cdots \otimes \Idem_k) } N$ a \emph{type $D$ multimodule over $\Alg$}.

A \emph{(right) $\Ainfty$ module (or type $A$ structure) over $\Alg$} is a vector space $M$ over $\FF_2$ with a right action of $\Idem$ and maps
\[ m_{k+1}: M \otimes_\Idem \underbrace{\Alg \otimes_\Idem \cdots \otimes_\Idem \Alg}_\text{k times} \rightarrow M \]
satisyfing the folowing $\Ainfty$ relation for any $x \in M$ and any $a_1, \ldots a_n \in \Alg$:
\begin{equation}
\begin{aligned}
0 &= \sum_{i=0}^n m_{n-i+1}\big( m_{i+1}(x, a_1, \ldots, a_i) , a_{i+1}, \ldots, a_n \big) \\
&+ \sum_{i=1}^{n-1} m_{n}\big( x, a_1, \ldots, a_{i-1} , \mu(a_i, a_{i+1}), a_{i+2}, \ldots, a_n \big) \\
&+ \sum_{i=1}^n m_{n+1}\big( x, a_1, \ldots, a_{i-1}, d(a_i), a_{i+1}, \ldots, a_n \big).
\end{aligned}
\end{equation}
An informal statement of the $\Ainfty$ relations may be easier to remember: for any ordered set of inputs, the sum of all ways of combining those inputs using two steps (where each step is $\mu$, $d$, or some $m_i$) is zero. We also require that $m_2(x, {\bf 1}) = x$ and $m_k(x, \ldots, {\bf 1}, \ldots) =0$ for all $k > 2$. If $m_k = 0$ for all sufficiently large $k$, we say that the $\Ainfty$ module $M$ is \emph{bounded}.

More generally, we can define an \emph{$\Ainfty$ multimodule $M$ over $\Alg_1, \ldots, \Alg_k$} as follows: $M$ is a vector space over $\FF_2$ with a right action of $\Idem_1 \otimes \cdots \otimes \Idem_k$. $M$ is also equipped with maps
\[ m_{1, i_1, \ldots, i_k}: M \bigotimes_{\Idem_1 \otimes \cdots \otimes \Idem_k} \Alg_1^{\otimes i_1} \otimes \cdots \otimes \Alg_k^{\otimes i_k} \rightarrow M \]
satisfying an appropriate version of the $\Ainfty$ relation (we will generally suppress the subscripts on $m$ from the notation). To define the relation, we introduce the following functions. For $\vec{a}_\ell = (a_\ell^1, \ldots, a_\ell^k) \in \Alg_\ell^{\otimes k}$ and $0 \le j \le k$, define
\begin{equation*}
\begin{aligned}
T_j( \vec{a}_\ell ) :=& \left(a_\ell^1, \ldots, a_\ell^j \right)   \in \Alg_\ell^{\otimes j}\\
T^j( \vec{a}_\ell ) :=& \left(a_\ell^{j+1}, \ldots, a_\ell^k \right) \in \Alg_\ell^{\otimes k-j}\\
\bar{\mu} ( \vec{a}_\ell ) :=& \sum_{j=1}^{k-1} \left(a_\ell^1, \ldots, a_\ell^{j-1}, a_\ell^j a_\ell^{j+1}, a_\ell^{j+2}, \ldots, a_\ell^k \right)     \in \Alg_\ell^{\otimes k-1}\\
\bar{d}( \vec{a}_\ell ) :=& \sum_{j=1}^k \left(a_\ell^1, \ldots, a_\ell^{j-1}, d(a_\ell^j), a_\ell^{j+1}, \ldots, a_\ell^k \right) \in \Alg_\ell^{\otimes k}.
\end{aligned}
\end{equation*}
Now we can write down the $\Ainfty$ relation for multimodules. For any $x \in M$ and any $\vec{a}_1, \ldots, \vec{a}_k$ in $\Alg^{\otimes i_1}, \ldots, \Alg^{\otimes i_k}$, we have
\begin{equation}
\label{Ainfty_relation}
\begin{aligned}
0 &= \sum_{j_1, \ldots, j_k} m \! \left( m \! \left( x, T_{j_1}( \vec{a}_1), \ldots, T_{j_k} (\vec{a}_k) \right) , T^{j_1}( \vec{a}_1), \ldots, T^{j_k}( \vec{a}_k )  \right)\\
&+  \sum_{j=1}^k m \! \left( x, \vec{a}_1, \ldots, \vec{a}_{j-1}, \bar{\mu}(\vec{a}_j), \vec{a}_{j+1}, \ldots, \vec{a_k} \right) \\
&+ \sum_{j=1}^k m \! \left( x, \vec{a}_1, \ldots, \vec{a}_{j-1}, \bar{d}(\vec{a}_j), \vec{a}_{j+1}, \ldots, \vec{a_k} \right).
\end{aligned}
\end{equation}

It is possible to define combination multimodules, with some type $D$ actions and some type $A$ actions. Such a multimodule $N$ is equipped with maps
\[ \delta_1^{1,i_{k+1}, \ldots, i_\ell} : N \otimes \Alg_1^{\otimes i_{k+1}} \otimes \cdots \otimes \Alg_k^{\otimes i_\ell} \longrightarrow \Alg_{k+1} \otimes \cdots \otimes \Alg_{k+\ell} \otimes N \]
satisfying the appropriate versions of (\ref{delta1relation}) and (\ref{Ainfty_relation}). Type $DD$, $AA$, and $DA$ bimodules are discussed in \cite{LOT:Bimodules}, and the generalization to more algebra actions is straightforward.

If $M$ is an $\Ainfty$ module over $\Alg$ and $N$ is a type $D$ module over $\Alg$, and if at least one of them is bounded, we may define the \emph{box tensor product} $M \boxtimes N$ to be the vector space $M \otimes_\Idem N$ equipped with the differential
\[ \partial^\boxtimes (x \otimes y) = \sum_{k=0}^\infty (m_{k+1} \otimes \id_N )( x \otimes \delta_k(y)). \]
If $M$ is a multimodule over $\Alg_1, \cdots, \Alg_k$ such that the action of $\Alg_k$ is type $A$, and $N$ is a multimodule over $\Alg_k, \Alg_{k+1}, \ldots, \Alg_{k+\ell}$ such that the action of $\Alg_k$ is type $D$, and either $M$ or $N$ is bounded, then a box tensor product with respect to $\Alg_k$ can be defined in a similar way (see \cite[Section 2.3.2]{LOT:Bimodules} for the case when $N$ and $M$ are bimodules). $M \boxtimes_{\Alg_k} N$ is a multimodule over $\Alg_1, \ldots, \Alg_{k-1}, \Alg_{k+1}, \ldots, \Alg_{k+l}$, and the operations on $M \boxtimes N$ are determined by pairing operations on $M$ with sequences of operations in $N$ such that the $A_k$ outputs of the operations on $N$ match the $A_k$ inputs of the operation on $M$.

\begin{remark}
\label{edge_reduction}
We will often represent a $k$-fold type $D$ multimodule $M$ as a labeled, directed graph, where vertices correspond to the generators of $M$, and there is an arrow from $x_i$ to $x_j$ labeled by $a_{ij}$ if $a_{ij} \ne 0$ is the coefficient of $ x_j$ in $\partial(x_i)$. Here $a_{ij}$ is an element of $\Alg_1 \otimes \cdots \otimes \Alg_k$, the tensor product of $k$ copies of the torus algebra. We omit the edge label when $a_{ij} = 1$. We sometimes refer to an unlabeled arrow from $x_i$ to $x_j$ as a \emph{differential} from $x_i$ to $x_j$. Graphs with unlabeled edges can be simplified by a well known edge reduction algorithm \cite[Section 2.6]{Levine}: we eliminate the endpoints $x_i$ and $x_j$ of the unlabeled edge and all edges attached to these two vertices, and for each ``zig-zag"
$$x_k \overset{a_{kj}}\longrightarrow x_j \longleftarrow x_i \overset{a_{i\ell}}\longrightarrow x_\ell$$
we add an edge
$$x_k \overset{ a_{kj} a_{i\ell} }\longrightarrow x_\ell, $$
or if there is already an edge from $x_k$ to $x_\ell$ we add $a_{kj} a_{i\ell}$ to the label of that edge. The resulting graph represents a type $D$ multimodule that is quasi-isomorphic to $M$.
\end{remark}

\subsection{The Torus Algebra}
\bigbreak

To define bordered Heegaard Floer invariants, we associate a differential algebra to each boundary component of a 3-manifold with boundary. The algebra associated to the torus splits into a direct sum
\[ \Alg(T^2) = \Alg(T^2, -1) \oplus \Alg(T^2, 0) \oplus \Alg(T^2, 1). \]
$\Alg(T^2, -1)$ is $\FF_2$, and $\Alg(T^2, 1)$ is quasi-isomorphic to $\FF_2$, so we need only discuss $\Alg(T^2, 0)$.

The algebra $\Alg(T^2, 0)$ is generated as a vector space over $\FF_2$ by eight elements: two idempotents, $\iota_0$ and $\iota_1$, and six Reeb elements $\rho_1, \rho_2, \rho_3, \rho_{12}, \rho_{23}$, and $\rho_{123}$. The idempotents satisfy $\iota_i \iota_j = \delta_{ij} \iota_i$, and the identity element is ${\bf 1} = \iota_0 + \iota_1$. The Reeb elements interact with idempotents on either side as follows:
\[ \iota_0 \rho_1 = \rho_1 \iota_1 = \rho_1, \quad \iota_1 \rho_2 = \rho_2 \iota_0 = \rho_2, \quad \iota_0 \rho_3 = \rho_3 \iota_1 = \rho_3, \]
\[ \iota_0 \rho_{12} = \rho_{12} \iota_0 = \rho_{12}, \quad \iota_1 \rho_{23} = \rho_{23} \iota_1 = \rho_{23}, \quad \iota_0 \rho_{123} = \rho_{123} \iota_1 = \rho_{123}. \]
The only nonzero products of Reeb elements are
$\rho_1 \rho_2 = \rho_{12}$, $\rho_2 \rho_3 = \rho_{23}$, and ${ \rho_1 \rho_{23} = \rho_{12} \rho_3 = \rho_{123} }$. Although $\Alg(T^2)$ is a differential algebra, the differential on $\Alg(T^2, 0)$ is zero. For more on the torus algebra and how it arises in bordered Heegaard Floer homology, see \cite[Sec 11.1]{LOT:Bordered}.

\subsection{Bordered manifolds and bordered diagrams}

A bordered 3-manifold with $k$ torus boundary components is an oriented 3-manifold $Y$ with $\partial Y$ a disjoint union of $k$ tori $F_1, \ldots, F_k$, along with diffeomorphisms $\phi_i: T^2 \rightarrow F_i$. If $\phi_i$ is orientation reversing, then the corresponding boundary component is said to be type $D$; otherwise it is said to be type $A$. In this paper, we will deal almost exclusively with type $D$ boundaries.

A bordered 3-manifold can be represented by an arced bordered Heegaard diagram.

\begin{definition}
An \emph{arced bordered Heegaard diagram with k (torus) boundary components} is a quadruple $(\Sigma, {\boldsymbol \alpha}, {\boldsymbol \beta}, {\bf z})$, where

\begin{itemize}
\item $\Sigma$ is a compact surface of genus $g$ with $k$ boundary components;

\item ${\boldsymbol \alpha} = \{ \alpha_1^1, \alpha_2^1, \alpha_1^2, \alpha_2^2, \ldots, \alpha_1^k, \alpha_2^k, \alpha_1, \alpha_2, \ldots, \alpha_{g-k} \}$, where $\alpha_1^i$ and $\alpha_2^i$ are arcs embedded in $\Sigma$ with boundary on the $i$th component of $\partial \Sigma$ and $\alpha_j$ is an embedded circle in $\Sigma$, and the $\alpha$ circles/arcs are pairwise disjoint;

\item ${\boldsymbol \beta}$ is $g$-tuple of disjoint circles in $\Sigma$;

\item ${\bf z}$ is a basepoint $z$ in $\Sigma \backslash ({\boldsymbol \alpha} \cup {\boldsymbol \beta})$ together with arcs in $\Sigma \backslash ({\boldsymbol \alpha} \cup {\boldsymbol \beta})$ connecting $z$ to each boundary component of $\Sigma$.
\end{itemize}
We also require that $\boldalpha$ and $\boldbeta$ intersect transversely and $\Sigma \backslash \boldalpha$ and $\Sigma \backslash \boldbeta$ are connected.
\end{definition}

An arced bordered Heegaard diagram gives rise to a bordered 3-manifold by attaching 2-handles to a thickened version of the Heegaard surface $\Sigma$. The one and two boundary cases are described in Constructions 5.3 and 5.6 of \cite{LOT:Bimodules}, and the construction for more boundary components is completely analogous.

\begin{figure}

\begin{center}
Type $A$ boundary \qquad \quad Type $D$ boundary 

\begin{overpic}[scale = 1]{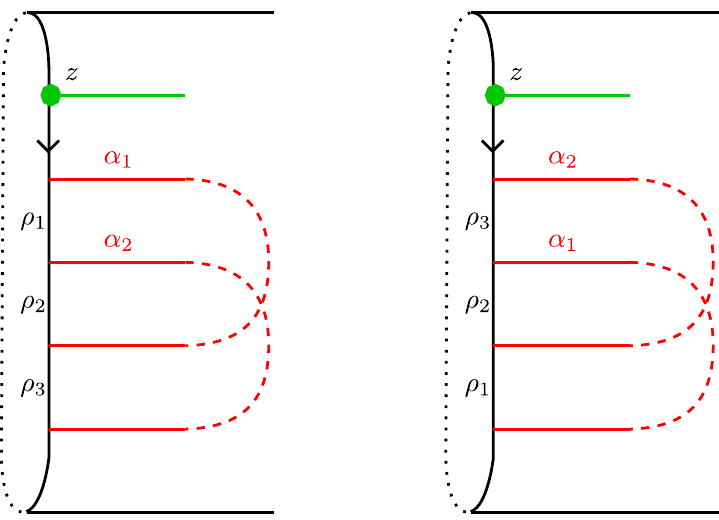}

\end{overpic}
\end{center}

\caption{Boundary markings for type $A$ and $D$ boundaries on a bordered Heegaard diagram.}
\label{boundary_labels}
\end{figure}

To define bordered invariants, we will also need to equip a bordered Heegaard diagram with labels on the boundary, as in Figure \ref{boundary_labels}. Each component of $\partial \Sigma$ is divided into four segments by the arcs $\alpha^i_1$ and $\alpha^i_2$, with one containing a basepoint, an endpoint of an arc in ${\bf z}$. Progressing from the basepointed segment in the direction which agrees with the boundary orientation on $\partial \Sigma$, we label the three remaining segments on the $i$th boundary component by $\rho_1^i, \rho_2^i$, and $\rho_3^i$ for type $A$ boundaries, or by $\rho_3^i, \rho_2^i,$ and $\rho_1^i$ for type $D$ boundaries. In each case, $\rho^i_{12}, \rho^i_{23}$, and $\rho^i_{123}$ refer to the relevant concatenations. We call these oriented arcs \emph{Reeb chords}. The assumption that $\Sigma \backslash \boldalpha$ is connected implies that the endpoints of $\alpha^i_1$ and $\alpha^i_2$ alternate. We assume that the first endpoint after the basepoint (following the boundary orientation) is $\alpha^i_1$ for type $A$ boundaries and $\alpha^i_2$ for type $D$ boundaries.

We can associate a copy of the torus algebra to each boundary component, so that the Reeb chords on $\partial \Sigma$ correspond directly to the Reeb elements of the algebra. By abuse of notation, we often use $\rho^i_I$ to refer both to the Reeb chord on $\partial \Sigma$ and the corresponding algebra element in the corresponding copy of the torus algebra.

\subsection{Type $D$ invariants}

Let $Y$ be a borderd 3-manifold with $k$ boundary components, and let $\HeegDiag$ be an arced bordered Heegaard diagram representing $Y$ which is provincially admissible in the sense of \cite[Definition 4.23]{LOT:Bordered}. Choose a complex structure $J$ on $\Sigma \times [0,1] \times \RR$. To ensure transversality, the choice of $J$ must be generic; however, for the computations in this paper we may assume that $J$ splits as $J_\Sigma \times J_D$, where $J_\Sigma$ is a generic complex structure on $\Sigma$ and $J_D$ is a generic complex structure on $[0,1] \times \RR$. Split complex structures provide enough flexibility for transversality when the projections to $\Sigma$ of all curves being considered are somewhere injective \cite[Proposition 3.9]{Lipshitz:cylindrical}. Given these choices, we will associate to $\HeegDiag$ a type $D$ multimodule $\CFD{k} (\HeegDiag, J)$ over $k$ copies of the torus algebra $\Alg(T^2)$. We will often suppress $J$ from the notation.

Let $\Gen(\HeegDiag)$ be the set of unordered $g$-tuples ${\bf x} = (x_1, \ldots, x_g)$ which contain exactly one point on each $\beta$ curve, exactly one point on each $\alpha$ curve, and at most one point on each $\alpha$ arc. Elements of $\Gen(\HeegDiag)$ fall into different spin$^c$-structures according to how many $\alpha$ arcs are occupied on each boundary. As a vector space over $\FF_2$, $\CFD{k} (\HeegDiag)$ is generated by $\Gen (\HeegDiag)$, and it splits as a direct sum over spin$^c$-structures on $Y$ \cite[Lemma 4.21]{LOT:Bordered}. Each generator $\x \in \Gen(\HeegDiag)$ comes equipped with an idempotent in the algebra associated to each boundary component; if $\x$ has exactly one $\alpha^i$ arc occupied, then the corresponding idempotent in $\Alg_i = \Alg(T^2)$ is $\iota^i_1$ if $\x$ contains a point on $\alpha^i_1$ and $\iota^i_0$ if $\x$ contains a point on $\alpha^i_2$.

The differential on $\CFD{k}(\HeegDiag)$ counts $J$-holomorphic curves in $\Sigma \times [0,1] \times \RR$ with appropriate boundary conditions (for precise statements of these conditions, see \cite[Section 5.2]{LOT:Bordered}). These curves can be sorted into relative homology classes. For any $\x, \y \in \Gen( \HeegDiag )$, let $\pi_2(\x, \y)$ denote the set of homology classes of curves in  $\Sigma \times [0,1] \times \RR$ with boundary conditions consistent with a differential connecting $\x$ to $\y$. Computing the differential involves counting the holomorphic representatives for each homology class.

Under the projection $\Sigma \times [0,1] \times \RR \rightarrow \Sigma$, a homology class $B \in \pi_2(\x, \y)$ projects to an element of $H_2(\Sigma, \boldalpha \cup \boldbeta \cup \partial \Sigma)$, called the \emph{domain} of $B$. $B$ is determined by its domain. A domain is a linear combination of components of $\Sigma \backslash (\boldalpha \cup \boldbeta)$, which we call \emph{regions}. Furthermore, the domain of any $B \in \pi_2(\x, \y)$ must satisfy the following conditions:

\begin{itemize}
\item the multiplicity of the region containing the basepoint $z$ is 0;
\item at each $p \in \boldalpha \cap \boldbeta$, let $n_1(p), \ldots, n_4(p)$ be the multiplicities of the four regions with corners at $p$, counting counterclockwise starting from an $\alpha$. Then
\begin{equation}
\label{corner_condition}
n_1(p) - n_2(p) + n_3(p) - n_4(p) = \begin{cases}
1 & p \in \x \backslash \y  \\
-1 & p \in \y \backslash \x  \\
0 & \text{else.}
\end{cases}
\end{equation}
\end{itemize}

A domain is called \emph{positive} if every region has non-negative multiplicity. Only positive domains can support holomorphic representatives. Because the Heegaard diagram $\HeegDiag$ is provincially admissible, there are a finite number of positive domains with multiplicity at most $1$ in the regions adjacent to $\partial \Sigma$ (we will see that only these are relevant for computing $\CFD{k}$). Finding them is a simple matter of linear algebra.

In addition to its domain $B$, a holomorphic curve that contributes to the differential of $\x$ also specifies a sequence of Reeb chords $\vecrho = (\vecrho[1], \ldots, \vecrho[k] )$, where $\vecrho[i] = (\rho^i_{I_1}, \ldots, \rho^i_{I_n} )$ are sequences corresponding to each boundary component of $\Sigma$. For each boundary component with one $\alpha^i$ arc occupied by $\x$, the pair $(B, \vecrho[i])$ will satisfy the following conditions:

\begin{itemize}
\item the initial point (with respect to the boundary orientation) of $\rho^i_{I_1}$ lies on the same $\alpha^i$ arc as $\x$;
\item for each $m > 1$, the initial point of $\rho^i_{I_m}$ lies on the same $\alpha^i$ arc as the terminal point of $\rho^i_{I_{m-1}}$.
\end{itemize}
A pair $(B, \vecrho[i])$ satisfying the above conditions is called \emph{strongly boundary monotonic}. For each boundary component with zero or two $\alpha^i$ arcs occupied by $\x$, we may assume that $\vecrho[i] = ()$. The pair $(B, \vecrho)$ coming from a holomorphic curve will also satisfy the following property:
\begin{itemize}
\item the intersection of $B$ with the $i$th component of $\partial \Sigma$ is equal to the sum of the Reeb chords in $\vecrho[i]$ as elements of $H_1(\partial \Sigma, \boldalpha \cap \partial \Sigma)$.
\end{itemize}
We say that the pair $(B, \vecrho)$ is \emph{compatible} if it satisfies this condition and each $(B, \vecrho[i])$ is strongly boundary monotonic (compare \cite[Definition 5.61]{LOT:Bordered}).

Given generators $\x, \y \in \Gen( \HeegDiag )$, a homology class $B \in \pi_2(\x, \y)$, and a sequence of Reeb chords $\vecrho$ such that $(B, \vecrho)$ is compatible, we can define $\mathcal{M}^B(\x, \y, \vecrho)$ to be the moduli space of $J$-holomorphic curves in $\Sigma \times [0,1] \times \RR$ with domain $B$ and whose asymptotics specify the initial generator $\x$, the final generator $\y$, and the sequence of Reeb chords $\vecrho$ (for the full definition, see \cite[Section 5]{LOT:Bordered}).

The dimension of the moduli space $\mathcal{M}^B(\x, \y, \vecrho)$ is one less than the index $\ind(B, \vecrho)$, defined in \cite[Definition 5.61]{LOT:Bordered}. In the special case of toroidal boundary, the index is given by
\begin{equation}
\label{index1}
\ind(B, \vecrho) = e(B) + n_\x(B) + n_\y(B) +\sum_{1 \le i \le k} \left[ \frac{1}{2} \left| \vecrho[i] \right| +\sum_{j < l} L(\rho_{I_j}^i, \rho_{I_l}^i) \right],
\end{equation}
where
\[ e(B) = \chi(B) -\frac{ \text{\# acute corners} }{4} + \frac{ \text{\# obtuse corners} }{4} \]
is the Euler measure of $B$, $n_\x(B)$ (respectively $n_\y(B)$) is the sum over $x_i \in \x$ (respectively $y_i \in \y$) of the average multiplicity in $B$ of the four regions incident to $x_i$ (respectively $y_i$), $\left| \vecrho[i] \right|$ is the number of Reeb chords in the part of $\vecrho$ associated to the $i$th component of $\partial \Sigma$, and $L(\rho_{I_j}, \rho_{I_l})$ is a linking term for Reeb chords defined as follows:
\begin{eqnarray*}
L(\rho_1, \rho_2) = L(\rho_2, \rho_3) = L(\rho_{12}, \rho_3) = L(\rho_1, \rho_{23}) &=& \tfrac{1}{2} \\
L(\rho_2, \rho_1) = L(\rho_3, \rho_2) = L(\rho_3, \rho_{12}) = L(\rho_{23}, \rho_1) &=& -\tfrac{1}{2} \\
L(\rho_{12}, \rho_{23}) &=& 1 \\
L(\rho_{23}, \rho_{12}) &=& -1 \\
L(\rho_{I_j}, \rho_{I_l}) &=& 0 \text{ for all other pairs of } I_j \text{ and } I_l.
\end{eqnarray*}
The differential counts J-holomorphic curves in moduli spaces with dimension 0, so we only need to consider domains and Reeb chords with $\ind(B, \vecrho) = 1$.

To define the differential we need one more piece of notation. If $\rho^i_{I}$ represents a Reeb chord on the $i$th boundary component of $\Sigma$, let $a(\rho^i_{I})$ denote the corresponding element of $\torusalgebra_i$, the copy of the torus algebra associated to the $i$th boundary component. If $\vecrho[i] = (\rho^i_{I_1}, \ldots, \rho^i_{I_n} )$ is a sequence of Reeb chords on the $i$th boundary, let $a( \vecrho[i] )$ denote the element $\rho^i_{I_1} \rho^i_{I_2} \cdots \rho^i_{I_n} \in \torusalgebra_i$, and if $\vecrho = (\vecrho[1], \ldots, \vecrho[k] )$, let $a(\vecrho)$ denote the element $a(\vecrho[1] ) \otimes \cdots \otimes a(\vecrho[k]) \in \torusalgebra_1 \otimes \cdots \otimes \torusalgebra_k$. We now define the differential on $\CFD{k}$ as follows: For any $\x \in \Gen( \HeegDiag )$,
\[ \partial(\x) = \sum_{\y \in \Gen( \HeegDiag )}       \sum_{B \in \pi_2(\x, \y)}    \!\!    \sum_{\substack{ \{ \vecrho | \ind(B, \vecrho) = 1, \\ (B, \vecrho) \text{is compatible} \} }}           \!\!\!\!\!\!\!               \Big( \# \mathcal{M}^B(\x, \y, \vecrho) \Big) a(\vecrho) \otimes \y , \]
where the count of a moduli space is taken mod 2. 

Note that the multimodule $\CFD{k}(\HeegDiag, J)$ depends on the choices of $\HeegDiag$ and $J$. However, its quasi-isomorphism type is an invariant of the bordered manifold $Y$. We denote this quasi-isomorphism class by $\CFD{k}(Y)$. We will deviate slightly from the notation introduced here when $k \le 2$ in order to agree with existing notation. That is, we will omit the superscript in $\CFD{1}$, and we will write $\CFDD$ instead of $\CFD{2}$.

\subsection{Type $A$ Invariants}

Let $Y$ be a bordered 3-manifold with $k$ boundary components and let $\HeegDiag$ be a provincially admissible arced bordered Heegaard diagram representing $Y$ and $J$ a chosen complex structure, as before. We can define a type $A$ multimodule over $k$ copies of the torus algebra, denoted $\CFA{k} (\HeegDiag)$. In this paper, we will never need to compute $\CFA{k}(\HeegDiag)$. However, as a computational trick we will make use of the relationship between $\CFA{k}(\HeegDiag)$ and $\CFD{k}(\HeegDiag)$, so it will be helpful to state the definition.

$\CFA{k}(\HeegDiag)$ is generated by the same set $\Gen(\HeegDiag)$ that generates $\CFD{k}(\HeegDiag)$. The differential and higher multiplications are defined by counting the same $J$-holomorphic curves that appear in the definition of $\CFD{k}(\HeegDiag)$. We will assume for the sake of comparison that the Reeb chords on the boundary are labeled the same as if we were computing $\CFD{k}(\HeegDiag)$, so that for a given domain the compatible sequences of Reeb chords $\vec\rho$ and the moduli spaces $\mathcal{M}^B(\x, \y, \vecrho)$ are exactly the same. However, with this convention we must change the algebra elements in the $\Ainfty$ operation, since normally the Reeb chords are labeled in the opposite order for type $A$ Heegaard diagrams. Let the function $\bar{a}$ be the same as $a$ except that it also interchanges $\rho^i_1$ with $\rho^i_3$ and $\rho^i_{12}$ with $\rho^i_{23}$. Then given a generator $\x \in \Gen(\HeegDiag)$ and sequence of Reeb chords $\vecrho = (\vecrho[1], \ldots, \vecrho[k] )$, 
\[ m( \x, \bar{a}(\vecrho[1]), \ldots, \bar{a}(\vecrho[k]) ) = \sum_{\y \in \Gen( \HeegDiag )}   \!\!\!\!\!     \sum_{\substack{ \{ B \in \pi_2(\x, \y) | \\ \ind(B, \vecrho) = 1, \\  (B, \vecrho) \text{ compatible} \} }}       \!\!\!\!\!\! \!       \Big( \# \mathcal{M}^B(\x, \y, \vecrho) \Big) \y, \]
where we think of $\bar{a}(\vecrho[i])$ as an element of $\torusalgebra_i^{\otimes | \vecrho[i] |}$, and where the moduli space counts are taken mod 2.

\subsection{Tensor Products and the Pairing Theorem}

For a bordered manifold with many boundary components, we can define bordered invariants which are type $D$ with respect to some boundaries and type $A$ with respect to others. These invariants can be obtained from $\CFD{k}$ by taking the box tensor product with the bimodule $\widehat{\mathit{CFAA}}(\mathbb{I})$, which can be found in \cite[Figure 21]{LOT:Bimodules}. An alternative shorthand algorithm for converting to type $D$ boundaries to type $A$ is described in \cite[Section 2.3]{splicing}.

Bordered invariants satisfy a pairing theorem \cite[Theorem 11]{LOT:Bimodules}. Given a bordered invariant for $Y_1$ which is type $A$ with respect to the $i$th boundary component and a bordered invariant for $Y_2$ which is type $D$ with respect to the $j$th boundary component, we can compute the box tensor product of the two multimodules with respect to the corresponding copies of the torus algebra, assuming the modules are appropriately bounded. The pairing theorem states that up to $\Ainfty$-homotopy equivalence, the result is the bordered invariant for the manifold obtained by gluing the $i$th boundary of $Y_1$ to the $j$th boundary of $Y_2$.

In this paper, we will work primarily with type $D$ modules, and convert only one boundary component at a time to type $A$ in order to tensor with another type $D$ module.

\subsection{Useful Results for Computation} This section collects a handful of results that are useful when explicitly computing a type $D$ bordered invariant.

The first is a slight rephrasing of \cite[Proposition 2.1]{Levine}:

\begin{prop}
\label{idempotent_proposition}
\quad

(a) For a given boundary component, the only non-empty sequences of Reeb chords which can contribute nonzero terms to the differential in $\CFD{k}$ are $(\rho^i_1)$, $(\rho^i_2)$, $(\rho^i_3)$, $(\rho^i_1, \rho^i_2)$, $(\rho^i_2, \rho^i_3)$, $(\rho^i_1, \rho^i_2, \rho^i_3)$, and $(\rho^i_{123})$.

(b) Furthermore, if $B \in \pi_2(\x, \y)$ contributes with $\vecrho[i] = (\rho^i_2)$ or $\vecrho[i] = (\rho^i_1, \rho^i_2)$, then $\y$ contains a point on $\alpha^i_2$. If $B$ contributes and $\vecrho[i]$ is $(\rho^i_1)$, $(\rho^i_3)$, $(\rho^i_{123})$, or $(\rho^i_2, \rho^i_3)$,  then $\y$ contains a point on $\alpha^i_1$.

\end{prop}

In particular, this proposition implies that only domains with multiplicity 0 or 1 in every region that intersects $\partial \Sigma$ can contribute nontrivially to the differential in $\CFD{k}$. For provincially admissible Heegaard diagrams this ensures that there is a finite number of positive domains to consider.

Another implication of Proposition \ref{idempotent_proposition} is that Equation (\ref{index1}) can be simplified for type $D$ computations.

\begin{lem}
\label{index_lemma}
If the pair $(B, \vecrho)$ contributes a nonzero term to the differential of $\CFD{k}$, then the index of the pair is given by
\begin{equation}
\label{index2}
\ind(B, \vecrho) = \ind(B) = e(B) + n_\x(B) + n_\y(B) + \frac{\# \{ Z \in \pi_0(\partial\Sigma) | Z \cap B \ne \emptyset\} }{2}.
\end{equation}
In particular the index depends only on $B$.
\begin{proof}
We examine the term in brackets in Eq. \ref{index1}. For the $i$th component of $\partial\Sigma$, there is a contribution to the index of \[ \left[ \frac{1}{2} \left| \vecrho[i] \right| +\sum_{j < l} L(\rho_{I_j}^i, \rho_{I_l}^i) \right]. \]
We can evaluate this term for each of the sequences of Reeb chords allowed by Proposition \ref{idempotent_proposition}. If $\vecrho[i]$ is $(\rho^i_1)$, $(\rho^i_2)$, $(\rho^i_3)$, or $(\rho^i_{123})$, then $\left| \vecrho[i] \right| = 1$ and there are no linking terms. If $\vecrho[i]$ is $(\rho^i_1, \rho^i_2)$ or $(\rho^i_2, \rho^i_3)$, then $\left| \vecrho[i] \right| = 2$, and there is one linking term, with a value of $-\tfrac{1}{2}$. If $\vecrho[i]$ is $(\rho^i_1, \rho^i_2, \rho^i_3)$, then $\left| \vecrho[i] \right| = 3$, and the two nonzero linking terms $L(\rho^i_1, \rho^i_2)$ and $L(\rho^i_2, \rho^i_3)$ evaluate to $-\tfrac{1}{2}$. In any of these cases, the total contribution to the index is $\tfrac{1}{2}$. The only other possibility is that $\vecrho[i] = ()$, which happens when $B$ does not contain any regions adjacent to the $i$th boundary component of $\Sigma$. In this case, the contribution of $\vecrho[i]$ to the index is 0. Summing over all boundary components yields

\[ \sum_{1 \le i \le k} \left[ \frac{1}{2} \left| \vecrho[i] \right| +\sum_{j < l} L(\rho_{I_j}^i, \rho_{I_l}^i) \right] = \frac{\# \{ Z \in \pi_0(\partial\Sigma) | Z \cap B \ne \emptyset\} }{2}. \]

\end{proof}
\end{lem}

Lemma \ref{index_lemma} allows us to exclude a domain $B$ from consideration in computing $\CFD{k}$ if $\ind(B) \ne 1$, without needing to consider all sequences of Reeb chords compatible with $B$. 

In practice, computing $\CFD{k}$ from a Heegaard diagram begins by writing down all positive domains $B \in \pi_2(\x, \y)$ for each pair of generators $\x$ and $\y$, and then eliminating as many domains as possible using Proposition \ref{idempotent_proposition} and Lemma \ref{index_lemma}. At some point, however, it is necessary to prove that a given domain/Reeb chord pair \emph{does} contribute to the differential. The following proposition asserts that a domain which can be realized as an immersed polygon always contributes.

\begin{prop}
\label{polygons}
Let $P$ be a $2n$-gon, with edges numbered consecutively, and suppose that there is map $P \xrightarrow{u} \Sigma$  satisfying the following conditions:
\begin{itemize}
\item
$u|_{\partial P}$ takes even edges of $P$ to $\boldbeta$, odd edges of $P$ to $\boldalpha \cup \partial\Sigma$, and corners to acute corners;

\item
$u$ is an immersion, except at the preimages of $\boldalpha \cap \partial \Sigma$;

\item
for each boundary component of $\Sigma$, at most one edge of $P$ maps to $\alpha_1^i \cup \alpha_2^i$, and for each $\beta \in \boldbeta$, at most one edge of $P$ maps to $\beta$.

\end{itemize}
The image of $u$ covers each region in $\Sigma$ with a certain multiplicity; let $B(u)$ be the corresponding positive domain. The image of $\partial P$ determines a sequence of Reeb chords $\vecrho(u)$, with the chords in the image of a single edge ordered according to the boundary orientation on $\partial P$. If $B(u) \in \pi_2(\x, \y)$ for some generators $\x$ and $\y$ in the middle spin$^c$-structure, then $(B(u), \vecrho(u))$ is compatible, and $\Big( \# \mathcal{M}^B(\x, \y, \vecrho) \Big) \equiv 1$ (mod 2).

\begin{proof}
A holomorphic curve in $\Sigma \times [0,1] \times \RR$ is equivalent to a holomorphic map of a Riemann surface with boundary into $\Sigma$ along with a branched covering map of that surface over the unit disk $D^2 \subset \CC$ (see \cite[Lemma 3.6]{OzSz:first}). For a specific domain, we look at Riemann surfaces which map onto the given domain in $\Sigma$, such that the preimages of the $\alpha$ arcs (together with boundary Reeb chords) and $\beta$ arcs map to the right and left boundaries, respectively, in the projection to $D^2$, and the preimages of the $\x$ and $\y$ corners map to $-i$ and $i$, respectively.

In this case, we already have a map from the polygon $P$ to $\Sigma$. There is a unique choice of complex structure on $P$ that makes $u$ holomorphic (induced by pulling back the complex structure on $\Sigma$). So we need to show that with this fixed complex structure, there is a unique $n$-fold branched covering map to $D^2$ up to an $\RR$ action.

First choose a biholomorphic map from $P$ to the upper half plane $\HH$, which takes one of the $\y$ corners to $\infty$, and the other corners to points $x_1, x_2, \ldots, x_{2n-1}$ along the real axis. We now want to find a degree $n$ map $\HH \rightarrow \HH$ which takes $x_i$ to 0 for $i$ odd and to $\infty$ for $i$ even, and takes $\infty$ to $\infty$. Such a map is given by
\[ z \rightarrow \frac{ (z-x_1)(z-x_3)\cdots(z-x_{2n-1}) }{ (z-x_2)(z-x_4)\cdots(z-x_{2n-2}) }. \]
This map is unique up to scaling. Finally we can find a biholomorphic map from $\HH$ to $D^2$ which takes $0$ to $-i$ and $\infty$ to $+i$. Composing these three maps gives the desired $k$-fold branched cover $P \rightarrow D^2$.
\end{proof}

\end{prop}

\begin{figure}
\centering
\bigbreak
\includegraphics[scale = .8]{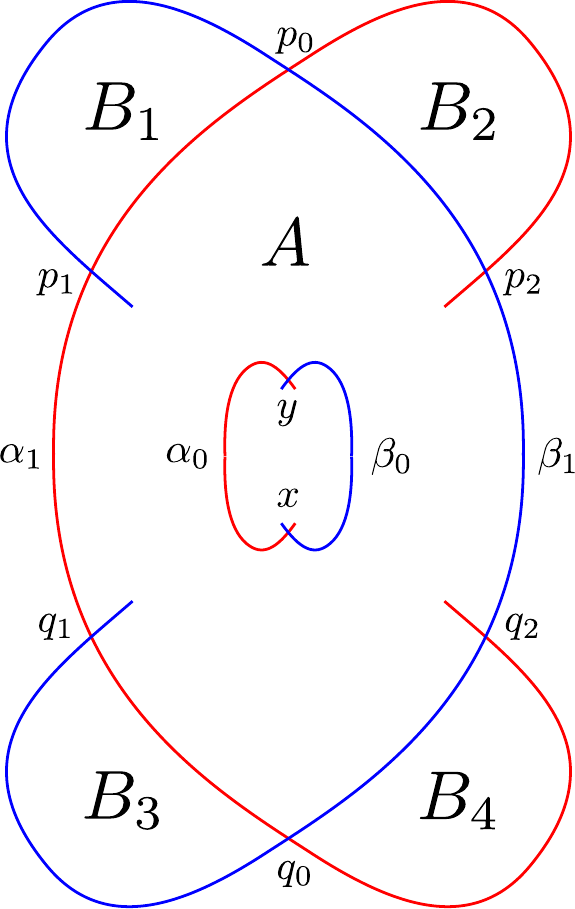}
\caption{}
\label{index0domain}
\end{figure}

Another common situation in which the moduli space of holomorphic curves can be understood is pictured in Figure \ref{index0domain}. The following is \cite[Lemma 3.2]{Levine}, but we recall the proof here in order to introduce notation and reasoning that will be useful later. 

\begin{prop}
\label{Levine_prop}
Suppose a Heegaard diagram contains an annulus $A$ as in Figure \ref{index0domain} and one or more of the bigons $B_1, \ldots, B_4$, where $\alpha$ arcs may contain segments of $\partial\Sigma$, and where the ends of $\alpha_1$ and $\beta_1$ leave $A$ through the opposite boundary component. Let $D_i$ denote the domain corresponding to the union of $A$ and $B_i$. Then either $D_1$ and $D_3$ count toward the differential and $D_2$ and $D_4$ do not, or vice versa, depending on the choice of complex structure $J$ on $\Sigma$.

\begin{proof}
Let $A_r$ denote the standard annulus, $S^1 \times [0,r]$ with a fixed complex structure. For a unique positive number $r$ there is a holomorphic map $u: A_r \rightarrow A$ taking $S^1 \times \{0\}$ to the inner boundary of $D$, $D \cap ( \alpha_0 \cup \beta_0)$, and taking $S^1 \times \{r\}$ to the outer boundary $D \cap (\alpha_1 \cup \beta_1)$. This map is unique up to rotation in the $S^1$ factor. Let $a_0$ and $b_0$ denote the inverse images in $S^1 \times \{0\}$ of $\alpha_0$ and $\beta_0$, respectively. Let $a_1$ and $b_1$ denote the respective inverse images of $\alpha_1$ and $\beta_1$ in $S^1 \times \{r\}$. Define $\Theta_A^{x,y}$ to be $l(a_0) / l(b_0)$, the ratio of the lengths of the preimages of the $\alpha$ and $\beta$ arcs on the boundary of $A$ which contains $x$ and $y$. Similarly, define $\Theta_A^{p_0, q_0}$ to be $l(a_1) / l(b_1)$. The domain $A$ will have a holomorphic representative if $\Theta_A^{x,y} = \Theta_A^{p_0, q_0}$ \cite[Lemma 9.3]{OzSz:first}, but for a generic choice of complex structure this will not be the case.

Now consider the domain $D_1$. It is an annulus with one obtuse corner at $p_1$. There is a one parameter family of conformal structures depending on how much we cut into the annulus along the $\alpha$ or $\beta$ arc at the obtuse corner. We specify the length of this cut by a parameter $c$, where $c <0$ corresponds to cutting along $\alpha_1$ and $c>0$ corresponds to cutting along $\beta_1$. The cut approaches $p_0$ as $c \to -\infty$ and it approaches $\alpha_0$ as $c \to \infty$. For any value of $c$ there is a holomorphic map $U^c: A_r \rightarrow D_1$, unique up to rotation in the $S^1$ factor of $A_r$. We can define $\theta_{D_1}^{x, y}(c)$ and $\theta_{D_1}^{p_1, q_0}(c)$ analogously to $\Theta_A^{x,y}$ and $\Theta_A^{p_0, q_0}$, as the ratio of the lengths of the $\alpha$ and $\beta$ components on the corresponding boundary of $D_1$.

As the cutting parameter varies, $D_1$ will have a holomorphic representative each time $\theta_{D_1}^{x, y}(c) =\theta_{D_1}^{p_1, q_0}(c)$, and thus the number of holomorphic representatives is determined by the number of zeros of $\theta_{D_1}^{x, y} - \theta_{D_1}^{p_1, q_0}$. The mod 2 count of these zeros is determined by the end behavior of $\theta_{D_1}^{x, y} - \theta_{D_1}^{p_1, q_0}$. As $c$ approaches $\infty$, the cut along $\beta_1$ from $p_1$ approaches $\alpha_0$. In this limit $\theta_{D_1}^{x, y}$ becomes very large and $\theta_{D_1}^{p_1, q_0}$ becomes very small, so $\theta_{D_1}^{x, y} - \theta_{D_1}^{p_1, q_0} \to +\infty$.  In the other extreme, we cut along $\alpha_1$ from $p_1$ to $p_0$. The limit is a broken flow where the bigon $B_1$ is pinched off from the annulus $A$. In this limit $\theta_{D_1}^{x, y} = \Theta_{A}^{x, y}$ and $\theta_{D_1}^{p_1, q_0} = \Theta_{A}^{p_0, q_0}$. Therefore, the domain $D_1$ will contribute to the differential in $\CFD{k}$ if and only if $\Theta_A^{p_0, q_0} > \Theta_{A}^{x, y}$.

The domains $D_2$, $D_3$, and $D_4$ can be analyzed in the same way. For $D_3$ the results are the same: cutting along $\beta_1$ from $q_1$ to $\alpha_0$ makes $\theta_{D_3}^{x, y} - \theta_{D_3}^{p_0, q_1}$ approach $+\infty$, and cutting along $\alpha_1$ from $q_1$ to $q_0$ yields $\Theta_{D_3}^{x, y} - \theta_{D_3}^{p_0, q_1} = \Theta_{A}^{x, y} - \Theta_{A}^{p_0, q_0}$, so $D_3$ contributes if and only if $\Theta_A^{p_0, q_0} > \Theta_{A}^{x, y}$. The domains $D_2$ and $D_4$, on the other hand, contribute if and only if $\Theta_{A}^{x, y} > \Theta_A^{p_0, q_0}$.
\end{proof}
\end{prop}

We will often encounter annular domains which fit the form of the annuli in Proposition \ref{Levine_prop} except that one boundary component has more than one $\alpha$ segment and more than one $\beta$ segment. For instance, the bigon $B_1$ might be replaced with a quadrilateral. In practice, quadrilaterals behave like bigons in this context, but it is not immediately apparent how to extend the proof of Proposition \ref{Levine_prop} for more general annuli. Instead, we will use the following proposition to simplify a domain by pinching off an extra $\alpha$ or $\beta$ arc.

\begin{prop}
\label{pinching}
Let $\gamma$ be an arc in a domain $D$ which is a small pushoff of one of the $\beta$ segments or one of the $\alpha$ segments (possibly containing Reeb chords) in $\partial D$, as pictured below. Assume that $\gamma$ only passes through regions with multiplicity 1 in $D$. Let $D'$ be the domain which results from collapsing $\gamma$ to a point and removing the bigon on the left. Then for an appropriate choice of complex structure, $D$ contributes to $\CFD{k}$ if and only if  $D'$ would contribute.

\begin{figure}[htbp]
\begin{center}
\begin{overpic}{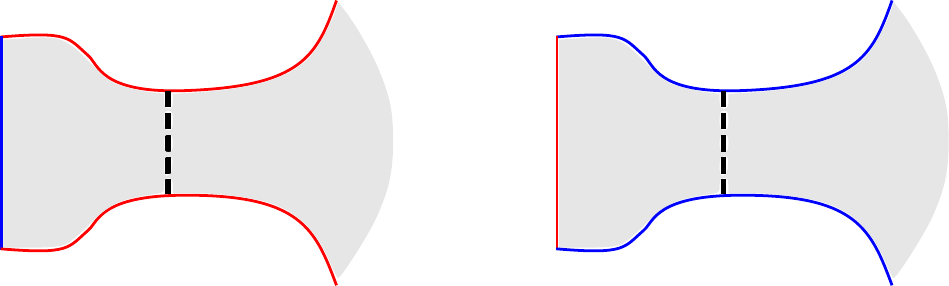}
\put(15, 14){\scriptsize $\gamma$}
\put(73, 14){\scriptsize $\gamma$}

\put(-4, 14){$\color{blue} \beta$}
\put(30, 28){\color{red} $\alpha_1$}
\put(30, 0){\color{red} $\alpha_2$}

\put(54, 14){\color{red} $\alpha$}
\put(88, 28){$\color{blue} \beta_1$}
\put(88, 0){$\color{blue} \beta_2$}

\end{overpic}
\end{center}
\end{figure}

\begin{proof}
Given a complex structure $J$ on $D$, identify a neighborhood of $\gamma$ with $[0, 1] \times (-\epsilon, \epsilon)$. Consider the one parameter family of complex structures $J_t$, $t > 0$, such that the same neighborhood of $\gamma$ is identified with $[0,1] \times (-t,t)$ and $J_t$ agrees with $J$ outside that neighborhood. As $t$ goes to infinity, the neighborhood of $\gamma$ is pinched and stretched---effectively the arc $\gamma$ becomes shorter. The limiting complex structure $J_\infty$ corresponds to $\gamma$ being pinched to a single point, resulting in a bigon $B$ and the the domain $D'$ joined at a point. A $J_\infty$-holomorphic curve with domain $D$ splits as a holomorphic map to the bigon and a holomorphic curve with domain $D'$. By the Riemann Mapping Theorem, there is an $\RR$ family of holomorphic maps from the standard bigon to $B$, and precisely one once the point of contact with the preimage of $D'$ is determined. Therefore the existence of $J_\infty$-holomorphic curves with domain $D$ is equivalent to the existence of $J_\infty$-holomorphic curves with domain $D'$. 

$J_\infty$ is not a valid complex structure to choose when computing $\CFD{k}$, but we can choose $J_t$ for arbitrarily large $t$, and standard compactness and gluing arguments show that for $t$ sufficiently large, $D$ has a $J_t$-holomorphic representative if and only if it has a $J_\infty$-holomorphic representative. Thus for a complex structure with the arc $\gamma$ sufficiently pinched, the statement of the proposition holds.

\end{proof}
\end{prop}

Finally, we discuss how $\Ainfty$ relations can be used to deduce which domains count toward the differential, even if we are computing $\CFD{k}$. The key is the following observation:

\begin{lem}
\label{Ainfty_trick}
A sequence of Reeb chords $\vecrho = (\vecrho[1], \ldots, \vecrho[k])$ contributes $a(\vecrho) \otimes \y$ to the differential of $\x$ in $\CFD{k}$ if and only if $m(\x, \bar{a}( \vecrho[1] ), \ldots, \bar{a}( \vecrho[k] )) = \y$ in $\CFA{k}$.

\begin{proof}
This follows directly from the definitions of $\CFD{k}$ and $\CFA{k}$, since both involve counts of the same moduli spaces. For a given domain $B \in \pi_2(\x, \y)$ that is compatible with $\vecrho$, the pair $(B, \vecrho)$ may contribute $a(\vecrho) \otimes \y$ to $\partial \x$ in $\CFD{k}$, and it may contribute the operation $m(\x, \bar{a}(\vecrho[1]), \ldots, \bar{a}(\vecrho[k])) = \y$ to $\CFA{k}$. In both cases, the pair contributes if and only if $\#(\moduli) \equiv 1$ (mod 2).

\end{proof}
\end{lem}

Here we say that $\vecrho$ \emph{contributes} to $\CFD{k}$ if the relevant counts of moduli spaces are nonzero, even if the contribution $a(\vecrho) \otimes \y$ may be zero. In comparison, notice that Proposition \ref{idempotent_proposition} and Lemma \ref{index_lemma} discuss when a domain \emph{contributes a nonzero term} to $\CFD{k}$. Thus a domain that is ruled out from consideration for $\CFD{k}$ by Proposition \ref{idempotent_proposition} or Lemma \ref{index_lemma} might still contribute to $\CFA{k}$.

Lemma \ref{Ainfty_trick} is most useful for checking if domains contribute to $\CFD{k}$ when $\vecrho$ contains the long chord $\rho^i_{123}$ for some boundary component. For example, suppose in the one boundary case that $\vecrho = (\rho_{123})$ is compatible with a domain $B$ from $\x$ to $\y$. If the domain is too complicated to understand the moduli space $\moduli$ directly, we can instead ask whether $(B, \vecrho)$ contributes the operation $m(\x, \rho_{123}) = \y$ to $\CFA{}$. To answer this, we consider the $\Ainfty$ relation (Equation \ref{Ainfty_relation}) corresponding to $\x$ and $\vec{\rho '} = (\rho_{12}, \rho_3)$. The relation says that
\[ 0 = m\bigg(\x, \mu(\rho_{12}, \rho_3) \bigg) + m\bigg( m(\x, \rho_{12}) , \rho_3 \bigg) . \]
The first term is just $m(\x, \rho_{123})$, the operation we are interested in, and the second term might be easier to analyze. If $m(\x, \rho_{12}) = 0$, for instance, then the second term in the relation is 0, and thus $(B, \vecrho)$ does not contribute to $\CFA{}$ or to $\CFD{}$.

\subsection{Gradings}
\label{sec:gradings}

Bordered Heegaard Floer invariants can be equipped with a relative grading on each spin$^c$-structure as described in \cite[Chapter 10]{LOT:Bordered} and \cite[Section 6.5]{LOT:Bimodules}. We recall here the construction of these gradings for manifolds with only torus boundary components. We will only discuss the refined grading.

Let $Y$ be a bordered manifold represented by a bordered Heegaard diagram $\HeegDiag$. Let $\CFhat(\HeegDiag)$ denote the relevant bordered Heegaard Floer invariant. The gradings for $\CFhat(\HeegDiag)$ lie in a noncommutative group which depends on the number and type of boundary components. We will denote this group $G_{n, m}$ where $n$ is the number of type $D$ boundary components of $Y$ and $m$ is the number of type $A$ boundary components. $G_{n,m}$ is generated by tuples $(j; a_1, b_1; a_2, b_2; \ldots; a_{n+m}, b_{n+m})$, where every entry is in $\frac{1}{2}\Z$, and $a_i + b_i \in \Z$ for each $i$. $j$ is referred to as the \emph{Maslov component} of the grading. Multiplication on this group is defined as follows:

$$(j; a_1, b_1; a_2, b_2; \ldots; a_{n+m}, b_{n+m}) \cdot (j'; a'_1, b'_1; a'_2, b'_2; \ldots; a'_{n+m}, b'_{n+m}) = $$
$$( j + j' +  C; a_1+a'_1, b_1+b'_1; a_2+a'_2, b_2+b'_2; \ldots; a_{n+m}+a'_{n+m}, b_{n+m}+b'_{n+m} ),$$
where the correction term $C$ is given by 
$$C = \begin{array}{|cc|} a_1 & a'_1 \\ b_1 & b'_1\end{array} + \ldots + \begin{array}{|cc|} a_n & a'_n \\ b_n & b'_n\end{array} +  \begin{array}{|cc|} a'_{n+1} & a_{n+1} \\ b'_{n+1} & b_{n+1}\end{array} + \ldots +  \begin{array}{|cc|} a'_{n+m} & a_{n+m} \\ b'_{n+m} & b_{n+m}\end{array} \,.$$

Given generators $\x$ and $\y$, a domain in $B \in \pi_2(\x, \y)$ can be given a grading in $G_{n,m}$ \cite[Definition 10.1]{LOT:Bordered}. The Maslov component of $gr(B)$ is given by
\[ -e(B)-n_\x(B)-n_\y(B), \]
where $e(B)$, $n_\x(B)$, and $n_\y(B)$ are are the same quantities that appear in the index formula, Equation \ref{index1}. For each $1 \le i \le n+m$, let $\gamma_i$ denote the intersection of $\partial B$ with the $i$th boundary component of $\HeegDiag$, which can be thought of as a linear combination of the Reeb chords $\rho^i_1$, $\rho^i_2$, and $\rho^i_3$. If $\gamma_i = c_1 \rho^i_1 + c_2 \rho^i_2 + c_3 \rho^i_3$, then the $i$th pair of coefficients in $gr(B)$ is given by
\[ a_i = \frac{ c_1 + c_2 - c_3}{2}, \qquad  b_i = \frac{ -c_1 + c_2 + c_3}{2}. \]

To define the gradings on a bordered multimodule in a given spin$^c$-structure, we choose a base generator $\x$ in that spin$^c$-structure. Let $\mathcal{P}(\x)$ be the subgroup of $G_{n,m}$ generated by $\{ gr(B) | B\in \pi_2(\x, \x) \}$. $\CFhat(\HeegDiag)$ then has a well defined grading by the set $G_{n,m} / \mathcal{P}(\x)$. Up to canonical isomorphism, this grading set does not depend on the choice of $\x$. We define the relative grading by the following rule: if $\y$ is generator in the same spin$^c$-structure as $\x$ and $B$ is a domain connecting $\x$ to $\y$, then $ gr(\y) = gr(\x) gr(B)$.

In many cases, gradings can be computed directly from the labeled graph representing $\CFhat(\HeegDiag)$, without reference to the Heegaard diagram. To do this, we use the fact that elements of the torus algebra have gradings in $G_{n,m}$. Recall that $\CFhat(\HeegDiag)$ is a module over $n+m$ copies of the torus algebra, one for each boundary of $Y$, and $\rho^i_I$ denotes an element of the torus algebra associated to the $i$th boundary. The Maslov component of $gr(\rho^i_I)$ is $-\frac{1}{2}$ and the coefficients $a_j$ and $b_j$ are zero for all $j \ne i$. The coefficients $a_i$ and $b_i$ are determined by $I$ as follows:
\begin{center}
\begin{tabular}{|c|c|c |}
\hline
$I$ & $a_i$ & $b_i$ \\
\hline
1 & $\frac{1}{2}$ & $-\frac{1}{2}$ \\
2 & $\frac{1}{2}$ & $\frac{1}{2}$ \\
3 & $-\frac{1}{2}$ & $\frac{1}{2}$ \\
12 & $1$ & $0$ \\
23 & $0$ & $1$ \\
123 & $\frac{1}{2}$ & $\frac{1}{2}$ \\
\hline
\end{tabular}
\end{center}
This grading respects the algebra product in the sense that $gr(\rho_{I_1} \rho_{I_2}) = gr(\rho_{I_1}) gr(\rho_{I_2})$. The grading on $\CFhat(\HeegDiag)$ also respects the module structure in the sense that $gr(\rho_I x) = gr(\rho_I) gr(x)$, where the product on the right refers to the left action of the group $G_{n,m}$ on the set $G_{n,m}/ \mathcal{P}(\x)$. Finally, the grading on $\CFhat(\HeegDiag)$ satisfies the following relation  \cite[Definition 2.5.2]{LOT:Bimodules}:
\begin{equation}
\label{grading_relation}
gr( m_{k+1}(x, \rho_{I_1}, \ldots, \rho_{I_k}) ) = \lambda^{k-1} gr(\rho_{I_k}) \cdots gr(\rho_{I_1}) gr(x) .
\end{equation}
Here $\lambda = (1; 0,0; \ldots; 0, 0)$ is the preferred central element of $G_{n,m}$. The same relation applies for both type $D$ and type $A$ modules if we think of the differential $\partial$ as an $m_1$ map. Thus $gr(\partial x) = \lambda^{-1} gr(x)$.

To compute the relative grading from the graph representing $\CFhat(\HeegDiag)$, we choose a base generator $\x$ and assign it an arbitrary grading. The gradings of the remaining generators can be determined using Equation \ref{grading_relation}, as long as each generator is connected by $\x$ by a path of arrows (that is, as long as the graph is connected). A loop in the graph representing $\CFhat(\HeegDiag)$, along with Equation \ref{grading_relation}, gives rise to a value for $gr(\x)$ which may not be equal to the value initially chosen for $gr(\x)$. The difference is $gr(B)$ for some periodic domain $B \in \pi_2(\x, \x)$. If there are enough loops in the graph (there must be one independent periodic domain for each boundary component of $Y$), then we can determine $\mathcal{P}(\x)$.

\section{Direct Computation of $\CFD{k}(\Y)$}
\label{sec:trimodule}

\begin{figure}
\label{deriving_Heeg_diag}
\caption{Constructing the Heegaard diagram for $\Y$.}
\centering

\bigbreak
\includegraphics[scale = .7]{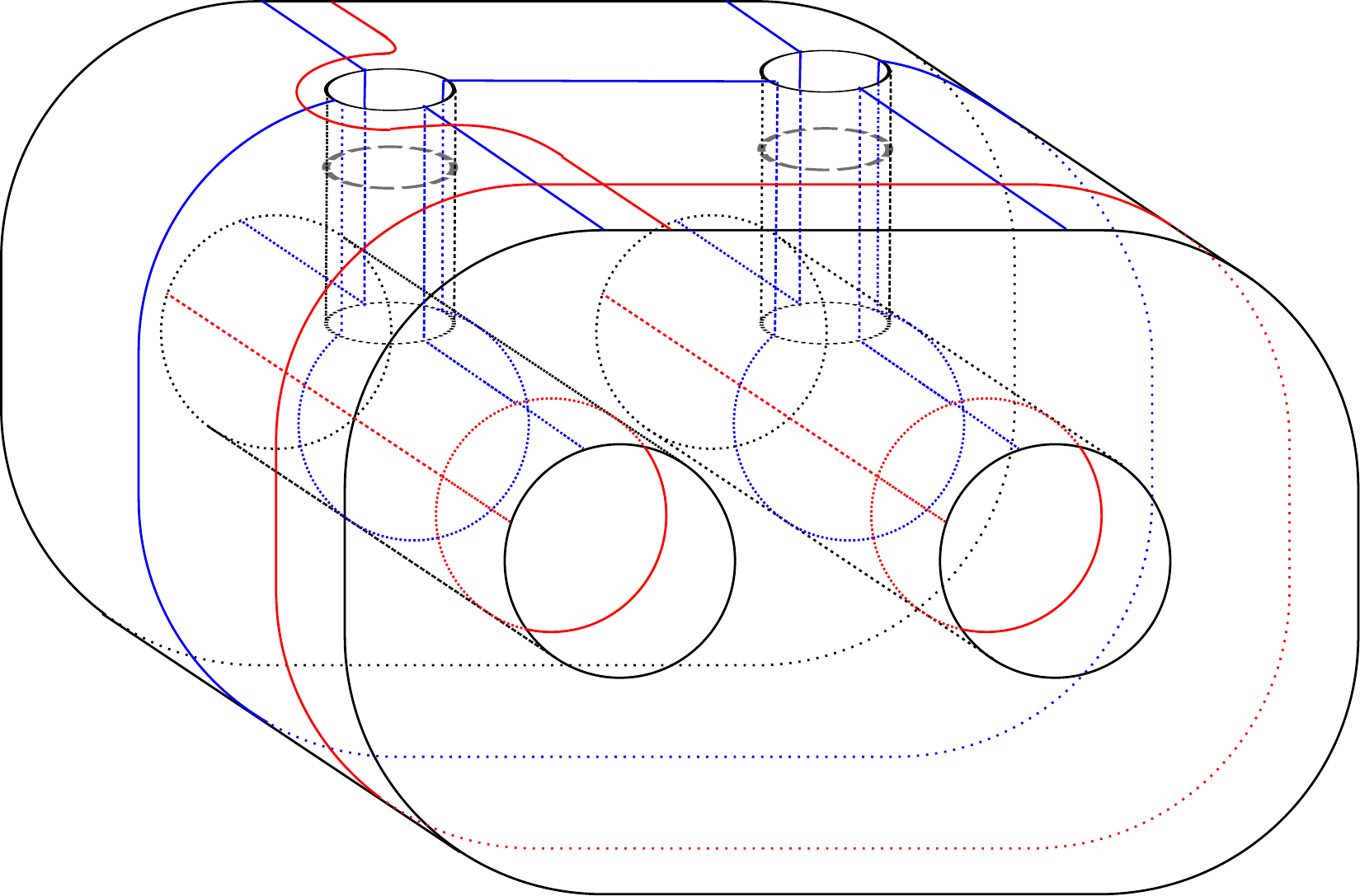}

(a) The front and back faces are identified. Attaching 2-handles to the $\beta$ curves (blue), capping off the drilling tubes along the dotted gray lines, and adding a 3-handle yields $\Y$. Intersecting pairs of $\alpha$ curves (red) specify a parametrization of each boundary component of $\Y$. 
\bigbreak

\includegraphics[scale = .6]{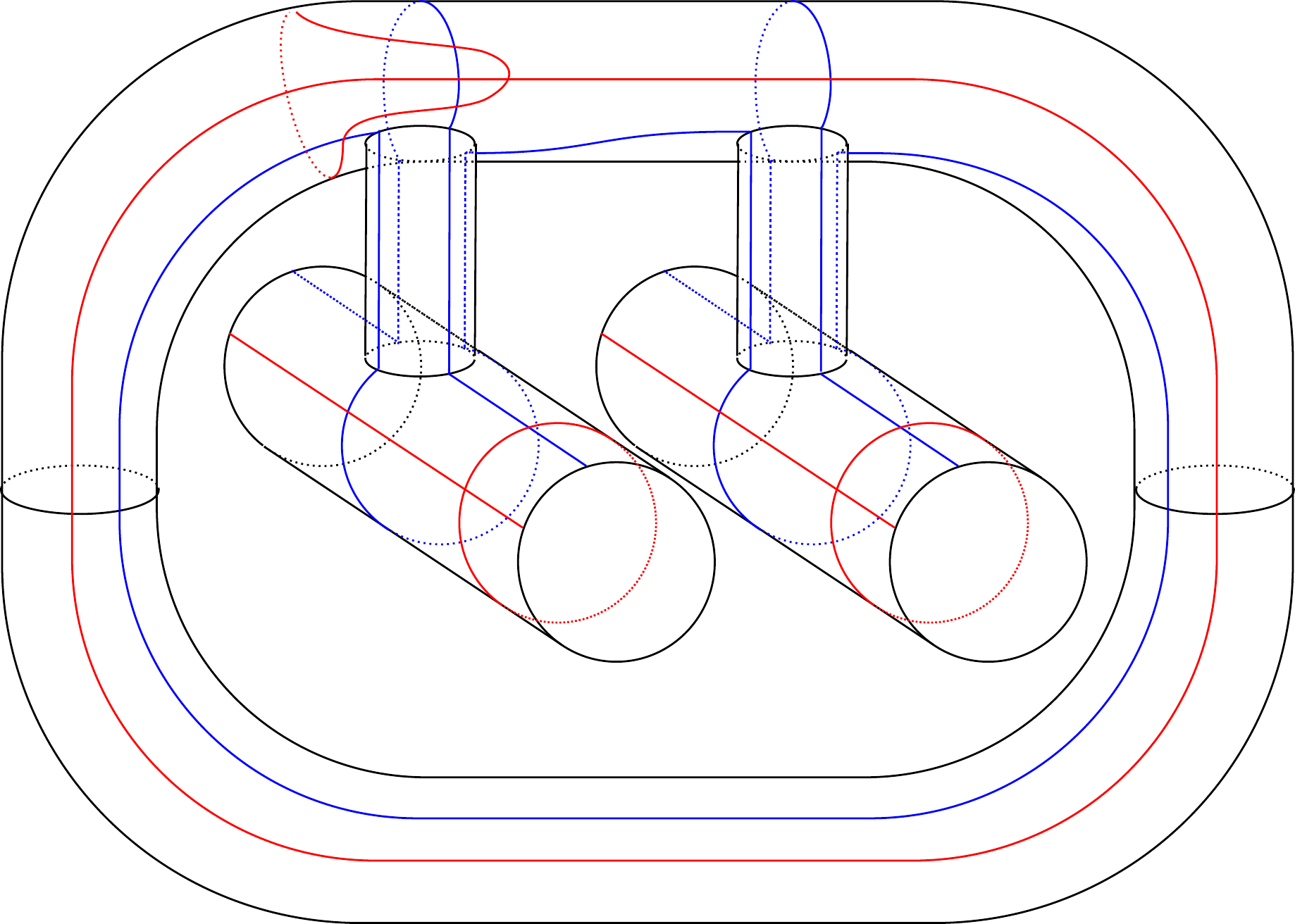}

(b) The diagram is easier to read and manipulate if we redraw the outer torus.
\end{figure}

\begin{figure}
\centering

\begin{overpic}[scale = 1]{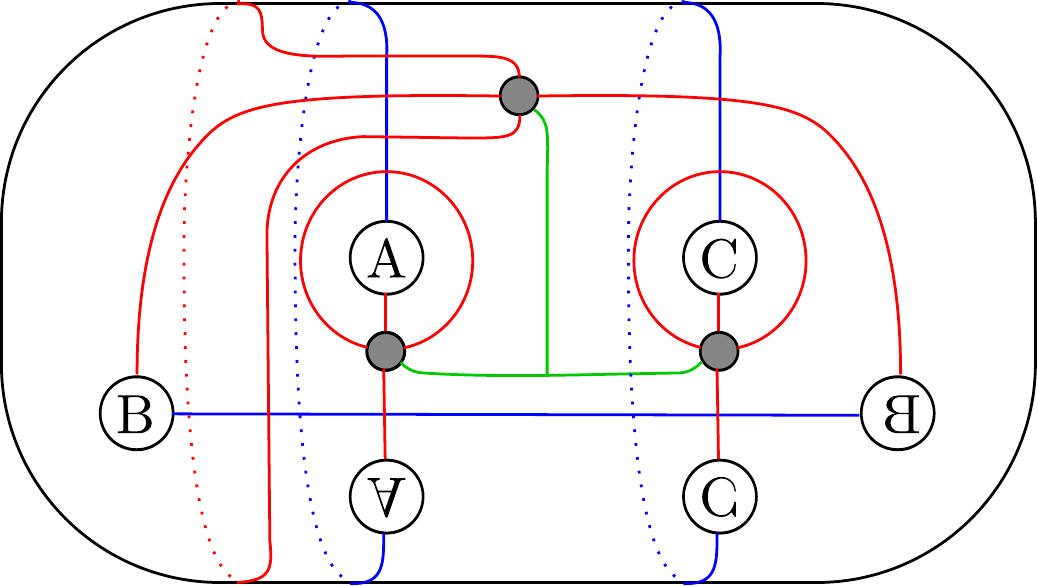}
\end{overpic}

(c) An isotopy simplifies the diagram. Pairs of circles labeled by letters signify 1-handle attachment. We remove a small disk (shaded) around each intersection of $\alpha$ curves, resulting in the genus 3 Heegaard surface $\Sigma$ with three boundary components. There is a basepoint $z$ connected by arcs (green) to each boundary component.
\end{figure}

\begin{figure}
\centering
\begin{overpic}[scale = 1.1]{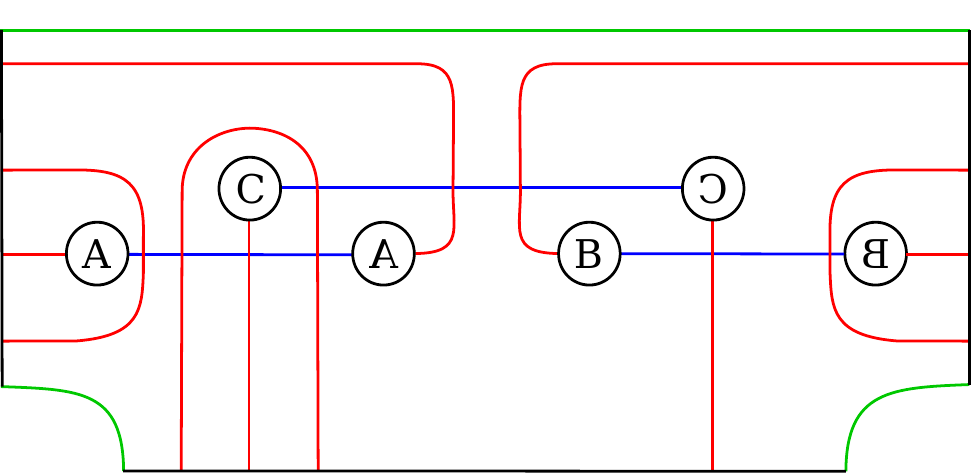}

\put(-4, 17){$\rho_1$}
\put(-4, 26){$\rho_2$}
\put(-4, 36){$\rho_3$}

\put(100, 36){$\sigma_1$}
\put(100, 26){$\sigma_2$}
\put(100, 17){$\sigma_3$}

\put(21, -3){$\tau_3$}
\put(28, -3){$\tau_2$}
\put(50, -3){$\tau_1$}

\put(1, 15){\color{red} \scriptsize $\alpha^\rho_1$}
\put(1, 24){\color{red} \scriptsize $\alpha^\rho_2$}
\put(95, 15){\color{red} \scriptsize $\alpha^\sigma_2$}
\put(95, 24){\color{red} \scriptsize $\alpha^\sigma_1$}
\put(33, 5){\color{red} \scriptsize $\alpha^\tau_2$}
\put(74, 5){\color{red} \scriptsize $\alpha^\tau_1$}

\put(33, 20){\color{blue} \tiny $\beta_1$}
\put(63, 30){\color{blue} \tiny $\beta_2$}
\put(78, 20){\color{blue} \tiny $\beta_3$}

\end{overpic}
\label{HeegDiag_simple}
\caption{A bordered Heegaard diagram $\HeegDiag$ for $\Y$, with type $D$ boundaries.}
\end{figure}

\begin{figure}
\centering
\begin{overpic}[scale = 1.1, tics = 5]{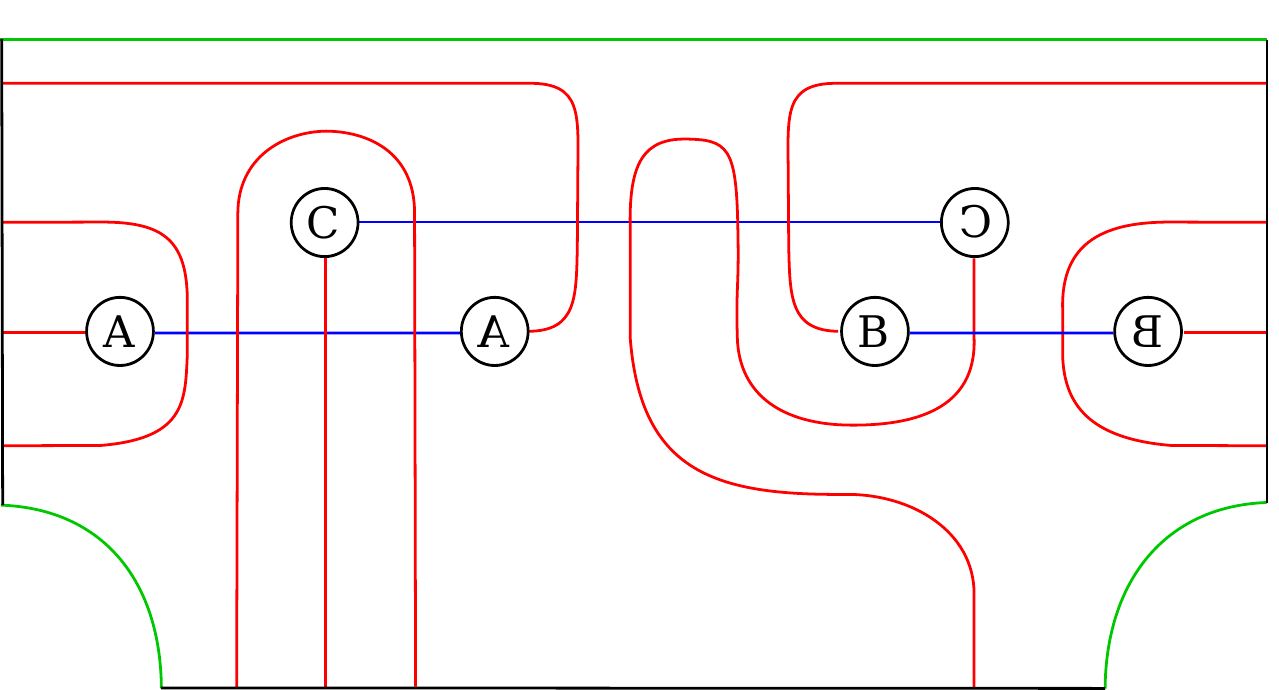}
\put(-3, 23){$\rho_1$}
\put(-3, 32){$\rho_2$}
\put(-3, 41){$\rho_3$}
\put(100, 41){$\sigma_1$}
\put(100, 32){$\sigma_2$}
\put(100, 23){$\sigma_3$}
\put(21, -2){$\tau_3$}
\put(28, -2){$\tau_2$}
\put(50, -2){$\tau_1$}

\put(13, 41){$R_1$}
\put(4, 32){$R_2$}
\put(4, 22){$R_3$}
\put(20, 15){$R_4$}
\put(27, 15){$R_5$}
\put(41, 15){$R_3$}
\put(20, 31){$R_7$}
\put(27, 31){$R_6$}
\put(37, 33){$R_2$}
\put(67, 33){$R_6$}
\put(83, 41){$R_7$}
\put(92, 32){$R_6$}
\put(92, 22){$R_8$}
\put(61, 24){$R_8$}
\put(52, 39){$R_9$}

\put(15 , 26){\scriptsize $a$}
\put(19 , 26){\scriptsize $b$}
\put(26 , 26){\scriptsize $c$}
\put(33 , 26){\scriptsize $d$}
\put(33 , 37.5){\scriptsize $e$}
\put(45.5 , 37.5){\scriptsize $f$}
\put(50 , 37.5){\scriptsize $j$}
\put(58.5 , 37.5){\scriptsize $k$}
\put(62.5, 37.5){\scriptsize $g$}
\put(77 , 26){\scriptsize $h$}
\put(84 , 26){\scriptsize $i$}

\put(11, 18){\color{red} \scriptsize $\alpha^\rho_1$}
\put(11, 48.5){\color{red} \scriptsize $\alpha^\rho_2$}
\put(85, 18){\color{red} \scriptsize $\alpha^\sigma_2$}
\put(85, 48.5){\color{red} \scriptsize $\alpha^\sigma_1$}
\put(15, 10){\color{red} \scriptsize $\alpha^\tau_2$}
\put(77, 10){\color{red} \scriptsize $\alpha^\tau_1$}

\put(34, 29){\color{blue} \tiny $\beta_1$}
\put(71, 37.5){\color{blue} \tiny $\beta_2$}
\put(85, 29){\color{blue} \tiny $\beta_3$}

\end{overpic}
\caption{A slightly modified arced Heegaard diagram $\Y$, $\HeegDiag'$.  }
\label{second}
\end{figure}

In this section we explicitly compute the type $D$ trimodule associated to $\Y$, the trivial $S^1$-bundle over the pair of paints $\pants$.

\subsection{Choosing a bordered Heegaard diagram for $\Y$}
We obtain the Heegaard surface $\Sigma$ from the boundary of $\Y$ by drilling through $\Y$ to connect each inner torus boundary component with the outer boundary component. This surface is pictured in Figure \ref{deriving_Heeg_diag}(a), where the front and back faces are identified by the identity map. To obtain the 3-manifold $\Y$ from this surface, we attach three 2-handles to the inside along the $\beta$ curves, and fill in the drilling tubes by attaching disks along the thick grey dotted lines. Filling in the interior with a 3-ball yields $\Y$.

We decorate each boundary component of $\Y$ with a pair of intersecting $\alpha$ curves to specify a parametrization of the boundary. There are many choices of parametrization, but it is natural and convenient to choose one $\alpha$ curve to lie in the base surface $\pants$ and the other to be an $S^1$ fiber.

To finish the construction of $\HeegDiag$ we must also remove a disk around each $\alpha-\alpha$ intersection point, turning the $\alpha$ curves into arcs and giving the surface $\Sigma$ three boundary components, and we must chose a basepoint $z$ connected by arcs to each component of $\partial\Sigma$. These features are not shown until after the surface has been simplified by isotopy in Figure \ref{deriving_Heeg_diag}(b) and (c). Notice the placement of the $\alpha$ curves relative to the drilling tunnels in Figure \ref{deriving_Heeg_diag}(a). This was to ensure that there is one component of $\Sigma \backslash (\boldalpha \cup \boldbeta)$ that meets all three components of $\partial\Sigma$.

The borderd Heegaard diagram $\HeegDiag$ in \ref{deriving_Heeg_diag}(c) can be represented in the plane (except for the handles) by cutting $\Sigma$ open along the arcs in $\z$. The result is shown in Figure \ref{HeegDiag_simple}, and some relevant labels have been added. The Reeb chords along the three boundary components are labeled in the order consistent with type $D$ boundaries, and they are denoted by $\rho$, $\sigma$, and $\tau$ instead of $\rho^1$, $\rho^2$, and $\rho^3$. The $\alpha$ arcs are also labeled to correspond to type $D$ boundaries. Tracing through the sequence in Figure \ref{deriving_Heeg_diag} with $\alpha$ labels in mind, note that the arcs $\alpha^\rho_1$, $\alpha^\sigma_2$, and $\alpha^\tau_1$ represent curves in the base surface $\pants$ of the $S^1$-bundle $\Y$, and the arcs $\alpha^\rho_2$, $\alpha^\sigma_1$, and $\alpha^\tau_2$ represent $S^1$ fibers.

Before computing $\CFD{3} (\Y)$, we make one final adjustment to the Heegaard diagram $\HeegDiag$. Computing directly from $\HeegDiag$ would involve a few large and complicated domains which are difficult to analyze. It will be convenient to perform an isotopy to produce the new Heegaard diagram $\HeegDiag'$ shown in Figure \ref{second}. This change introduces a few extra generators, but it eliminates the trickiest domains and overall makes the computation easier.

\subsection{Complex Structure}

To compute $\CFD{3}$, we also must fix a generically chosen complex structure $J$ on $\Sigma$. We collect here some relevant choices about $J$ that will be used in the computation. Note that if $J$ were chosen with different properties we would get a different, but quasi-isomorphic, trimodule.

As in the proof of Proposition \ref{Levine_prop}, if the annulus $A$ has one $\alpha$ and one $\beta$ segment on each boundary component, let $\Theta_A^{p,q}$ denote the ratio of the lengths of the $\alpha$ and $\beta$ segments on the boundary component which contains $p$ and $q$. As shorthand we will use, for example, the subscript $``67"$ for the annulus $R_6 R_7$. We will assume that $J$ satisfies:
\begin{itemize}
\item
$\Theta_{67}^{d, b} > \Theta_{67}^{i, i}$
\item
$\Theta_{1267}^{a, a} > \Theta_{1267}^{i, i}$
\end{itemize}
We will also assume that the following arcs are arbitrarily short, as in Proposition \ref{pinching}:
\begin{itemize}
\item
an arc in $R_3$ from $\alpha^\rho_2$ to $\alpha^\tau_1$;
\item
an arc in $R_8$ from $\alpha^\sigma_1$ to $\alpha^\tau_1$;
\item
an arc in $R_1 \cup R_2$ parallel to $\alpha^\tau_2$, from $\beta_1$ to itself.
\end{itemize}

It is straightforward to check that these choices are consistent with each other and that there are suitable complex structures satisfying all of these conditions.

\subsection{Generators}

The chain complex $\CFD{3}(\HeegDiag')$ is generated by the set $\Gen(\HeegDiag')$ consisting of triples of intersection points with one point on each $\beta$ circle and at most one point on each $\alpha$ arc. In total, there are 23 generators. These generators fall into seven different spin$^c$-structures, corresponding to how many $\alpha$ arcs are occupied on each boundary component. 

We begin by computing the summand of $\CFD{3}(\HeegDiag')$ corresponding to the middle spin$^c$-structure, with exactly one $\alpha$ arc occupied on each boundary (the other spin$^c$-structures are much easier and will be addressed at the end of this section). There are seven generators in the middle spin$^c$-structure:  $x = (a,e,i)$,  $y = (a, g, h)$,  $z = (c, f, i)$,  $v = (b, f, i)$,  $w = (d, f, i)$,  $s = (a, j, i)$,  and $t = (a, k, i)$.

\subsection{Possible Domains}

We begin by listing domains in $\pi_2(\x, \y)$ for any pair of generators $\x$ and $\y$. Recall that a domain connecting $\x$ and $\y$ is a linear combination of the regions $R_1, \ldots, R_9$ in Figure \ref{second} with the correct multiplicity at each corner (that is, satisfying Equation \ref{corner_condition}). We do not need to list every domain in $\pi_2(\x, \y)$, since only positive domains can have holomorphic representatives, and by Proposition \ref{idempotent_proposition} we can assume that regions which intersect $\partial\Sigma$ have multiplicity $0$ or $1$. The only region which does not intersect the boundary is $R_9$. The multiplicity of $R_9$ is also limited; in fact, $R_9$ can not combine with any other regions because domains with holomorphic representatives must be connected in $\Sigma \backslash (\alpha \cap \beta)$. Thus we can restrict to linear combinations of $R_1, \ldots, R_8$ with multiplicity 0 or 1 for each region, and the single domain $R_9$. All such domains which connect two generators are listed in Table \ref{table1}.

\begin{table}
\begin{center}
\caption{Domains which potentially contribute to the differential. All subsets of the regions $\{ R_1, \ldots, R_8 \}$ with the proper corner multiplicities, and the single domain $R_9$. We omit the $R$ for the purposes of this table. Thus 56 refers to the domain $R_5 R_6$, which connects $y$ to $x$.}
\label{table1}
\begin{tabular}{| p{30 px} | p{30 px} | p{25 px} | p{30 px} | p{30 px} | p{40 px} | p{30 px} | p{30 px} |}

\hline
\hspace{2px} from \bigbreak to& $x$ & $y$ & $z$ & $w$ & $v$ & $s$ & $t$ \\  \hline

$x$ & 
   $1 2 6 7$   $4 5 6 7$  &
$5 6$  &
$1 5 6 7$ &
$1 6 7$ &
$1$  \quad  $1 4 5 6 7$ &
 - &
$5 6 8$ \\  \hline
 
$y$ & 
$4 7$ &
   $1 2 6 7$   $4 5 6 7$  &
$1 7$ &
 -  &
$1 4 7$ &
 - &
$8$   $1 2 6 7 8$   $4 5 6 7 8 $ \\  \hline
 
$z$ & 
$2 4 6 7$ &
$2 6$ &
   $1 2 6 7$   $4 5 6 7$  &
$4 6 7$ &
$4$ \hspace{40 px}  $1 2 4 6 7$ &
 - &
$2 6 8$ \\  \hline

$w$ & 
$2$   $2 4 5 6 7$ &
$2 5 6$ &
$5$   $1 2 5 6 7$  &
   $1 2 6 7$   $4 5 6 7$  &
$1 2$ \hspace{40 px}  $4 5$   $1 2 4 5 6 7$ &
 - &
$2 5 6 8$ \\  \hline

$v$ & 
$2 6 7$ &
 - &
$5 6 7$ &
$6 7$ &
   $1 2 6 7$   $4 5 6 7$  &
 - &
- \\  \hline

$s$ & 
$2 3$   $2 3 4 5 6 7$ &
$2 3 5 6$ &
$ 3 5$   $1 2 3 5 6 7$  &
$3$   $1 2 3 6 7$   $3 4 5 6 7$ &
$1 2 3$ \hspace{30 px}  $3 4 5$   $1 2 3 4 5 6 7$ &
   $1 2 6 7$   $4 5 6 7$  &
$2 3 5 6 8$  $9$ \\  \hline

$t$ & 
 -  &
 -  &
 -  &
 -  &
 -  &
 -  &
   $1 2 6 7$   $4 5 6 7$    \\  \hline

\end{tabular}
\end{center}
\end{table}

\subsection{Compatibility and Idempotents} Several domains in Table \ref{table1} can be ruled out using Proposition \ref{idempotent_proposition}. Consider for example the domain $R_2 R_3$, which potentially contributes to the differential from $x$ to $s$. By Proposition \ref{idempotent_proposition}, this domain can only contribute with the Reeb chords $(\rho_1, \rho_2)$, and then the contribution $a(\rho_{12}) \otimes s$ is zero unless $s$ contains a point on $\alpha^\rho_2$. Since $s$ does not contain a point on  $\alpha^\rho_2$, the domain $R_2 R_3$ has zero contribution to the differential.

In general, for a differential ending in $s$ to be nontrivial, the algebra element for the $\rho$ boundary can not be $\rho_2$ or $\rho_{12}$. This means that the domain associated with such a differential can not contain $R_2$ without containing $R_1$. In addition to $R_2 R_3$, this line of reasoning eliminates the domains $R_2 R_3 R_4 R_5 R_6 R_7$, $R_2 R_3 R_5 R_6$, and $R_2 R_3 R_5 R_6 R_8$. A similar analysis on the other boundaries shows that domains contributing nontrivial differentials ending in $s$ cannot contain $R_5$ without $R_4$, $R_7$ without $R_6$, or $R_8$. This further rules out the domains $R_3 R_5$ and $R_1 R_2 R_3 R_5 R_6 R_7$. Finally, applying  the same technique to differentials ending in other generators rules out the following domains:

\begin{center}
\begin{tabular}{ r  l }
\text{domains to} x : & 4567, 14567, 568 \\
\text{domains to} y : & 1267, 4567\\
\text{domains to} z : & 1267, 12467, 268 \\
\text{domains to} w : & 24567, 12567, 1267, 4567, 12, 45,  124567, 2568 \\
\text{domains to} v : & 1267, 4567\\
\end{tabular}
\end{center}

\subsection{Polygons}

Of the remaining domains from Table \ref{table1}, many are immersed polygons and therefore contribute to the differential by Proposition \ref{polygons}. The proposition depends on the sequence of Reeb chords $\vecrho$, but each of the following domains has only one compatible sequence of Reeb chords, so Proposition \ref{polygons} tells us the entire contribution of the domain to the differential:

\begin{center}
\begin{tabular}{c >{contributes}c c p{50px} c >{contributes}c c}
$R_1$ & & $v \overset{\rho_3}\longrightarrow x$ & &$R_2$ & & $x \overset{\rho_2}\longrightarrow w$ \\
$R_3$ & & $w \overset{\rho_1 \tau_1}\longrightarrow s$ & & $R_4$ & & $v \overset{\tau_3}\longrightarrow z$ \\
$R_5$ & & $z \overset{\tau_2}\longrightarrow w$ & & $R_8$ & & $t \overset{\sigma_3}\longrightarrow y$ \\
$R_9$ & & $t \overset{1}\longrightarrow s$ & & $R_1 R_7$ & & $z \overset{\rho_3 \sigma_1}\longrightarrow y$ \\
$R_2 R_6$ & & $y \overset{\rho_2 \sigma_2}\longrightarrow z$ & & $R_4 R_7$ & & $x \overset{\tau_3 \sigma_1}\longrightarrow y$ \\
$R_5 R_6$ & & $y \overset{\tau_2 \sigma_2}\longrightarrow x$ & & & \multicolumn{1}{c}{} &  
\end{tabular}
\end{center}

Here the notation $\x \overset{a}\longrightarrow \y$ means that there is an $a \otimes \y$ term in $\partial \x$.

\subsection{Index}

At this point there are 22 domains in Table \ref{table1} whose contribution to $\CFD{k}$ remain unknown. Of these, 11 can be ruled out by showing that $\ind(B) \ne 1$. In general computing the index is a good task for a computer, but because we have narrowed the list of domains down so much we can work out the index computations by hand.

The quantities $e(B)$ and $n_\x(B)$ for any generator $\x$ are additive, so it is helpful to record their values for individual regions (see Table \ref{index_by_region}). For instance, region $R_1$ has euler measure $e(R_1) = -\tfrac{1}{2}$, because it has euler characteristic $\chi(R_1) = 1$ and six acute corners. $R_1$ has a corner at the point $a$, which means that the average multiplicity of $R_1$ near $a$ is $\tfrac{1}{4}$. $R_1$ also has corners at $b$, $e$, and $f$. For the generator $x = (a, e, i)$ we find that $n_x(R_1) = \tfrac{1}{4} + \tfrac{1}{4} + 0 = \tfrac{1}{2}$, and for $y = (a, g, h)$ we get $n_y(R_1) = \tfrac{1}{4} + 0 + 0 = \tfrac{1}{4}$. It is straightforward to fill in the rest of Table \ref{index_by_region}.

From this information, it is easy to compute the index as in Table \ref{index_table}. We add the euler measures of all of the regions in a given domain $B$ to find $e(B)$. Similarly, we add the values of $n_\x$ and $n_\y$ for each region for the relevant generators $\x$ and $\y$ to find $n_\x(B)$ and $n_\y(B)$. Finally, we count how many of the three components of $\partial \Sigma$ are touched by $B$ (that is, we find $\# \{ Z \in \pi_0(\partial\Sigma) | Z \cap B \ne \emptyset\} $), and add half of this number to $e(B) + n_\x(B) + n_\y(B)$. By Equation \ref{index2}, the result is $\ind(B)$. Table \ref{index_table} only shows the computation for regions that are ruled out by this method. The index can be computed in the same way for the remaining 11 domains, but they all have $\ind(B) = 1$, so more work is needed to determine if they contribute to $\CFD{k}$.

\begin{table}
\caption{}
\label{index_by_region}
\begin{center}
\begin{tabular}{| c | rrrrrrrr |}
\hline
               & \quad $e$\quad & \quad $n_x$ & \quad $n_y$\quad & \quad $n_z$\quad & \quad $n_v$\quad & \quad $n_w$\quad & \quad $n_s$\quad & \quad $n_t$\quad \\ [1mm] \hline
$R_1$  & $-\frac{1}{2}$	& $\frac{1}{2}$ & $\frac{1}{4}$	& $\frac{1}{4}$	& $\frac{1}{2}$ & $\frac{1}{4}$	& $\frac{1}{4}$	& $\frac{1}{4}$ \\  [2mm]
$R_2$  & $-\frac{1}{2}$	& $\frac{1}{2}$ & $\frac{1}{4}$ 	& $\frac{1}{4}$	& $\frac{1}{4}$ & $\frac{1}{2}$	& $\frac{1}{4}$	& $\frac{1}{4}$	\\  [2mm]
$R_3$  & $-1$			& $\frac{1}{4}$ & $\frac{1}{4}$	& $\frac{1}{4}$	& $\frac{1}{4}$	& $\frac{1}{2}$	& $\frac{1}{2}$	& $\frac{1}{4}$ \\ [2mm]
$R_4$  & $0$			& $ 0$ 		& $ 0$ 		& $\frac{1}{4}$	& $\frac{1}{4}$	& $ 0$		& $ 0$		& $ 0$		\\  [2mm]
$R_5$  & $0$			& $ 0$ 		& $ 0$ 		& $\frac{1}{4}$	& $ 0$		& $\frac{1}{4}$	& $ 0$		& $ 0$ 		\\  [2mm]
$R_6$  & $-1$			& $\frac{1}{2}$ 	& $\frac{1}{2}$ 	& $\frac{1}{2}$	& $\frac{1}{4}$	& $\frac{1}{2}$	& $\frac{1}{4}$	& $\frac{1}{4}$	\\  [2mm]
$R_7$  & $-1$			& $\frac{1}{2}$ 	& $\frac{1}{2}$ 	& $\frac{1}{2}$	& $\frac{1}{2}$	& $\frac{1}{4}$	& $\frac{1}{4}$	& $\frac{1}{4}$ \\  [2mm]
$R_8$  & $-\frac{1}{2}$	& $\frac{1}{4}$ 	& $\frac{1}{2}$ 	& $\frac{1}{4}$	& $\frac{1}{4}$	& $\frac{1}{4}$	& $\frac{1}{4}$	& $\frac{1}{2}$ \\ [2mm] \hline
\end{tabular}
\end{center}
\end{table}

\begin{table}
\caption{}
\label{index_table}
\begin{center}
\begin{tabular}{| cc | cccc | c |}
\hline
Domain	& $\x \to \y$ 	& $e(B)$ 		& $n_\x$ 		& $n_\y$		 & bdys hit/2 & $\ind(B)$ \\ [1mm] \hline
1267 	& $x \to x$   	& $-3$ 		& 2 			& 2 			& 1		 	& 2 \\ [2mm]
2467		& $x \to z$   	& $-\frac{5}{2}$	& $\frac{3}{2}$	& $\frac{3}{2}$	& $\frac{3}{2}$	& 2 \\ [2mm]
256		& $y \to w$   	& $-\frac{3}{2}$	& $\frac{3}{4}$	& $\frac{5}{4}$	& $\frac{3}{2}$	& 2 \\ [2mm]
1567		& $z \to x$   	& $-\frac{5}{2}$	& $\frac{3}{2}$	& $\frac{3}{2}$	& $\frac{3}{2}$	& 2 \\ [2mm]
4567		& $z \to z$   	& $-2$		& $\frac{3}{2}$	& $\frac{3}{2}$	& 1			& 2 \\ [2mm]
147		& $v \to y$   	& $-\frac{3}{2}$	& $\frac{5}{4}$	& $\frac{3}{4}$	& $\frac{3}{2}$	& 2 \\ [2mm]
67		& $w \to v$   	& $-2$		& $\frac{3}{4}$	& $\frac{3}{4}$	& $\frac{1}{2}$	& 0 \\ [2mm]
1267		& $s \to s$   	& $-3$		& $1$		& $1$		& $1$		& 0 \\ [2mm]
4567		& $s \to s$   	& $-2$		& $\frac{1}{2}$	& $\frac{1}{2}$	& $1$		& 0 \\ [2mm]
1267		& $t \to t$   	& $-3$		& $1$		& $1$		& $1$		& 0 \\ [2mm]
4567		& $t \to t$   	& $-2$		& $\frac{1}{2}$	& $\frac{1}{2}$	& $1$		& 0 \\ [2mm]
\hline
\end{tabular}
\end{center}
\end{table}

\subsection{Index Zero Annulus}

\begin{figure}
\begin{center}
\includegraphics[scale = .8]{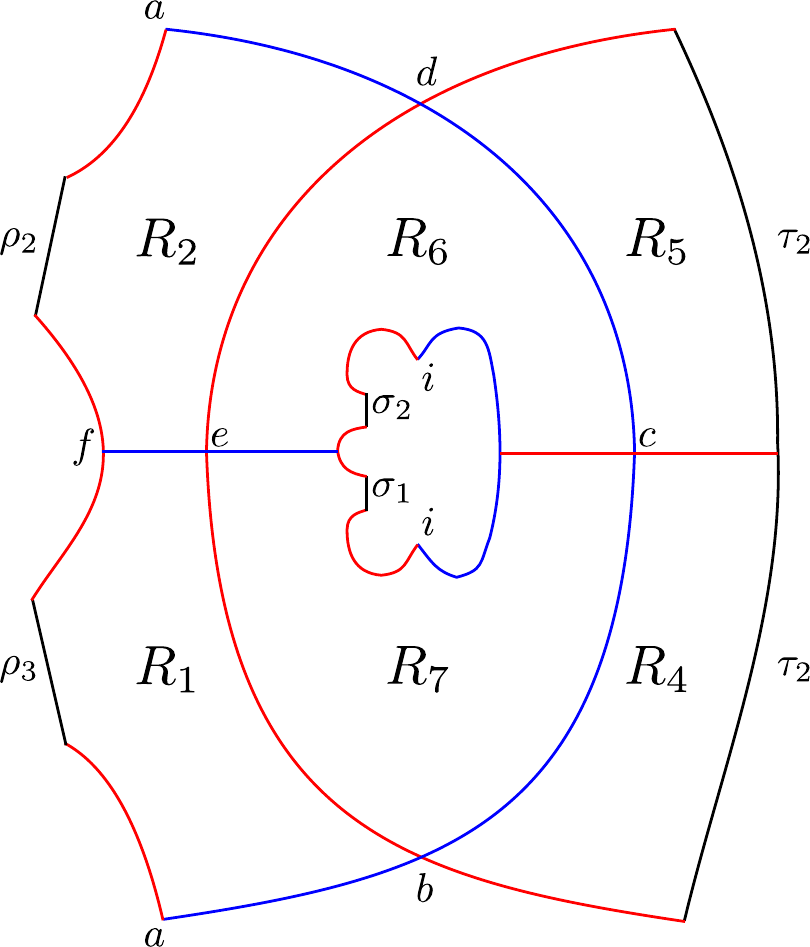}
\caption{The index 0 annulus $R_6 R_7$, with four surrounding regions.}
\label{indexfig}
\end{center}
\end{figure}

$R_6 R_7$ is an index 0 annulus of the same form as $A$ in Proposition \ref{Levine_prop}. The domains $R_2 R_6 R_7$, $R_5 R_6 R_7$, $R_1 R_6 R_7$, and $R_4 R_6 R_7$ in $\HeegDiag'$ correspond to the domains $D_1$, $D_2$, $D_3$, and $D_4$, respectively, in Figure \ref{index0domain}.  By Proposition \ref{Levine_prop}, the contribution of these four domains is determined by the choice of complex structure $J$ on $\Sigma$, and in particular on the resulting ratios of lengths $\Theta_{6 7}^{i, i}$ and $\Theta_{6 7}^{d, b}$. Recall that we chose the complex structure on $\Sigma$ so that $\Theta_{67}^{d, b} > \Theta_{67}^{i, i}$.

It follows directly from the proof of Proposition \ref{Levine_prop} that $R_4 R_6 R_7$ and $R_5 R_6 R_7$ do not contribute to the differential for our choice of $J$. It is also true that $R_1 R_6 R_7$ contributes $\rho_3 \sigma_{12} \otimes x$ to $\partial w$ and $R_2 R_6 R_7$ contributes $\rho_2 \sigma_{12} \otimes v$ to $\partial x$, but for these domains the outer boundary of the annulus has too many $\alpha$ and $\beta$ segments to apply Proposition \ref{Levine_prop} directly. First we apply Proposition \ref{pinching} and pinch along the arcs in $R_1$ or $R_2$ that are parallel to $\alpha^\tau_2$; recall that the complex structure was chosen to be consistent with pinching these arcs. The annuli that result from pinching the arcs completely have holomorphic representatives by Proposition \ref{Levine_prop}, and so $R_1 R_6 R_7$ and $R_2 R_6 R_7$ contribute to the differential.

\begin{figure}
\begin{center}
\begin{overpic}[scale = .8]{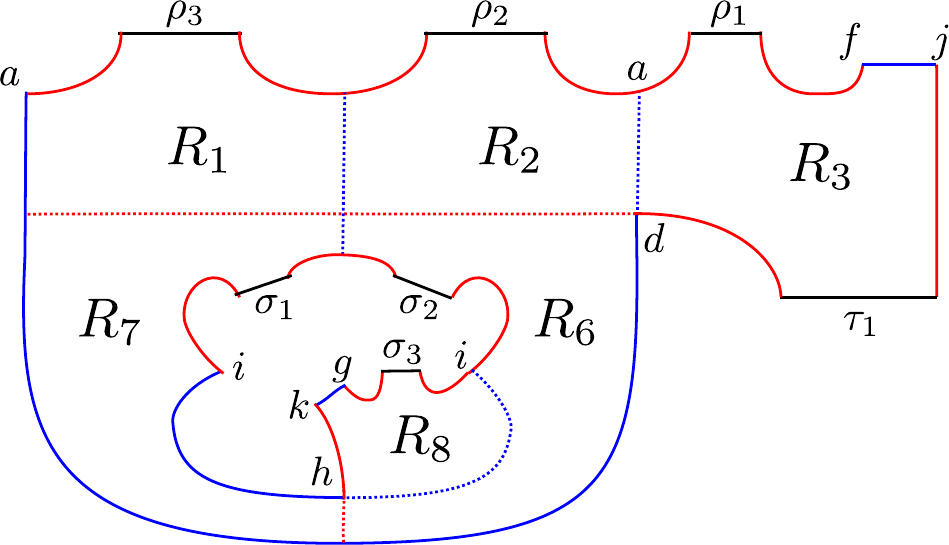}
\end{overpic}
\caption{The annulus $R_1 R_2 R_6 R_7 R_8$ contributes to the differential, while $R_1 R_2 R_3 R_6 R_7$ does not.
}
\label{regions_123678}
\end{center}
\end{figure}

We need to perform a similar analysis on two more domains, which are obtained from adding regions to the index 0 annulus  $R_1 R_2 R_6 R_7$ (see Figure \ref{regions_123678}). If we let the inner boundary in Figure \ref{regions_123678} correspond to the outer boundary in Figure \ref{index0domain}, then $R_8$ is analogous to $B_4$. $R_8$ is not a bigon, however, so we must first use Proposition \ref{pinching} to notice that the contribution of $R_1 R_2 R_6 R_7 R_8$ is the same as the contribution of the annulus which would result from collapsing an arc in $R_8$ which connects $\alpha^\sigma_1$ to $\alpha^\tau_1$. This pinched annulus would contribute \emph{with the Reeb chords $(\rho_2, \rho_3, \sigma_1, \sigma_2, \sigma_3)$} by Proposition \ref{Levine_prop}, using the fact that $\Theta_{1267}^{a, a} > \Theta_{1267}^{i, i}$. Thus the pair $(R_1 R_2 R_6 R_7, (\rho_2, \rho_3, \sigma_1, \sigma_2, \sigma_3))$ contributes to $\CFD{k}$. We emphasize however that this domain has a second compatible sequence of Reeb chords, about which Proposition \ref{Levine_prop} says nothing. The contribution of this domain with $\vecrho = (\rho_2, \rho_3, \sigma_{123})$ will be discussed in section \ref{decomposable_boundaries}.

The domain $R_1 R_2 R_3 R_6 R_7$, with Reeb chords $(\rho_1, \rho_2, \rho_3, \sigma_1, \sigma_2, \tau_1)$, can be analyzed in a similar way. By Propostion \ref{pinching} we will treat $R_3$ as a bigon attached to $R_1 R_2 R_6 R_7$, since the extra $\beta$ segment can be removed by pinching along $\gamma_3$. In this case Proposition \ref{Levine_prop} does not apply, because the arc $\alpha^\tau_2$ cutting into the annulus from the obtuse corner does not leave the annulus on the opposite boundary component, but the reasoning is similar. There is a one parameter family of cuts starting at $d$. We can define the ratios $\theta_{12367}^{d, a}$ and $\theta_{12367}^{i, i}$, which depend on the cutting parameter, as in the proof of Proposition \ref{Levine_prop}. There is a holomorphic representative for each zero of $\theta_{12367}^{d, a} - \theta_{12367}^{i, i}$. Cutting along $\alpha_2^\tau$ from $d$ to $b$ pinches off the annulus $R_6 R_7$. In this extreme, $\theta_{12367}^{d, a} - \theta_{12367}^{i, i}$ approaches $\Theta_{67}^{d, b} - \Theta_{67}^{i, i}$, which is positive for our choice of complex structure. In the other extreme, cutting along $\beta_1$ from $d$ to $a$ pinches off the annulus $R_1 R_2 R_6 R_7$, and $\theta_{12367}^{d, a} - \theta_{12367}^{i, i}$ becomes $\Theta_{1267}^{a, a} - \Theta_{1267}^{i, i} > 0$. Since $\theta_{12367}^{d, a} - \theta_{12367}^{i, i}$ has the same sign at both extremes, the number of zeros is even and the pair ($R_1 R_2 R_3 R_6 R_7, (\rho_1, \rho_2, \rho_3, \sigma_1, \sigma_2, \tau_1))$ does not contribute to the differential.

\subsection{Decomposable Boundaries}
\label{decomposable_boundaries}

We have determined the contribution of all domains in Table \ref{table1} except for the following seven:

\begin{center}
\begin{tabular}{>$c<$ >{ from }c >$c<$ >{to}c >$c<$}
R_1 R_2 R_3 & & v & & s \\
R_3 R_4 R_5 & & v & & s \\
R_1 R_2 R_3 R_6 R_7& & w & & s \\
R_3 R_4 R_5 R_6 R_7 & & w & & s \\
R_1 R_2 R_6 R_7 R_8 & & t & & y \\
R_4 R_5 R_6 R_7 R_8 & & t & & y \\
R_1 R_2 R_3 R_4 R_5 R_6 R_7 & & v & & s
\end{tabular}
\end{center}

Each of these domains is compatible with multiple Reeb chord sequences. The contribution of each domain/Reeb chord pair must be considered separately.

\begin{figure}
\begin{center}
\begin{tabular}{c c c}
\includegraphics[scale = 1]{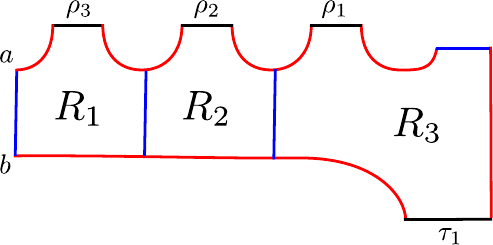}  & \qquad  &
\includegraphics[scale = 1]{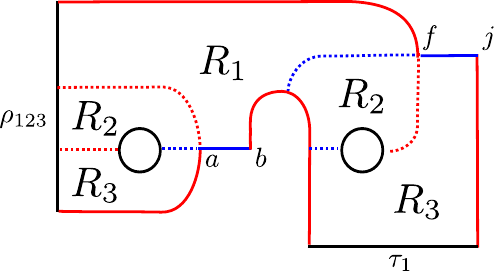} \\
$(a)$ & & $(b)$ \end{tabular}
\caption{Two versions of the domain $R_1 R_2 R_3$. $(a)$ represents the sequence of Reeb chords $\vec\eta_1 = (\rho_3, \rho_2, \rho_1, \tau_1)$. This domain/Reeb chord pair is a polygon and contributes to the differential. $(b)$ represents the sequence of Reeb chords $\vec\eta_2 = (\rho_{123}, \tau_1)$. This genus one domain also contributes to the differential}
\label{123versions}
\end{center}
\end{figure}

\bigbreak
${\bf R_1 R_2 R_3}$: 
The domain $R_1 R_2 R_3$ is compatible with both $\vec\eta_1 = (\rho_1, \rho_2, \rho_3, \tau_1)$ and $\vec\eta_2 = (\rho_{123}, \tau_1)$. By cutting along $\alpha$ arcs, the domain can be represented differently for each Reeb chord sequence (see Figure \ref{123versions}). Figure \ref{123versions}(a) is an immersed polygon; it is clear that the conditions of Propositions \ref{polygons} are satisfied, and so the pair $(R_1, R_2, R_3, \vec\eta_1)$ contributes to the differential.

For $\vec\eta_2$, we use Lemma \ref{Ainfty_trick} and consider the $\Ainfty$ relation for $(v, \rho_1, \rho_{23}, \tau_3)$:
\[ 0 = m\left( v, \mu(\rho_1, \rho_{23}), \tau_3 \right) + m \left( m (v, \rho_1), \rho_{23}, \tau_3 \right) . \]
There are no other nonzero terms in the relation. Note that it is impossible to have an $\Ainfty$ operation involving $\tau_3$ and not $\rho_3$, since both Reeb chords are on the same region $R_3$. Thus the term $m\left( m(v, \rho_1, \tau_3), \rho_{23} \right)$ does not appear in the $\Ainfty$ relation. Since we use $\FF_2$ coefficients, the relation above can be rewritten as
\[ m(v, \rho_{123}, \tau_3) = m \left( m (v, \rho_1), \rho_{23}, \tau_3 \right).  \]
The inner operation on the right, $m(v, \rho_1)$, records the contribution of the domain $R_1$ with $\vecrho = (\rho_1)$. We showed that this pair contributes in $\CFD{3}$, and so by Lemma \ref{Ainfty_trick} it also contributes to $\CFA{3}$. Thus $m(v, \rho_1) = x$. The outer operation is determined by the contribution of of the domain $R_2 R_3$. This domain was eliminated from consideration for $\CFD{3}$, but it may still contribute to $\CFA{3}$. To find out if it does we use another $\Ainfty$ relation, this time for $(x, \rho_2, \rho_3, \tau_3)$. The relation implies that
\[ m(x, \rho_{23}, \tau_3) = m \left( m(x, \rho_2), \rho_3, \tau_3 \right). \]
Since $R_2$ and $R_3$ are known to contribute in $\CFD{k}$ (and thus in $\CFA{k}$), we find that $m(x, \rho_{23}, \tau_3) = m(w, \rho_3, \tau_3) = s$, and $m(v, \rho_{123}, \tau_3) = s$. By Lemma \ref{Ainfty_trick} the pair $(R_1 R_2 R_3), \vec\eta_2)$ contributes to $\CFD{k}$.

$(R_1 R_2 R_3, \vec\eta_1)$ and $(R_1 R_2 R_3, \vec\eta_2)$ both contribute the term $\rho_{123} \tau_1 s$ to $\partial(v)$ in $\CFD{k}$. Over $\FF_2$, these contributions cancel, so the total contribution of $R_1 R_2 R_3$ to the differential is zero.

\bigbreak
${\bf R_3 R_4 R_5}$:
The two compatible Reeb chord sequences are $\vec\eta_1 = (\rho_1, \tau_1, \tau_2, \tau_3)$ and $\vec\eta_2 = (\rho_1, \tau_{123})$. The first does not contribute, because a holomorphic representative consistent with $\vec\eta_1$ would project to $R_3 R_4 R_5$ with cuts along $\alpha^\tau_1$ and $\alpha^\tau_2$, making the domain disconnected. If we do not cut the domain along $\alpha^\tau_1$ and $\alpha^\tau_2$, we see that it is an immersed polygon compatible with $\vec\eta_2$, and Proposition \ref{polygons} tells us that $(R_3 R_4 R_5, \vec\eta_2)$ contributes to the differential. Overall, the domain $R_3 R_4 R_5$ contributes $\rho_1 \tau_{123} s$ to $\partial(v)$.

\bigbreak
${\bf R_1 R_2 R_3 R_6 R_7}$:
This domain is compatible with $\vec\eta_1 = (\rho_1, \rho_2, \rho_3, \tau_1, \sigma_1, \sigma_2)$ and $\vec\eta_2 = (\rho_{123}, \tau_1, \sigma_1, \sigma_2)$. It was shown in the previous section that there is a contribution from $\vec\eta_1$. For the other case consider the $A_\infty$ relation for $(w, \rho_1, \rho_{23}, \tau_3, \sigma_3, \sigma_2)$, which says that $m \left( w, \mu(\rho_1, \rho_{23}), \tau_3, \sigma_2, \sigma_1\right) = m \left( m(w, \rho_1, \sigma_3, \sigma_2), \rho_{23}, \tau_3\right)$. The inner operation on the right is nontrivial because we have shown that the domain $R_1 R_6 R_7$ contributes to $\CFD{3}$. The inner operation evaluates to $x$. The outer operation, $m( x, \rho_{23}, \tau_3)$, evaluates to $s$ as shown above. Thus $m( w, \rho_{123}, \tau_3, \sigma_2, \sigma_1) = s$, and by Lemma \ref{Ainfty_trick} the pair $(R_1 R_2 R_3 R_6 R_7, \vec\eta_2)$ contributes to $\CFD{3}$. The total mod 2 contribution of $R_1 R_2 R_3 R_6 R_7$ is zero.

\bigbreak
${\bf R_3 R_4 R_5 R_6 R_7}$:
This domain is compatible with $\vec\eta_1 = (\rho_1, \tau_1, \tau_2, \tau_3, \sigma_1, \sigma_2)$ and $\vec\eta_2 = (\rho_1, \tau_{123}, \sigma_1, \sigma_2)$. $\vec\eta_1$ does not contribute because realizing the domain with boundary Reeb chords $\vec\eta_1$ would involve cutting along $\alpha^\tau_2$, which leaves the domain disconnected. For the contribution of $\vec\eta_2$, consider the $A_\infty$ relation for $(w, \rho_3, \tau_{12}, \tau_3, \sigma_3, \sigma_2)$. To find all terms of the relation, first note that any type $A$ operation that involves $m(\x, \ldots, \tau_{12}, \tau_3, \ldots)$ will be trivial for boundary monotonicity reasons. As a result, any term in the $\Ainfty$ relation splits as an operation involving $\tau_{12}$ and an operation involving $\tau_3$. We look in Table \ref{table1} for domains that connect $w$ to another generator and involve $R_4$ and $R_5$, but not $R_3$, and some subset of $\{R_6, R_7\}$. The only option is the domain $R_4 R_5 R_6 R_7$ connecting $w$ to itself. The relation can be written as
\[ m(w, \rho_3, \tau_{123}, \sigma_3, \sigma_2) = m\left( m( w, \tau_{12}, \sigma_3, \sigma_2), \rho_3, \tau_3 \right). \]
To compute the inner operation on the right, we can use another $\Ainfty$ relation to show that
\[ m( w, \tau_{12}, \sigma_3, \sigma_2) = m\left( m(w, \tau_1, \sigma_3, \sigma_2), \tau_2 \right) = m( 0, \tau_2) = 0, \]
where $m(w, \tau_1, \sigma_3, \sigma_2) = 0$ because the domain $R_4 R_6 R_7$ was shown not to contribute to $\CFD{3}$ based on the choice of complex structure $J$. Thus $(R_3 R_4 R_5 R_6 R_7, \vec\eta_2)$ does not contribute to $\CFA{3}$ or $\CFD{3}$, and overall the domain $R_3 R_4 R_5 R_6 R_7$ does not contribute.

\bigbreak
${\bf R_1 R_2 R_6 R_7 R_8}$:
The compatible Reeb chord sequences are $\vec\eta_1 = (\rho_2, \rho_3, \tau_1, \sigma_1, \sigma_2, \sigma_3)$ and $\vec\eta_2 = (\rho_2, \rho_3, \tau_1, \sigma_{123})$. It has already been shown that there is no contribution with $\vec\eta_1$. For the contribution of $\vec\eta_2$, use the $\Ainfty$ relation for $(t, \rho_2, \rho_1, \sigma_{12}, \sigma_3)$. Look for in Table \ref{table1} for domains starting at $t$ which involve $R_6$ and $R_8$ but not $R_7$, and some subset of $\{ R_1, R_2 \}$; there is only one option. The $\Ainfty$ relation becomes
\[ m(t, \rho_2, \rho_1, \mu(\sigma_{12}, \sigma_3) ) = m( m(t, \rho_2, \sigma_{12}), \rho_1, \sigma_3) = m(z, \rho_1, \sigma_3) = y. \]
Thus by Lemma \ref{Ainfty_trick}, the domain $R_1 R_2 R_6 R_7 R_8$ contributes $\rho_{23} \sigma_{123} y$ to $\partial( t )$.

\bigbreak
${\bf R_4 R_5 R_6 R_7 R_8}$:
This domain is compatible with $\vec\eta_1 = (\tau_2, \tau_3, \sigma_1, \sigma_2, \sigma_3)$ and $\vec\eta_2 = (\tau_2, \tau_3, \sigma_{123})$. With $\vec\eta_1$, the $R_4 R_5 R_6 R_7 R_8$ is realized as an immersed polygon, so the pair contributes to the differential. With $\vec\eta_2$, the contribution of this domain is determined by the $\Ainfty$ relation
\[ m( t, \tau_2, \tau_1, \mu(\sigma_{12}, \sigma_3) ) = m( m( t, \tau_2, \sigma_{12}), \tau_1, \sigma_3) = m(x, \tau_1, \sigma_3) = y \]
and Lemma \ref{Ainfty_trick}. The domain contributes with $\vec\eta_2$, and the total mod 2 contribution of the domain is zero.

\bigbreak
${\bf R_1 R_2 R_3 R_4 R_5 R_6 R_7}$:
This domain has four compatible sequences of Reeb chords: $\vec\eta_1 = (\rho_1, \rho_2, \rho_3, \tau_1, \tau_2, \tau_3, \sigma_1, \sigma_2)$, $\vec\eta_2 = (\rho_1, \rho_2, \rho_3, \tau_{123}, \sigma_1, \sigma_2)$, $\vec\eta_3 = (\rho_{123}, \tau_1, \tau_2, \tau_3, \sigma_1, \sigma_2)$, and $\vec\eta_4 = (\rho_{123}, \tau_{123}, \sigma_1, \sigma_2)$. To obtain a boundary with Reeb chords $\vec\eta_1$, we must cut along all $\alpha$ arcs. This produces an immersed polygon, but there are too many edges and corners for Proposition \ref{polygons} to apply. For instance, cutting along $\alpha^\tau_1$ produces corners at $h$, but the generators $v$ and $s$ do not contain the point $h$. Therefore the domain does not contribute with $\vec\eta_1$. For the other sequences of Reeb chords, we can use Lemma \ref{Ainfty_trick} and appropriate $\Ainfty$ relations. We find that the domain contributes with $\vec\eta_2$ and $\vec\eta_4$, and not with $\vec\eta_3$. Overall with $\FF_2$ coefficients the contribution of this domain is zero.

\subsection{Canceling Differentials}

\begin{figure}

\begin{center}
\begin{tikzpicture}[>=latex, scale=.9]

\node (y)  at (0,0) {$y$};
\node (v)  at (-5,0) {$v$};
\node (w) at (5,0) {$w$};
\node (x)  at (0,-5) {$x$};
\node (z)  at (0,5) {$z$};

\node (t) at (-6, 5) {$t$};
\node (s) at (6, 5) {$s$};

\draw[->, bend right = 20] (v) to node[pos = .6, below = 4pt] {\scriptsize $\rho_1 \tau_{123}$} (s);

\draw[->, bend left=8] (v) to node[pos = .75, above left = -2pt] {\scriptsize $\tau_3$} (z);

\node[fill = white] at (-1, .6) {};
\node[fill = white] at (-1.3, .5) {};
\node[fill = white] at (-.1, .9) {};
\node[fill = white] at (.1, 1) {};
\node[fill = white] at (3, 2.4) {};

\node[fill = white] at (-3.1, 2.1) {};
\node[fill = white] at (-3.4, 1.8) {};

\draw[->, bend right=8] (t) to node[pos = .25, above right = -2pt] {\scriptsize $\rho_{23} \sigma_{123}$} (y);
\draw[->, bend right=18] (t) to node[pos = .25, below left = -2pt] {\scriptsize $\sigma_3$}  (y);

\draw[->, bend right=8] (v) to node[below left = -2pt] {\scriptsize $\rho_3$} (x);
\draw[->, bend right=8] (x) to node[above right = -2pt] {\scriptsize $\rho_2 \sigma_{12}$}  (v);

\draw[->, bend right=10] (w) to node[right = -1pt] {\scriptsize $\rho_1 \tau_1$}  (s);

\draw[->, bend right=8] (w) to node[above left = -3pt] {\scriptsize $\rho_3 \sigma_{12}$} (x);
\draw[->, bend right=8] (x) to node[below right = -2pt] {\scriptsize $\rho_2$}  (w);

\draw[->, bend right=10] (y) to node[left = -2pt] {\scriptsize $\tau_2 \sigma_2$} (x);
\draw[->, bend right=10] (x) to node[right = -2pt] {\scriptsize $\tau_3 \sigma_1$}  (y);

\draw[->, bend right=10] (z) to node[left = -2pt] {\scriptsize $\rho_3 \sigma_1$} (y);
\draw[->, bend right=10] (y) to node[right = -2pt] {\scriptsize $\rho_2 \sigma_2$}  (z);

\draw[->, bend left=8] (z) to node[pos = .25, above right = -2pt] {\scriptsize $\tau_2$}  (w);

\draw[->, bend left = 20] (t) to (s);

\end{tikzpicture}

\end{center}

 \caption{$\CFD{k}(\HeegDiag')$ in the middle spin$^c$-structure, for the given choice of complex structure $J$.}
\label{answer}
\end{figure}
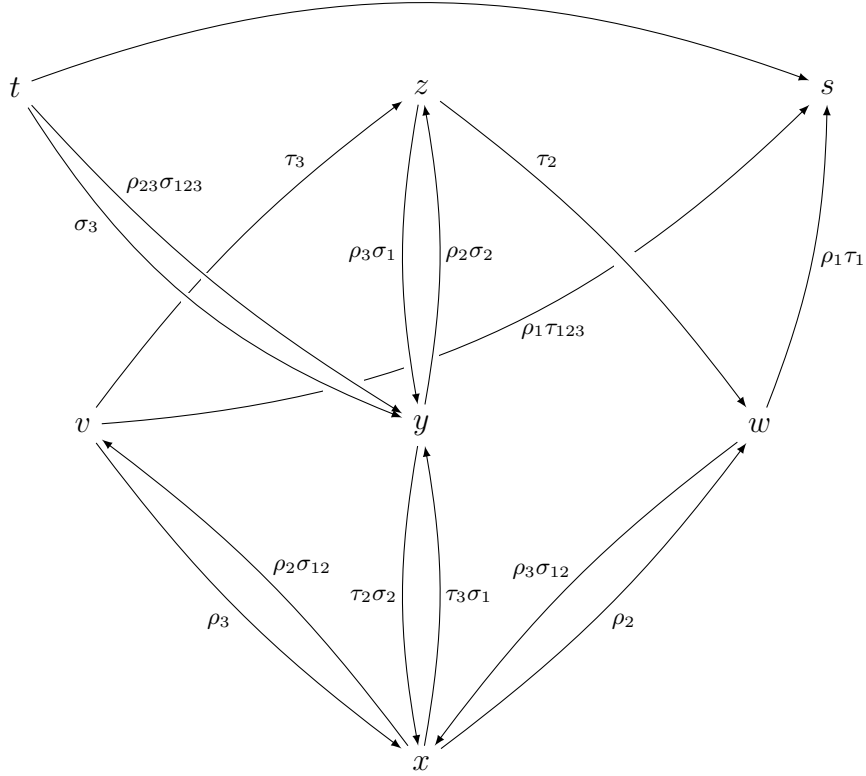

\begin{figure}

\begin{center}
\begin{tikzpicture}[>=latex, scale=.9]

\node (y)  at (0,0) {$y$};
\node (v)  at (-5,0) {$v$};
\node (w) at (5,0) {$w$};
\node (x)  at (0,-5) {$x$};
\node (z)  at (0,5) {$z$};

\draw[->, bend left=10] (v) to node[above] {\scriptsize $\rho_1 \tau_{123} \sigma_3$} (y);
\draw[->, bend right=10] (v) to node[below] {\scriptsize $\rho_{123} \tau_{123} \sigma_{123}$}  (y);

\draw[->, bend right=8] (v) to node[below left = -2pt] {\scriptsize $\rho_3$} (x);
\draw[->, bend right=8] (x) to node[above right = -2pt] {\scriptsize $\rho_2 \sigma_{12}$}  (v);

\draw[->, bend left=10] (w) to node[below] {\scriptsize $\rho_{123} \tau_1 \sigma_{123}$} (y);
\draw[->, bend right=10] (w) to node[above] {\scriptsize $\rho_1 \tau_1 \sigma_3$}  (y);

\draw[->, bend right=8] (w) to node[above left = -3pt] {\scriptsize $\rho_3 \sigma_{12}$} (x);
\draw[->, bend right=8] (x) to node[below right = -2pt] {\scriptsize $\rho_2$}  (w);

\draw[->, bend right=10] (y) to node[left = -2pt] {\scriptsize $\tau_2 \sigma_2$} (x);
\draw[->, bend right=10] (x) to node[right = -2pt] {\scriptsize $\tau_3 \sigma_1$}  (y);

\draw[->, bend right=10] (z) to node[left = -2pt] {\scriptsize $\rho_3 \sigma_1$} (y);
\draw[->, bend right=10] (y) to node[right = -2pt] {\scriptsize $\rho_2 \sigma_2$}  (z);

\draw[->, bend left=8] (v) to node[above left = -2pt] {\scriptsize $\tau_3$} (z);
\draw[->, bend left=8] (z) to node[above right = -2pt] {\scriptsize $\tau_2$}  (w);

\end{tikzpicture}

\caption{$\CFD{3}(\Y)$ in the middle spin$^c$-structure after canceling the differential from $t$ to $s$ in $\CFD{3}(\HeegDiag')$.}
\label{answer_simplified}

\end{center}
\end{figure}
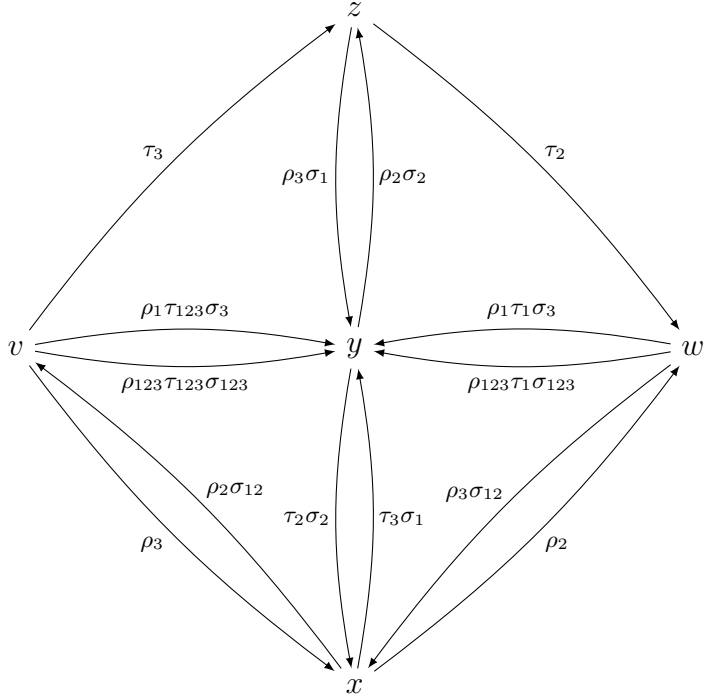

Putting everything together, the differential on $\CFD{k}(\Y)$ is recorded in Figure \ref{answer}. The unlabeled arrow from $t$ to $s$ is the differential corresponding to the bigon $R_9$. This unlabeled edge can be canceled using the edge reduction algorithm for type $D$ structures described in Remark \ref{edge_reduction}. We eliminate the arrow and the generators $t$ and $s$, and for every ``zig-zag''
\[ {\x} \overset{a_1}\longrightarrow s \longleftarrow t \overset{a_2}\longrightarrow \y \]
we introduce the new arrow ${\x} \overset{a_1 \cdot a_2}\longrightarrow {\y}$. The resulting simplified form of $\CFD{k}(\Y)$ (which is quasi-isomorphic to the first diagram) is given in Figure \ref{answer_simplified}.

\subsection{Extremal spin$^c$-structures}

To complete the computation of $\CFD{3}(\Y)$, we must compute $\CFD{3}(\Y, \spinc)$ for other spin$^c$-structures $\spinc$.

\bigbreak
{\bf (1, 2, 0):}
Consider first the spin$^c$-structure $\spinc$ that has 1, 2, and 0 $\alpha$ arcs occupied on the $\rho$, $\sigma$, and $\tau$ boundaries, respectively. The only generator in this spin$^c$-structure is $agi$, so $\CFD{3}(\Y, \spinc)$ has one generator and no differentials.

\bigbreak
{\bf (2, 1, 0) and (2, 0, 1):}
The spin$^c$-structures with two arcs occupied on the $\rho$ boundary each have only one generator, and no differentials. The generator with one $\alpha^\sigma$ arc is occupied is $afi$, and the generator with one $\alpha^\tau$  arc occupied is $afh$.

\bigbreak
{\bf (1, 0, 2):}
The spin$^c$-structure $\spinc$ that has one $\alpha$ arc occupied on the $\rho$ boundary, zero on the $\sigma$ boundary, and two on the $\tau$ boundary has 3 generators: $aeh$, $bfh$, and $dfh$. Recall that domains connecting these generators may not touch the $\sigma$ or $\tau$ boundaries, so we only need to consider the regions $R_1$, $R_2$, and $R_9$. It is easy to see that $R_1$ contributes a differential from $bfh$ to $aeh$ and $R_2$ contributes a differential from $aeh$ to $dfh$ (they are both polygons). None of the generators include the points $i$ or $j$, so $R_9$ is not involved in any differentials. Finally, $R_1 R_2$ does not contribute by Proposition \ref{idempotent_proposition}; the corners make $R_1 R_2$ a domain connecting $bfh$ to $dfh$, but neither of these generators contain a point on $\alpha^\rho_1$, which is required for a domain with Reeb chords $(\rho_2, \rho_3)$ to contribute nontrivially. For this spin$^c$-structure, $\CFD{3}(\Y, \spinc)$ is given by
\[ bfh  \overset{\rho_3}\longrightarrow aeh \overset{\rho_2}\longrightarrow dfh. \]

\bigbreak
{\bf (0, 2, 1):}
This spin$^c$-structure has 3 generators:  $bgi$, $dgi$, and $dgi$. Domains that contribute to the differential do not touch the $\rho$ or $\sigma$ boundaries. The only domains which connect two generators are $R_4$, $R_5$, and $R_4 R_5$. It is clear that the polygons $R_4$ and $R_5$ contribute, but the contribution of $R_4 R_5$ is zero. Thus $\CFD{3}(\Y, \spinc)$ for this spin$^c$-structure is given by
\[ bgi  \overset{\tau_3}\longrightarrow cgi \overset{\tau_2}\longrightarrow dgi. \]

\bigbreak
{\bf (0, 1, 2):}
This spin$^c$-structure has 7 generators:  $bgh$, $dgh$, $cei$, $bji$, $bki$, $djk$, and $dki$. $R_9$ contributes differentials from $bki$ to $bji$ and from $dki$ to $dji$. There can be no other differentials ending at $bji$ or $dji$, so the edge reduction algorithm lets us cancel these differentials and remove the generators $bki$, $bji$, $dki$, and $dji$ without adding new differentials (it is worth noting that we could also compute $\CFD{3}(\Y, \spinc)$ using the Heegaard diagram $\HeegDiag$ in Figure \ref{HeegDiag_simple} instead of $\HeegDiag'$, and we would not have to deal with these four generators at all). The only domains connecting the remaining three generators are $R_6$, $R_7$, and $R_6 R_7$. Once again, the individual regions contribute while $R_6 R_7$ has zero contribution for idempotent reasons, and $\CFD{3}(\Y, \spinc)$ is
\[ dgh  \overset{\sigma_2}\longrightarrow cei \overset{\sigma_1}\longrightarrow bgh. \]

\subsection{Gradings}

As described in Section \ref{sec:gradings}, $\CFD{3}(\Y)$ is graded by a set which is a quotient of the noncommutative group $G_{3,0}$. We will compute this (relative) grading for the middle spin$^c$-structure, using the form of $\CFD{3}(\Y)$ depicted in Figure \ref{answer_simplified}.

We choose $x$ to be the preferred generator and assign it the grading $\vec0 = (0; 0, 0; 0, 0; 0, 0)$. The arrow from $x$ to $w$ indicates that $\partial x$ contains the term $\rho_2 w$, and thus $gr(\partial x) = gr(\rho_2 w)$. By the relation \ref{grading_relation}, we have
\[ \lambda^{-1} gr(x) = gr(\partial x) = gr(\rho_2) gr(w), \]
and so
\[ gr(w) = \lambda^{-1} gr(\rho_2)^{-1} gr(x) = \lambda^{-1} \left( \frac{1}{2}; -\frac{1}{2}, -\frac{1}{2}; 0, 0; 0, 0 \right) \vec0 = \left( -\frac{1}{2}; -\frac{1}{2}, -\frac{1}{2}; 0, 0; 0, 0 \right). \]
Similarly, the arrow from $z$ to $w$ implies that $\lambda^{-1} gr(z) = gr(\partial z) = gr(\tau_2) gr(w)$, so 
\begin{eqnarray*}
gr(z) &=& \lambda \left( -\frac{1}{2}; 0, 0; 0, 0; \frac{1}{2}, \frac{1}{2} \right) \left( -\frac{1}{2}; -\frac{1}{2}, -\frac{1}{2}; 0, 0; 0, 0 \right) \\
&=& \left( 0;  -\frac{1}{2}, -\frac{1}{2}; 0, 0; \frac{1}{2}, \frac{1}{2} \right)
\end{eqnarray*}
The arrow from $v$ to $z$ implies that
\begin{eqnarray*}
gr(v) &=& \lambda gr(\tau_3) gr(z) \\
&=& \lambda \left( -\frac{1}{2}; 0, 0; 0, 0; -\frac{1}{2}, \frac{1}{2} \right) \left( 0;  -\frac{1}{2}, -\frac{1}{2}; 0, 0; \frac{1}{2}, \frac{1}{2} \right) \\
&=& \left( 0;  -\frac{1}{2}, -\frac{1}{2}; 0, 0; 0, 1 \right).
\end{eqnarray*}
Finally, the arrow from $y$ to $x$ implies that
\begin{eqnarray*}
gr(y) &=& \lambda gr(\sigma_2) gr(\tau_2) gr(x) \\
&=& \lambda \left( -\frac{1}{2}; 0, 0; \frac{1}{2}, \frac{1}{2}; 0, 0 \right) \left( -\frac{1}{2}; 0, 0; 0, 0; \frac{1}{2}, \frac{1}{2} \right) \vec0 \\
&=& \left( 0; 0, 0; \frac{1}{2}, \frac{1}{2}; \frac{1}{2}, \frac{1}{2} \right).
\end{eqnarray*}
We have now computed the gradings of each generator as elements of the group $G_{3,0}$. However, these gradings are only well defined modulo the action of $\mathcal{P}(x)$, the group generated by the gradings of periodic domains connecting $x$ to itself. To finish the computation, we need to find $\mathcal{P}(x)$.

Consider the arrow from $v$ to $x$, which implies that $ \lambda^{-1} gr(v) = gr(\rho_3) gr(x)$. It follows that
\begin{eqnarray*}
gr(x) &=& \lambda^{-1} gr(\rho_3)^{-1} gr(v) \\
&=& \lambda^{-1} \left( \frac{1}{2}; \frac{1}{2}, -\frac{1}{2}; 0, 0; 0, 0 \right)  \left( 0;  -\frac{1}{2}, -\frac{1}{2}; 0, 0; 0, 1 \right) \\
&=& \left( -1; 0, -1; 0, 1; 0, 0 \right).
\end{eqnarray*}
We have that $gr(x) = (-1; 0, -1; 0, 1; 0, 0)$, but also that $gr(x) = \vec0$. It follows that $(-1; 0, -1; 0, 1; 0, 0)$ and $\vec0$ are equivalent modulo $\mathcal{P}(x)$, and thus that $(-1; 0, -1; 0, 1; 0, 0) \in \mathcal{P}(x)$. In fact, since this nonzero value of $gr(x)$ was obtained from $\vec0$ by following a loop of edges with oriented labels $(\rho_2, -\tau_2, -\tau_3, \rho_3)$, the difference $(-1; 0, -1; 0, 1; 0, 0)$ corresponds to the grading of a periodic domain with boundary $\rho_{23} - \tau_{23}$.

Another value for $gr(x)$, and thus another element of $\mathcal{P}(x)$, can be found by considering the arrow from $x$ to $y$. We have that
\begin{eqnarray*}
gr(x) &=& \lambda gr(\sigma_1) gr(\tau_3) gr(y) \\
&=& \lambda \left( -\frac{1}{2}; 0, 0; \frac{1}{2}, -\frac{1}{2}; 0, 0 \right) \left( -\frac{1}{2}; 0, 0; 0, 0; -\frac{1}{2}, \frac{1}{2} \right) \left( 0; 0, 0; \frac{1}{2}, \frac{1}{2}; \frac{1}{2}, \frac{1}{2} \right) \\
&=& \left( 0; 0, 0; 1, 0; 0, 1 \right).
\end{eqnarray*}
So $(0; 0, 0; 1, 0; 0, 1)$ is an element of $\mathcal{P}(x)$, corresponding to a periodic domain with boundary $\sigma_{12} + \tau_{23}$.

Consider the loop formed by the arrow from $y$ to $x$, the $\rho_1 \sigma_3 \tau_{123}$ arrow from $w$ to $x$, and the arrow from $x$ to $w$. This loop corresponds to a periodic domain with boundary $\rho_{12} + \sigma_{23} + \tau_{12}$. As before, starting with $gr(x) = \vec0$ the arrow from $y$ to $x$ implies that $gr(y) = \left( 0; 0, 0; \frac{1}{2}, \frac{1}{2}; \frac{1}{2}, \frac{1}{2} \right)$. The arrow from $w$ to $y$ then implies that
\begin{eqnarray*}
gr(w) &=& \lambda gr(\rho_1) gr(\sigma_3) gr(\tau_1) gr(y) \\
&=& \lambda \left( -\frac{1}{2}; \frac{1}{2}, -\frac{1}{2}; 0, 0; 0, 0 \right) \left( -\frac{1}{2}; 0, 0; -\frac{1}{2}, \frac{1}{2}; 0, 0 \right) \left( -\frac{1}{2}; 0, 0; 0, 0; \frac{1}{2}, -\frac{1}{2} \right) gr(y) \\
&=&  \left( -\frac{1}{2}; \frac{1}{2}, -\frac{1}{2}; -\frac{1}{2}, \frac{1}{2}; \frac{1}{2}, -\frac{1}{2} \right)  \left( 0; 0, 0; \frac{1}{2}, \frac{1}{2}; \frac{1}{2}, \frac{1}{2} \right) \\
&=& \left( -\frac{1}{2}; \frac{1}{2}, -\frac{1}{2}; 0, 1; 1, 0 \right).
\end{eqnarray*}
The arrow from $x$ to $w$ then implies that
\begin{eqnarray*}
gr(x) &=& \lambda gr(\rho_2) gr(w) \\
&=& \lambda \left( -\frac{1}{2}; \frac{1}{2}, \frac{1}{2}; 0, 0; 0, 0 \right) \left( -\frac{1}{2}; \frac{1}{2}, -\frac{1}{2}; 0, 1; 1, 0 \right) \\
&=& \left( -\frac{1}{2}; 1, 0; 0, 1; 1, 0 \right)
\end{eqnarray*}
is an element of $\mathcal{P}(x)$.

Since $\Y$ has 3 boundary components, the space of periodic domains has dimension 3. Since the three elements of $\mathcal{P}(x)$ we have found are independent, they are enough to determine $\mathcal{P}(x)$:
\[ \mathcal{P}(x) = \left<  (-1; 0, -1; 0, 1; 0, 0), (0; 0, 0; 1, 0; 0, 1), \left( -\frac{1}{2}; 1, 0; 0, 1; 1, 0 \right) \right> \]

\section{Self Gluing}
\label{sec:self_gluer}
\FloatBarrier

Any graph manifold which is represented by a tree with only genus zero vertices can be obtained by gluing together copies of $\Y$, solid tori, and mapping cylinders of appropriate Dehn twists. Each time a new piece is glued on, the new bordered Heegaard Floer invariants can be obtained as a box tensor product by the pairing theorem.

\begin{figure}

\begin{center}

\begin{overpic}[scale = 1, tics = 10]{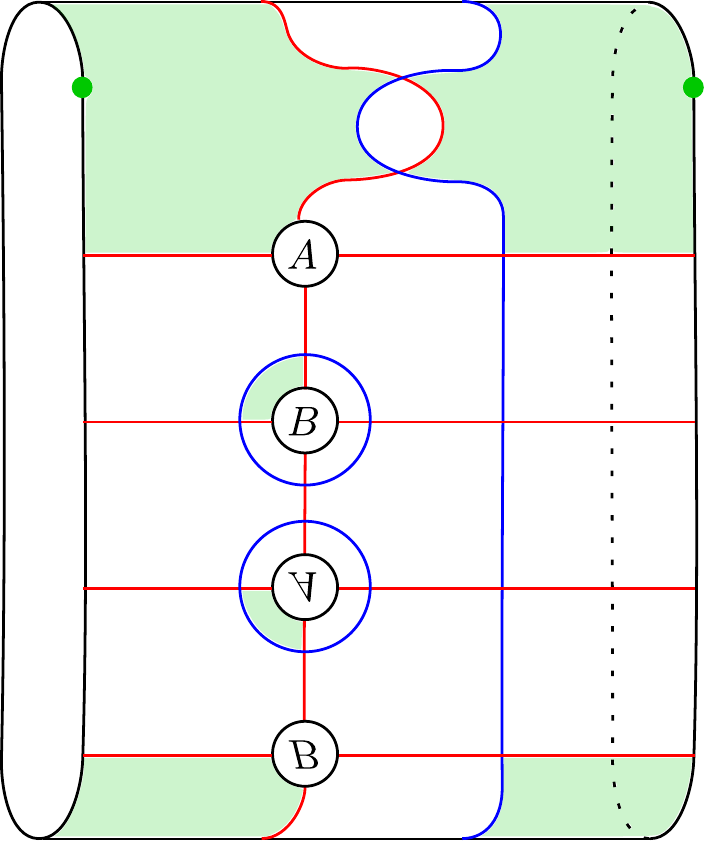}
\put(5, 19){$\rho_1$}
\put(5, 39){$\rho_2$}
\put(5, 59){$\rho_3$}
\put(84, 19){$\sigma_3$}
\put(84, 39){$\sigma_2$}
\put(84, 59){$\sigma_1$}

\put(20, 18){$R_1$}
\put(20, 38){$R_2$}
\put(20, 58){$R_3$}
\put(70, 58){$R_4$}
\put(70, 38){$R_5$}
\put(70, 18){$R_6$}

\put(45, 72){\small $R_7$}
\put(45, 58){$R_8$}
\put(45, 38){$R_9$}
\put(45, 18){$R_{10}$}
\put(45, 4){$R_{11}$}
\put(43, 96){\small $R_{11}$}
\put(44, 84){\footnotesize $R_{12}$}

\put(26, 51){\footnotesize $a$}
\put(37, 59){\footnotesize $b$}
\put(45, 51){\footnotesize $c$}
\put(37, 43){\footnotesize $d$}

\put(26, 31){\footnotesize $e$}
\put(37, 35){\scriptsize $f$}
\put(45, 31){\footnotesize $g$}
\put(37, 20){\footnotesize $h$}

\put(61, 71){\footnotesize $i$}
\put(61, 51){\footnotesize $j$}
\put(61, 31){\footnotesize $k$}
\put(61, 11){\footnotesize $\ell$}

\put(50.5, 79.5){\scriptsize $m$}
\put(42.5, 89.5){\scriptsize $n$}

\put(12, 12){\color{red} \scriptsize $\alpha^\rho_1$}
\put(12, 32){\color{red} \scriptsize $\alpha^\rho_2$}

\put(77, 66){\color{red} \scriptsize $\alpha^\sigma_1$}
\put(77, 46){\color{red} \scriptsize $\alpha^\sigma_2$}

\put(30, 93){\color{red} \scriptsize $\alpha_0$}

\put(60, 93){\color{blue} \scriptsize $\beta_0$}

\put(29, 41){\color{blue} \scriptsize $\beta_1$}
\put(29, 21){\color{blue} \scriptsize $\beta_2$}
\end{overpic}

\caption{A Heegaard diagram $\Hselfgluer$ for the self gluing bimodule. }
\label{self_gluer}
\end{center}
\end{figure}

To build up an arbitrary graph manifold from these building blocks, however, it is often necessary to glue two boundaries of one manifold together. The resulting bordered invariants are obtained by taking Hochschild homology. In this case, we must insert an additional bimodule, which corresponds to gluing a certain bordered Heegaard diagram $\Hselfgluer$ (Figure \ref{self_gluer}) between the two boundary components that are being glued. Strictly speaking, $\Hselfgluer$ is a bordered sutured diagram. This process is discussed in \cite[Section 4.4]{LipshitzTreumann}, and a Heegaard diagram isotopic to $\Hselfgluer$ is given there, but the bimodule associated to this Heegaard diagram is not computed. Let $\Yselfgluer$ denote the manifold represented by $\Hselfgluer$. The focus of the present section is to compute $\CFDD(\Yselfgluer)$.

\bigbreak

We first restrict our attention to the middle spin$^c$-structure, where exactly one $\alpha$ arc is occupied at each boundary component. In fact, this computation gives all of $\CFDD(\Hselfgluer)$; it will be shown at the end of this section that the other summands for $\CFDD(\Hselfgluer)$ are trivial.

\bigbreak
{\bf Complex structure:}
As always, the computation of $\CFDD(\Hselfgluer)$ depends on the complex structure $J$ chosen for the Heegaard surface $\Sigma$. We will make the following assumptions:

\begin{itemize}
\item
$\Theta_{8,9}^{d,b} < \Theta_{8,9}^{i,k}$;
\item
$\Theta_{9,10}^{h,f} < \Theta_{9,10}^{j,l}$;
\item
$\Theta_{2,3,8,9}^{a,a} < \Theta_{2,3,8,9}^{i,k}$;
\item
$\Theta_{1,2,9,10}^{e,e} < \Theta_{1,2,9,10}^{j,l}$.
\end{itemize}
We will also assume that the following arcs are sufficiently pinched to apply Proposition \ref{pinching} when necessary:

\begin{itemize}
\item
an arc through $R_1$, $R_2$ and $R_3$ parallel to $\alpha_0$ connecting $\beta_2$ in $R_1$ to $\beta_1$ in $R_3$;
\item
an arc in $R_7$ parallel to $\alpha_0$, connecting $\beta_0$ to $\beta_2$;
\item
an arc in $R_{11}$ parallel to $\alpha_0$, connecting $\beta_0$ to $\beta_1$.
\end{itemize}
Finally, we will also assume that $\Theta_{2,5,8,9,9,10}^{h, b} < \Theta_{2,5,8,9,9,10}^{i, \ell}$, where here $\Theta_{2,5,8,9,9,10}$ refers to the appropriate ratio of $\alpha$ and $\beta$ lengths for the annulus obtained from $R_2 R_5 R_8 R_9 R_9 R_{10}$ by pinching along the arc through $R_2$ mentioned above (see Figure \ref{regions_24589910}).

\bigbreak
{\bf Generators:}
$\CFDD(\Hselfgluer)$ has 20 generators in the middle spin$^c$ structure: $afi$, $afj$, $afk$, $afl$, $ahi$, $ahj$, $ahk$, $ahl$, $ang$, $amg$, $ebi$, $ebj$, $ebk$, $ebl$, $edi$, $edj$, $edk$, $edl$, $enc$, and $emc$.

\bigbreak
{\bf Domains:}
To list the domains that might contribute to $\CFDD(\Hselfgluer)$, note that the multiplicities of the regions on the boundary ($R_1, \ldots, R_6$) are at most 1, and the region $R_{12}$ cannot be combined with any other regions. Checking all positive connected domains which satisfy these conditions for appropriate corner multiplicity, we find 292 domains to consider. Of these, 200 can be eliminated by Proposition \ref{index_lemma} or Proposition \ref{idempotent_proposition} (though we should make note of these domains, in case they come up when checking $\Ainfty$ relations). All of these steps are easy to perform with a computer program. The remaining 92 domains are listed in Table \ref{selfgluer_domains}.

\begin{table}[htdp]
\begin{center}
\begin{tabular}{| lcc |}
\hline
Regions & $\x$ & $\y$ \\ \hline
1 & $edi$ & $ahi$\\
1 & $edj$ & $ahj$\\
1 & $edk$ & $ahk$\\
1 & $edl$ & $ahl$\\
2 & $afi$ & $edi$\\
2 & $afj$ & $edj$\\
2 & $afk$ & $edk$\\
2 & $afl$ & $edl$\\
3 & $ebi$ & $afi$\\
3 & $ebj$ & $afj$\\
3 & $ebk$ & $afk$\\
3 & $ebl$ & $afl$\\
4 & $afj$ & $afi$\\
4 & $ahj$ & $ahi$\\
4 & $ebj$ & $ebi$\\
4 & $edj$ & $edi$\\
5 & $afk$ & $afj$\\
5 & $ahk$ & $ahj$\\
5 & $ebk$ & $ebj$\\
5 & $edk$ & $edj$\\
6 & $afl$ & $afk$\\
6 & $ahl$ & $ahk$\\
6 & $ebl$ & $ebk$\\
6 & $edl$ & $edk$\\
7 & $amg$ & $ahi$\\
11 & $ebl$ & $enc$\\
12 & $amg$ & $ang$\\
12 & $emc$ & $enc$\\
1,2,3 & $ebi$ & $ahi$\\
1,2,3 & $ebj$ & $ahj$\\
1,2,3 & $ebk$ & $ahk$\\
\hline
\end{tabular}
\begin{tabular}{| lcc |}
\hline
Regions & $\x$ & $\y$ \\ \hline
1,2,3 & $ebl$ & $ahl$ \\
1,6,10 & $enc$ & $ang$\\
1,6,10 & $emc$ & $amg$\\
1,9,10 & $edj$ & $afl$\\
1,10,11 & $ebk$ & $ang$\\
2,5,9 & $ang$ & $enc$\\
2,5,9 & $amg$ & $emc$\\
2,8,9 & $afi$ & $ebk$\\
2,9,10 & $ahj$ & $edl$\\
3,4,8 & $enc$ & $ang$\\
3,4,8 & $emc$ & $amg$\\
3,7,8 & $emc$ & $ahj$\\
3,8,9 & $edi$ & $afk$\\
4,5,6 & $afl$ & $afi$\\
4,5,6 & $ahl$ & $ahi$\\
4,5,6 & $ebl$ & $ebi$\\
4,5,6 & $edl$ & $edi$\\
4,8,9 & $edj$ & $ebk$\\
4,8,11 & $afl$ & $ang$\\
5,8,9 & $edi$ & $ebj$\\
5,9,10 & $ahk$ & $afl$\\
6,7,10 & $emc$ & $edi$\\
6,9,10 & $ahj$ & $afk$\\
1,2,3,8,9 & $edi$ & $ahk$\\
1,2,3,9,10 & $ebj$ & $afl$\\
1,2,5,9,10 & $ebk$ & $ebl$\\
1,2,6,9,10 & $ebj$ & $ebk$\\
1,2,9,10,11 & $ebj$ & $enc$\\
1,5,6,9,10 & $edi$ & $afi$\\
1,8,9,10,11 & $edi$ & $ang$\\
2,3,4,8,9 & $ahj$ & $ahk$ \\
\hline
\end{tabular}
\begin{tabular}{| lcc |}
\hline
Regions & $\x$ & $\y$ \\ \hline
2,3,5,8,9 & $ahi$ & $ahj$\\
2,3,7,8,9 & $amg$ & $ahk$\\
2,4,5,8,9 & $afl$ & $ebl$\\
2,5,6,9,10 & $ahi$ & $edi$\\
3,4,5,8,9 & $edl$ & $afl$\\
3,7,8,9,10 & $emc$ & $afl$\\
4,5,6,8,9 & $edl$ & $ebk$\\
4,5,6,9,10 & $ahj$ & $afi$\\
4,5,8,9,11 & $edl$ & $enc$\\
4,8,9,10,11 & $ahj$ & $ang$\\
5,6,7,9,10 & $amg$ & $afi$\\
6,7,8,9,10 & $emc$ & $ebk$\\
7,8,9,10,11 & $amg$ & $ang$\\
7,8,9,10,11 & $emc$ & $enc$\\
1,2,3,4,5,8,9 & $edl$ & $ahl$\\
1,2,3,5,6,9,10 & $ebi$ & $afi$\\
1,2,3,7,8,9,10 & $emc$ & $ahl$\\
1,2,3,8,9,10,11 & $ebi$ & $ang$\\
1,2,4,5,6,9,10 & $ebj$ & $ebi$\\
1,2,5,8,9,9,10 & $edi$ & $ebl$\\
2,3,4,5,6,8,9 & $ahl$ & $ahk$\\
2,3,5,8,9,9,10 & $ahi$ & $afl$\\
2,4,5,8,9,9,10 & $ahj$ & $ebl$\\
2,5,6,8,9,9,10 & $ahi$ & $ebk$\\
2,5,7,8,9,9,10 & $amg$ & $ebl$\\
2,5,8,9,9,10,11 & $ahi$ & $enc$\\
4,5,6,7,8,9,10 & $emc$ & $ebi$\\
4,5,6,8,9,10,11 & $ahl$ & $ang$\\
1,2,3,4,5,6,8,9,10 & $enc$ & $ang$\\
1,2,3,4,5,6,8,9,10 & $emc$ & $amg$ \\
& & \\
\hline
\end{tabular}
\caption{List of 92 domains that might contribute to $\CFDD(\Hselfgluer)$, with the corresponding initial generators $\x$ and final generators $\y$.}
\label{selfgluer_domains}
\end{center}
\end{table}

\bigbreak
{\bf Polygons:}
All of the single region domains are easily seen to be polygons, and thus contribute to the differential by Proposition \ref{polygons}.  In addition, it is easy to check that the following domains are polygons: $R_1 R_6 R_{10}$, $R_1 R_{10} R_{11}$, $R_2 R_5 R_9$, $R_3 R_4 R_8$, $R_3 R_7 R_8$, $R_4 R_8 R_{11}$, and $R_6 R_7 R_{10}$. Each domain has only one sequence of Reeb chords to consider. Thus by Proposition \ref{polygons} each of these domains contributes, and we have quickly dispatched 38 of the entries in Table \ref{selfgluer_domains}.

\begin{figure}[htbp]
\begin{center}
\begin{overpic}[scale = 1]{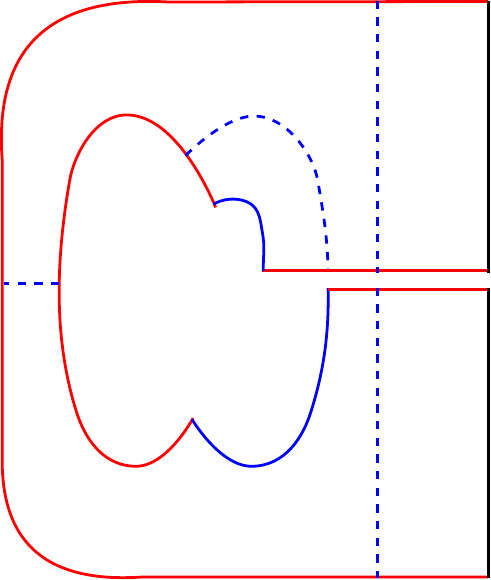}
\put(86, 75){$\sigma_1$}
\put(86, 25){$\sigma_2$}
\put(71, 75){$R_4$}
\put(71, 25){$R_5$}
\put(30, 85){$R_8$}
\put(30, 10){$R_9$}
\put(40, 69){$R_{11}$}

\put(53, 49){\scriptsize $c$}
\put(43, 49){\scriptsize $\ell$}
\put(33, 61){\scriptsize $n$}
\put(33, 29){\scriptsize $d$}

\end{overpic}
\qquad \qquad \qquad
\begin{overpic}[scale = 1]{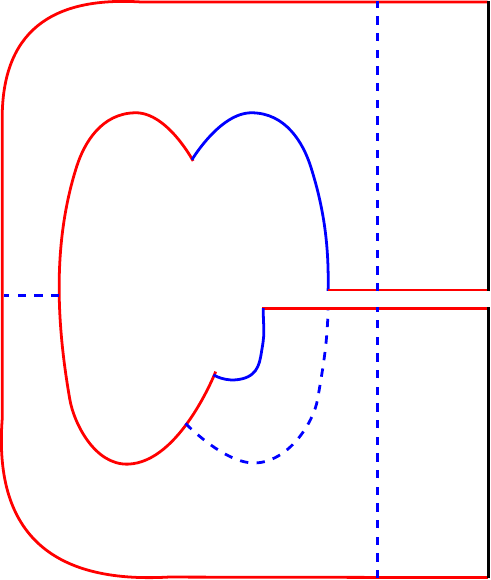}
\put(86, 75){$\sigma_2$}
\put(86, 25){$\sigma_3$}
\put(71, 75){$R_5$}
\put(71, 25){$R_6$}
\put(30, 85){$R_9$}
\put(30, 10){$R_{10}$}
\put(40, 27){$R_7$}

\put(53, 48){\scriptsize $g$}
\put(43, 48){\scriptsize $i$}
\put(33, 69){\scriptsize $f$}
\put(33, 37){\scriptsize $m$}
\end{overpic}
\caption{Regions $R_4 R_5 R_8 R_9 R_{11}$ and $R_5 R_6 R_7 R_9 R_{10}$ are both realized as immersed polygons when cut along $\alpha^\sigma_2$ and $\alpha^\sigma_1$, respectively.}
\label{regions458911}
\end{center}
\end{figure}
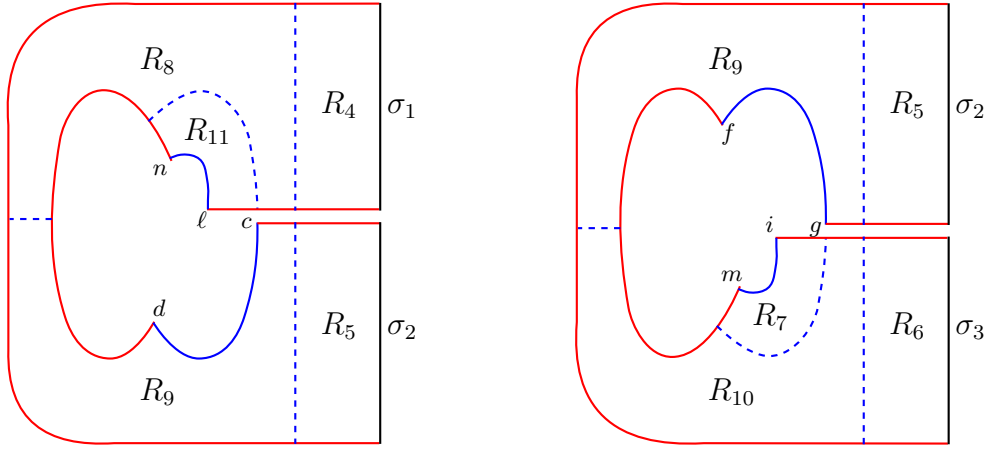

\begin{figure}[htbp]
\begin{center}
\includegraphics[scale = .8]{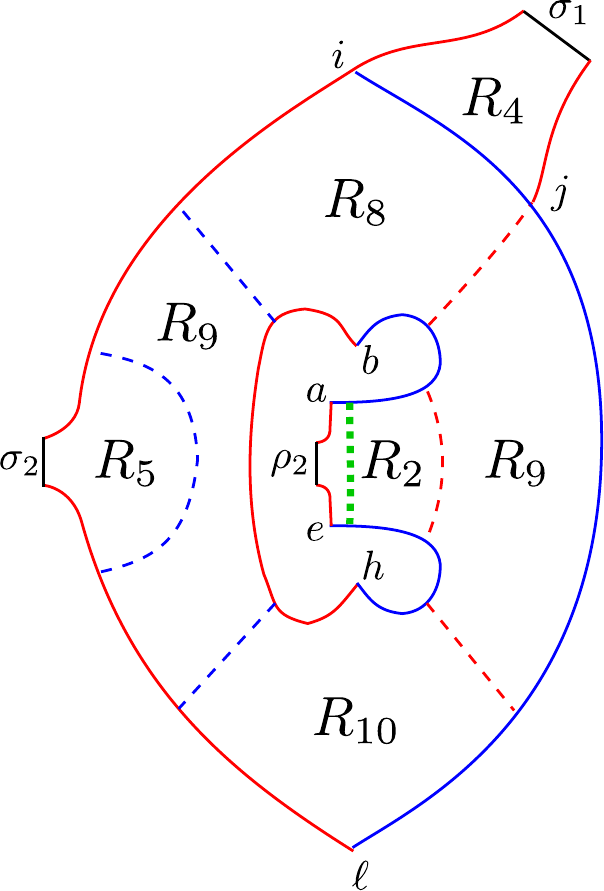}
\caption{The domain $R_2 R_4 R_5 R_8 R_9 R_9 R_{10}$. It is just like the annulus $A + B_2$ in Figure \ref{index0domain}, except that the inner boundary has too many $\alpha$ and $\beta$ segments. To fix this, we pinch along the green dotted arc in $R_2$.}
\label{regions_24589910}
\end{center}
\end{figure}

The domain $R_4 R_5 R_8 R_9 R_{11}$ is a polygon, though it may not be obvious at first glance. The only compatible sequence of Reeb chords is $(\sigma_1, \sigma_2)$. To realize the domain as an immersed surface with boundary Reeb chords $(\sigma_1, \sigma_2)$, we must cut along $\alpha^\sigma_2$, which produces a polygon (see Figure \ref{regions458911}). Similarly, the domain $R_5 R_6 R_7 R_9 R_{10}$ corresponds to a polygon with boundary $(\sigma_2, \sigma_3)$ after cutting along $\alpha^\sigma_1$. By Proposition \ref{polygons}, both of these domains contribute to the differential.

$R_1 R_2 R_3$ and $R_4 R_5 R_6$ can also be realized as polygons, and thus contribute with Reeb chords $(\rho_1,\rho_2, \rho_3)$ and $(\sigma_{123})$, respectively. Furthermore, $R_4 R_5 R_6$ can not contribute with its other sequence of compatible Reeb chords, $(\sigma_1, \sigma_2, \sigma_3)$, since cutting along the $\alpha^\sigma$ arcs from the boundary would produce a disconnected domain. Overall the domain $R_4 R_5 R_6$ contributes to the differential for each pair of generators it connects. The contribution of $R_1 R_2 R_3$ with $(\rho_{123})$ can be understood by examining the $\Ainfty$ relations for $(eb*, \rho_1, \rho_{23})$ and $(af*, \rho_2, \rho_3)$, where $*$ can be $i, j, k$, or $l$. The relations imply that 
\[ m(eb*, \rho_{123}) = m\left( m\big( m( eb*, \rho_1), \rho_2\big), \rho_3\right) = ah* . \]
Since the operation is nontrivial in $\CFAA$, $(R_1 R_2 R_3, (\rho_{123}))$ contributes to the differential from $eb*$ to $ah*$ in $\CFDD(\Hselfgluer)$ by Proposition \ref{Ainfty_trick}, and the total mod 2 contribution of $R_1 R_2 R_3$ is zero.

\bigbreak
{\bf Simple annuli:}
$R_8 R_9$ is an index zero annulus analogous to $A$ in Proposition \ref{Levine_prop}. The four domains obtained by adding $R_2$, $R_3$, $R_4$, or $R_5$ to this annulus may or may not contribute to $\CFDD(\Hselfgluer)$, depending on the choice of complex structure on $R_8 R_9$. Since we have chosen $J$ such that $\Theta_{8,9}^{d,b} < \Theta_{8,9}^{i,k}$, Proposition \ref{Levine_prop} asserts that none of these four domains contributes. Notice that since $R_2$ and $R_3$ are not bigons, we must first use Proposition \ref{pinching} to pinch off the extra $\alpha$ portion of the boundary, and then we can apply Proposition \ref{Levine_prop}. We specifically chose the complex structure $J$ to be consistent with pinching the appropriate arcs in $R_2$ and $R_3$.

Similarly, $R_9 R_{10}$ is an index zero annulus to which the regions $R_1$, $R_2$, $R_5$, or $R_6$ may be added. Given the choice that $\Theta_{9,10}^{h, f} < \Theta_{9,10}^{j, l}$, none of the four corresponding domains contribute.

Two more direct applications of Proposition \ref{Levine_prop} involve the annuli $R_2 R_3 R_8 R_9$ and $R_1 R_2 R_9 R_{10}$. Given that $\Theta_{2,3,8,9}^{a,a} < \Theta_{2,3,8,9}^{i, k}$, $R_2 R_3 R_7 R_8 R_9$ contributes to $\CFDD(\Hselfgluer)$, but $R_2 R_3 R_4 R_8 R_9$ and $R_2 R_3 R_5 R_8 R_9$ do not. The fact that $\Theta_{1,2,9,10}^{e,e} < \Theta_{1,2,9,10}^{j, l}$ implies that $R_1 R_2 R_9 R_{10} R_{11}$ contributes to $\CFDD(\Hselfgluer)$, but $R_1 R_2 R_5 R_9 R_{10}$ and $R_1 R_2 R_6 R_9 R_{10}$ do not. Note that for $R_2 R_3 R_7 R_8 R_9$ and $R_1 R_2 R_9 R_{10} R_{11}$ we make use of Proposition \ref{pinching} and the relevant assumptions about the complex structure on $R_7$ and $R_{11}$.

Finally, we will use Proposition \ref{Levine_prop} to account for the domain $R_2 R_4 R_5 R_8 R_9 R_9 R_{10}$, which connects $ahj$ to $ebl$. There is only one way to piece together these regions so that there are no unwanted corners, which is shown in Figure \ref{regions_24589910}. If we pinch $R_2$ along the arc connecting the two $\beta$ curves, then this domain has exactly the form of $D_2$ in Proposition \ref{Levine_prop}. Since we chose $J$ such that $\Theta_{2,5,8,9,9,10}^{h, b} < \Theta_{2,5,8,9,9,10}^{i, \ell}$, it follows that this domain does not count.

\bigbreak
{\bf More annuli:}
The domain $R_1 R_2 R_3 R_8 R_9$, with the Reeb chords $(\rho_1, \rho_2, \rho_3)$, is an annulus with one obtuse corner. If we pinch along the the arc through $R_1$, $R_2$, and $R_3$ parallel to $\alpha_0$, the modified annulus has one $\alpha$ and one $\beta$ segment on each boundary component. In this situation, we can apply the same reasoning as the proof of Proposition \ref{Levine_prop} (the only difference is that the cuts from the obtuse corner do not leave the annulus through to the opposite boundary component). Cutting along $\beta_1$ from $d$ makes the length of $\beta$ on the boundary component containing $d$ grow, so that as $c$ approaches $+\infty$, $\theta_{1,2,3,8,9}^{d, a}(c)$ approaches $-\infty$, and $\theta_{1,2,3,8,9}^{d, a}(c) - \theta_{1,2,3,8,9}^{i, k}(c)$ is negative. On the other hand, cutting along $\alpha_0$ from $d$ pinches off $R_1 R_2 R_3$ from the annulus $R_8 R_9$. In this extreme, $\theta_{1,2,3,8,9}^{d, a}(c) - \theta_{1,2,3,8,9}^{i, k}(c)$ approaches $\Theta_{8,9}^{d,b} - \Theta_{8,9}^{i, k} < 0$. Since the extremes are both negative, the mod 2 count of zeros, and thus the contribution of $R_1 R_2 R_3 R_8 R_9$ to $\CFDD(\Hselfgluer)$, is zero. An analogous argument shows that $R_1 R_2 R_3 R_9 R_{10}$ does not contribute with $(\rho_1, \rho_2, \rho_3)$.

$R_1 R_2 R_3 R_8 R_9$ may also contribute with the Reeb chord sequence $(\rho_{123})$. This contribution can be checked with the $\Ainfty$ trick, using the relation for $(edi, \rho_1, \rho_{23})$. This relation implies that $m(edi, \rho_{123}) = m \left( m( edi, \rho_1) , \rho_{23} \right)$. The only domain that could contribute a nontrivial operation $m(edi, \rho_1)$ is $R_3 R_8 R_9$. As discussed above, $R_3 R_8 R_9$ does not contribute for our choice of complex structure $J$, and therefore $R_1 R_2 R_3 R_8 R_9$ does not contribute with $\vecrho = (\rho_{123})$. An analogous argument also shows that $R_1 R_2 R_3 R_9 R_{10}$ does not contribute with $(\rho_{123})$.

The domain $R_4 R_5 R_6 R_8 R_9$ follows the same pattern. With Reeb chords $(\sigma_1, \sigma_2, \sigma_3)$ it is an annulus, and cuts in either direction split off the annulus $R_8 R_9$ or the annulus $R_4 R_5 R_8 R_9$. As we cut along $\beta_0$, $\theta_{4,5,6,8,9}^{l, k}(c) - \theta_{4,5,6,8,9}^{d, b}$ approaches $\Theta_{8,9}^{i, k} - \Theta_{8,9}^{d,b} > 0$. As we cut along $\alpha^\sigma_1$ toward the $\sigma$ boundary, $\theta_{4,5,6,8,9}^{l, k}(c)$ approaches $\infty$, and so $\theta_{4,5,6,8,9}^{l, k}(c) - \theta_{4,5,6,8,9}^{d, b}$ becomes positive. As a result, there is no contribution to the differential. The $\Ainfty$ relation for $(edl, \sigma_{12}, \sigma_3)$ reveals that $R_4 R_5 R_6 R_8 R_9$ also does not contribute with $(\sigma_{123})$. A similar argument shows that the domain $R_4 R_5 R_6 R_9 R_{10}$ does not contribute with either compatible sequence of Reeb chords.

\bigbreak
{\bf Corners: } Consider the domain $R_1 R_2 R_3 R_4 R_5 R_8 R_9$, which connects $edl$ to $ahl$. Any compatible sequence of Reeb chords must contain $(\sigma_1, \sigma_2)$, since $(\sigma_{12})$ would not be strongly boundary monotonic with respect to the $\sigma$ boundary. For the domain to have the chords $(\sigma_1, \sigma_2)$ along the $\sigma$ boundary, there must be a cut along $\alpha^\sigma_1$. However, such a cut would leave corners at the point $c$. Since neither the initial generator $edl$ nor the final generator $ahl$ contain $c$, it is impossible to have a corner at $c$. As a result, the domain $R_1 R_2 R_3 R_4 R_5 R_8 R_9$ can not contribute to $\CFDD(\Hselfgluer)$. The same reasoning applies to the domains $R_2 R_4 R_5 R_8 R_9$ and $R_3 R_4 R_5 R_8 R_9$.

Similarly, $R_1 R_2 R_3 R_5 R_6 R_9 R_{10}$ is only compatible with Reeb chord sequences containing $(\sigma_2, \sigma_3)$. This Reeb chord sequence requires a cut along $\alpha^\sigma_1$, which leaves corners at the point $g$. Since the initial generator $ebi$ and the final generator $afi$ do not contain $g$, this domain can not contribute to the differential. The same is true for the domains $R_1 R_5 R_6 R_9 R_{10}$ and $R_2 R_5 R_6 R_9 R_{10}$, so these also do not contribute.

The domain $R_2 R_5 R_8 R_9 R_9 R_{10} R_{11}$ connects the generators $ahi$ and $enc$. However, there is no way to piece together these seven regions without having corners at points other than $a$, $h$, $i$, $e$, $n$, and $c$. Therefore, this domain can not contribute to the differential.

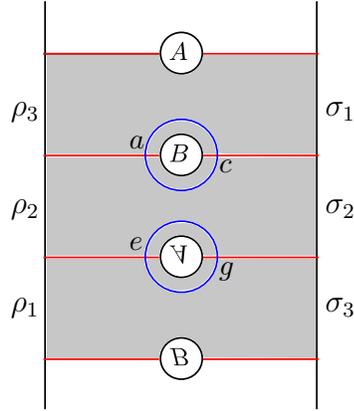
\begin{figure}[htbp]
\begin{center}
\begin{overpic}[scale = .8]{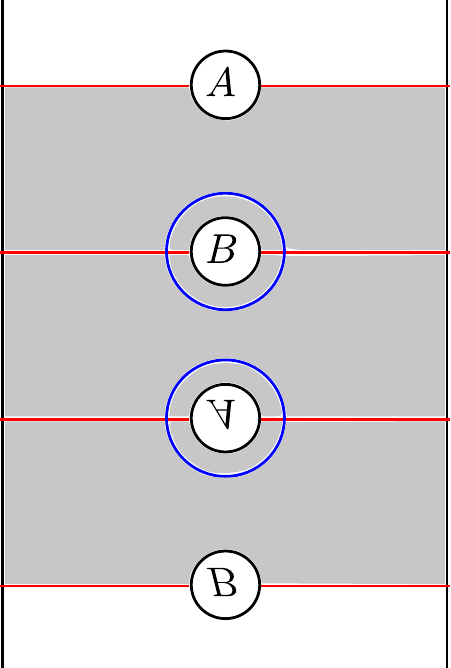}
\put(-8, 24){$\rho_1$}
\put(-8, 48){$\rho_2$}
\put(-8, 72){$\rho_3$}

\put(69, 72){$\sigma_1$}
\put(69, 48){$\sigma_2$}
\put(69, 24){$\sigma_3$}

\put(21, 64){\small $a$}
\put(43, 58){\small $c$}
\put(21, 39){\small $e$}
\put(43, 33){\small $g$}

\end{overpic}
\caption{A bordered Heegaard diagram for the identity bimodule, $\CFDD(\mathbb{I})$. $R_1 R_2 R_3 R_4 R_5 R_6 R_8 R_9 R_{10}$ contributes to $\CFDD(\Hselfgluer)$ if and only if the shaded region contributes to the differential from $ec$ to $ag$ in $\CFDD(\mathbb{I})$.}
\label{identity_diagram}
\end{center}
\end{figure}

\bigbreak
{\bf $ \bf R_1 R_2 R_3 R_4 R_5 R_6 R_8 R_9 R_{10}$:}
This domain has four compatible sequences of Reeb chords. It is possible to use $\Ainfty$ relations and analyze the contribution of each one. However, it is easier to notice that this domain contributes if an only if the shaded domain contributes in the Heegaard diagram for the mapping cylinder of the identity map shown in Figure \ref{identity_diagram}. The computation of $\CFDD(\mathbb{I})$ in \cite[Proposition 10.1]{LOT:Bimodules} reveals that this domain must contribute.

\bigbreak
{\bf Using $\bf \partial^2 = 0$:}
We can deduce the contribution of other domains using the fact that $\CFDD(\Hselfgluer)$ must satisfy $\partial^2 = 0$. At this point, we have enough information to deduce the contribution of all domains but one using $\partial^2$. Table \ref{differential_so_far} gives the differential on $\CFDD(\Hselfgluer)$ as computed so far, with coefficients $\lambda_{i_1, \ldots, i_k}$ representing unknown contributions. $\lambda_{i_1, \ldots, i_k}$ is 1 if the domain $R_{i_1} \cdots R_{i_k}$ contributes to the differential, and 0 otherwise.

\begin{table}
\begin{framed}
\begin{align*}
\partial(afi) &= \rho_2 (edi) \\
\partial(afj) &= \sigma_1 (afi) + \rho_2 (edj) \\
\partial(afk) &=\sigma_2 (afj) + \rho_2 (edk) \\
\partial(afl) &= \sigma_{123} (afi) + \sigma_3 (afk) + \sigma_1 ang + \rho_2 edl \\
\partial(ahi) &= \lambda_{2,3,5,8,9,9,10} \rho_{23}\sigma_2 (afl) +  \lambda_{2,5,6,8,9,9,10} \rho_2 \sigma_{23} (ebk) \\
\partial(ahj) &= \sigma_1(ahi) + \lambda_{4,8,9,10,11} \sigma_1 (ang) \\
\partial(ahk) &= \sigma_2(ahj) \\
\partial(ahl) &= \sigma_{123}(ahi) + \sigma_3(ahk) + \lambda_{2,3,4,5,6,8,9}\rho_{23}\sigma_{123}(ahk) + \lambda_{4,5,6,8,9,10,11} \sigma_{123}(ang) \\
\partial(ang) &= \rho_2 \sigma_2(enc) \\
\partial(amg) &= \sigma_{23}(afi) + (ahi) + \rho_{23}(ahk) + \lambda_{7,8,9,10,11}(ang) + (ang)\\
 &  \quad + \lambda_{2,5,7,8,9,9,10}\rho_2\sigma_2(ebl) + \rho_2\sigma_2(emc) \\
\partial(ebi) &= \rho_3(afi) + \lambda_{1,2,3,8,9,10,11} \rho_{123}(ang) \\
\partial(ebj) &= \rho_3(afj) +\sigma_1(ebi) + \lambda_{1,2,4,5,6,9,10} \rho_{12}\sigma_{123}(ebi) + \rho_{12}(enc) \\
\partial(ebk) &= \rho_3(afk) + \rho_1(ang) + \sigma_2(ebj) \\
\partial(ebl) &= \rho_3(afl) + \sigma_{123}(ebi) + \sigma_3(ebk) + (enc) \\
\partial(edi) &= \rho_1(ahi) + \lambda_{1,8,9,10,11} \rho_1(ang) + \lambda_{1,2,5,8,9,9,10} \rho_{12}\sigma_2(ebl) \\
\partial(edj) &= \rho_1(ahj) + \sigma_1(edi) \\
\partial(edk) &= \rho_1(ahk) + \sigma_2(edj) \\
\partial(edl) &= \rho_1(ahl) + \sigma_{123}(edi) + \sigma_3(edk) + \sigma_{12}(enc) \\
\partial(enc) &= \rho_3\sigma_1(ang) + \rho_1\sigma_3(ang) + \rho_{123}\sigma_{123}(ang) \\
\partial(emc) &= \lambda_{3,7,8,9,10}\rho_3(afl) + \rho_3(ahj) + \lambda_{1,2,3,7,8,9,10}\rho_{123}(ahl) \\
& \quad  + \rho_3\sigma_1(amg) + \rho_1\sigma_3(amg) + \rho_{123}\sigma_{123}(amg) + \lambda_{4,5,6,7,8,9,10}\sigma_{123}(ebi)\\
& \quad  + \lambda_{6,7,8,9,10} \sigma_3(ebk) + \sigma_3(edi) + \lambda_{7,8,9,10,11}(enc) + (enc)
\end{align*}
\end{framed}
\caption{The differential on $\CFDD(\Hselfgluer)$. $\lambda$ is used for coefficients that have yet to be determined; they are 0 or 1 depending on the contribution of the corresponding domain.}
\label{differential_so_far}
\end{table}

Consider the generator $ahj$. $\partial(ahj) = \sigma_1(ahi) + \lambda_{4,8,9,10,11} \sigma_1 (ang)$, and so
\begin{eqnarray*}
 0 = \partial^2(ahj) &=& \sigma_1 \big( \lambda_{2,3,5,8,9,9,10} \rho_{23}\sigma_2 (afl) +  \lambda_{2,5,6,8,9,9,10} \rho_2 \sigma_{23} (ebk) \big)  \\
& & +  \lambda_{4,8,9,10,11} \sigma_1 \big(\rho_2 \sigma_2(enc) \big).
\end{eqnarray*}
It follows that $\lambda_{2,3,5,8,9,9,10} = \lambda_{2,5,6,8,9,9,10} = \lambda_{4,8,9,10,11} = 0$. Given these values, we find that
\begin{eqnarray*}
0 = \partial^2(edi) &=& \rho_1\big(0 \big) + \lambda_{1,8,9,10,11} \rho_1\big(  \rho_2 \sigma_2(enc) \big) \\
& & + \lambda_{1,2,5,8,9,9,10} \rho_{12}\sigma_2 \bigg( \rho_3(afl) + \sigma_{123}(ebi) + \sigma_3(ebk) + (enc) \bigg) \\
&=& \big( \lambda_{1,8,9,10,11} + \lambda_{1,2,5,8,9,9,10} \big) \rho_{12}\sigma_2(enc) \\
& & + \lambda_{1,2,5,8,9,9,10}  \rho_{123}\sigma_2(afl) \\
& & + \lambda_{1,2,5,8,9,9,10} \rho_{12}\sigma_{23}(ebk).
\end{eqnarray*}
The coefficient of $afl$ implies that $\lambda_{1,2,5,8,9,9,10} = 0$, and the coefficient of $enc$ implies that $\lambda_{1,8,9,10,11} = 0$.

The coefficient of the $enc$ term of $\partial^2(amg)$ is $\lambda_{2,5,7,8,9,9,10} \rho_2 \sigma_2$, which implies that $\lambda_{2,5,7,8,9,9,10} = 0$. Then the $afl$ term of $\partial^2(amg)$ becomes $\lambda_{3,7,8,9,10} \rho_{23}\sigma_2 (afl)$, and the $ebk$ term becomes $\lambda_{6,7,8,9,10} \rho_2 \sigma_{23} (ebk)$, implying that $\lambda_{3,7,8,9,10} = \lambda_{6,7,8,9,10} = 0$. Similarly the $ang$ term of $\partial^2(edl)$ reveals that $\lambda_{4,5,6,8,9,10,11} = 1$ and the $ahk$ term of $\partial^2(edl)$ implies that $\lambda_{2,3,4,5,6,8,9} = 0$. The $afi$ term of $\partial^2(ebj)$ implies that $\lambda_{1,2,4,5,6,9,10} = 0$, and the $ang$ term implies that $\lambda_{1,2,3,8,9,10,11} = 1$. Finally, the $ang$ term of $\partial^2(emc)$ implies that $\lambda_{4,5,6,7,8,9,10} = 1$, and the $ahi$ term implies that $\lambda_{1,2,3,7,8,9,10} = 1$. The only coefficient in Table \ref{differential_so_far} that remains undetermined is $\lambda_{7,8,9,10,11}$.

\bigbreak
{$\bf R_7 R_8 R_9 R_{10} R_{11}$:}
We have determined that $\CFDD(\Hselfgluer)$ is one of two possibilities, depending on the value of $\lambda_{7,8,9,10,11}$. We will deduce the right choice by showing that one of these possible bimodules does not behave correctly under tensoring with type $A$ modules for the solid torus.

\begin{figure}[htbp]
\begin{center}
\begin{overpic}[scale = .6]{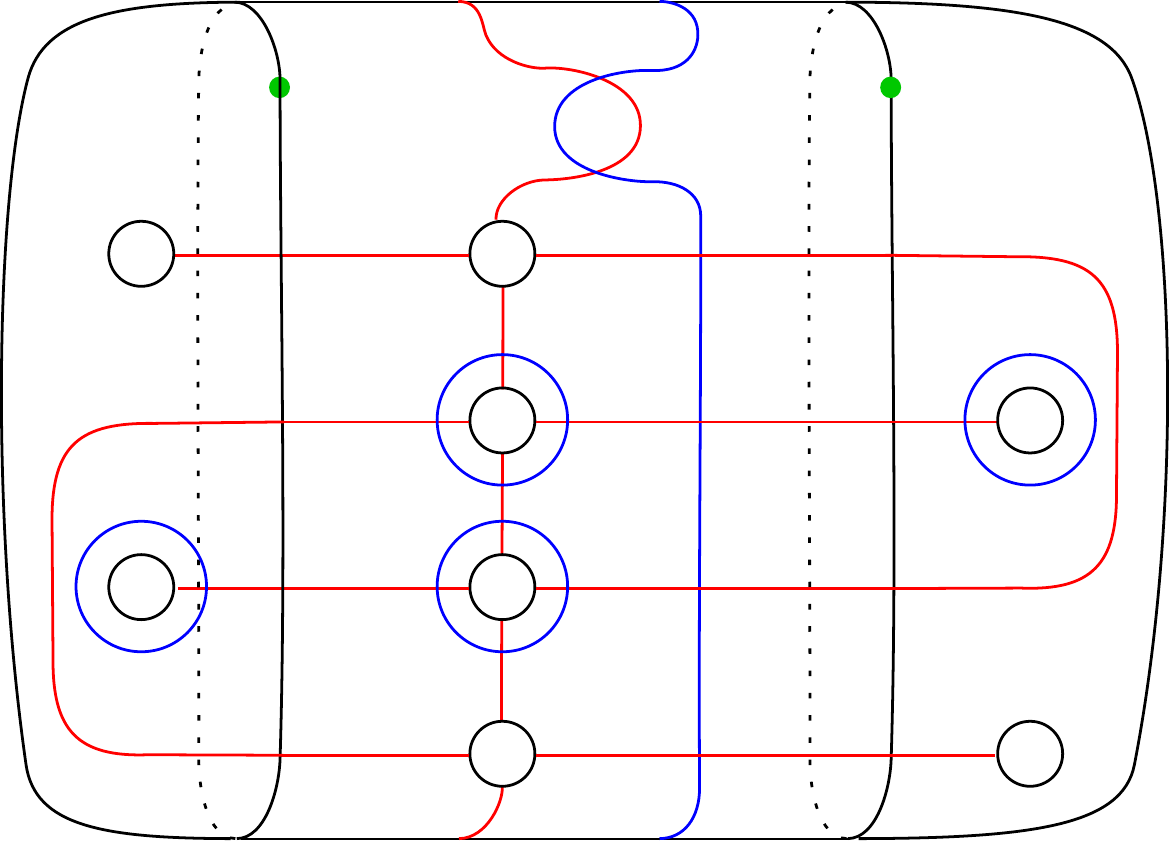}
\end{overpic}
\qquad \qquad
\begin{overpic}[scale = .7]{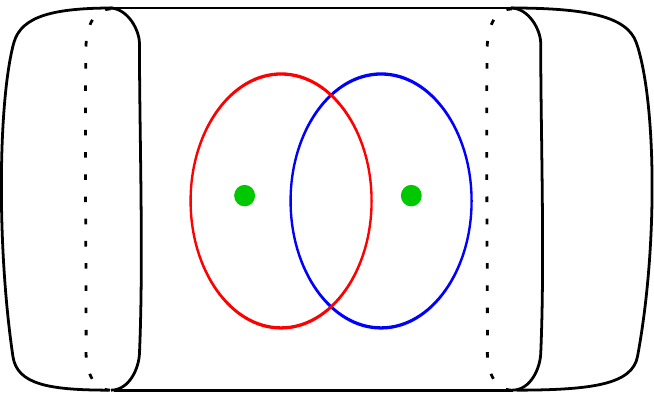}
\end{overpic}

\qquad \quad (a) \hspace{7 cm} (b)

\caption{The manifold obtained from the self-gluer by capping off both ends with identical 0-framed solid tori.}
\label{selfgluer_capped}
\end{center}
\end{figure}

Consider the closed, doubly basepointed Heegaard diagram in Figure \ref{selfgluer_capped}(a), which is obtained from the bordered Heegaard diagram $\Hselfgluer$ by gluing bordered Heegaard diagrams for solid tori to each boundary component. A sequence of isotopies and destabilizations leads to the diagram in Figure \ref{selfgluer_capped}{(b)}, so it is easy to check that $\HFhat$ of the manifold represented by this diagram has rank 2. $\HFhat$ can also be obtained by taking the box tensor product of $\CFDD(\Hselfgluer)$ with two copies of the type $A$ module for the solid torus. A bounded version of the solid torus module has three generators $x$, $y$, and $z$ and the following operations:
\[
m_1(x) = z, \quad m_2(x, \rho_1) = y, \quad m_2(x, \rho_2) = z, \quad m_2(x, \rho_{12}) = z.
\]
It is routine to perform the tensor products, and we find that the homology of the resulting chain complex is rank 2 if $\lambda_{7,8,9,10,11} = 1$, and rank 4 if $\lambda_{7,8,9,10,11} = 0$. Thus, the domain $R_7 R_8 R_9 R_{10} R_{11}$ must contribute, and we have completed the computation of $\CFDD(\Hselfgluer)$ in the middle spin$^c$-structure. The result is pictured in Figure \ref{selfgluer_answer}.

\begin{figure}[ps]

\begin{sideways}
\begin{tikzpicture}[>=latex, scale=.8]

\node (ebi)  at (0,0) {$ebi$};
\node (ebj)  at (0,3) {$ebj$};
\node (ebk)  at (0,6) {$ebk$};
\node (ebl)  at (0,9) {$ebl$};

\node (afi)  at (4,0) {$afi$};
\node (afj)  at (4,3) {$afj$};
\node (afk)  at (4,6) {$afk$};
\node (afl)  at (4,9) {$afl$};

\node (edi)  at (8,0) {$edi$};
\node (edj)  at (8,3) {$edj$};
\node (edk)  at (8,6) {$edk$};
\node (edl)  at (8,9) {$edl$};

\node (ahi)  at (12,0) {$ahi$};
\node (ahj)  at (12,3) {$ahj$};
\node (ahk)  at (12,6) {$ahk$};
\node (ahl)  at (12,9) {$ahl$};

\node (enc) at (0, 12) {$enc$};
\node (ang) at (-4, 9) {$ang$};
\node (amg) at (12, -3) {$amg$};
\node (emc) at (16, 0) {$emc$};

\draw[->] (ebi) to node[above = -2pt] {$\rho_3$} (afi);
\draw[->] (ebj) to node[above = -2pt] {$\rho_3$} (afj);
\draw[->] (ebk) to node[above = -2pt] {$\rho_3$} (afk);
\draw[->] (ebl) to node[above = -2pt] {$\rho_3$} (afl);

\draw[->] (afi) to node[above = -2pt] {$\rho_2$} (edi);
\draw[->] (afj) to node[above = -2pt] {$\rho_2$} (edj);
\draw[->] (afk) to node[above = -2pt] {$\rho_2$} (edk);
\draw[->] (afl) to node[above = -2pt] {$\rho_2$} (edl);

\draw[->] (edi) to node[above = -2pt] {$\rho_1$} (ahi);
\draw[->] (edj) to node[above = -2pt] {$\rho_1$} (ahj);
\draw[->] (edk) to node[above = -2pt] {$\rho_1$} (ahk);
\draw[->] (edl) to node[above = -2pt] {$\rho_1$} (ahl);

\draw[->] (ebl) to node[right = -2pt] {$\sigma_3$} (ebk);
\draw[->] (afl) to node[right = -2pt] {$\sigma_3$} (afk);
\draw[->] (edl) to node[right = -2pt] {$\sigma_3$} (edk);
\draw[->] (ahl) to node[right = -2pt] {$\sigma_3$} (ahk);

\draw[->] (ebk) to node[right = -2pt] {$\sigma_2$} (ebj);
\draw[->] (afk) to node[right = -2pt] {$\sigma_2$} (afj);
\draw[->] (edk) to node[right = -2pt] {$\sigma_2$} (edj);
\draw[->] (ahk) to node[right = -2pt] {$\sigma_2$} (ahj);

\draw[->] (ebj) to node[right = -2pt] {$\sigma_1$} (ebi);
\draw[->] (afj) to node[right = -2pt] {$\sigma_1$} (afi);
\draw[->] (edj) to node[right = -2pt] {$\sigma_1$} (edi);
\draw[->] (ahj) to node[right = -2pt] {$\sigma_1$} (ahi);

\draw[->, bend right = 20, pos = .6] (ebl) to node[left = -2pt] {$\sigma_{123}$} (ebi);
\draw[->, bend right = 20] (afl) to node[left = -2pt] {$\sigma_{123}$} (afi);
\draw[->, bend right = 20] (edl) to node[left = -2pt] {$\sigma_{123}$} (edi);
\draw[->, bend right = 20] (ahl) to node[left = -2pt] {$\sigma_{123}$} (ahi);

\draw[->] (ebl) to (enc);
\draw[->, bend left = 30] (ebj) to node[left = -2pt] {$\rho_{12}$} (enc);
\draw[->] (edl) to [out = 90, in = 0] node[right = 15pt] {$\sigma_{12}$} (enc);

\draw[->] (ebi) to [out = 180, in = -90] node[left = -2pt] {$\rho_{123}$} (ang);
\draw[->, bend left = 30] (ebk) to node[left = 2pt] {$\rho_1$} (ang);
\draw[->, bend right = 20, pos = .3] (afl) to node[above] {$\sigma_1$} (ang);
\draw[->] (ahl) to [out = 150, in = 30] node[below] {$\sigma_{123}$} (ang);

\draw[->] (amg) to (ahi);
\draw[->, bend right = 30] (amg) to node[right = -2pt] {$\rho_{23}$} (ahk);
\draw[->] (amg) to [out = 180, in = -90] node[left = 15pt] {$\sigma_{23}$} (afi);

\draw[->] (emc) to [out = 90, in = 0] node[right = -2pt] {$\rho_{123}$} (ahl);
\draw[->, bend right = 30] (emc) to node[right = 2pt] {$\rho_3$} (ahj);
\draw[->, bend left = 20, pos = .7] (emc) to node[below] {$\sigma_3$} (edi);
\draw[->] (emc) to [out = -150, in = -30] node[above] {$\sigma_{123}$} (ebi);

\draw[->] ([xshift = 5pt] ang.north) to [out = 90, in = 180] node[below right = -3pt] {$\rho_2 \sigma_2$} ([yshift = -5pt] enc.west);
\draw[->] ([yshift = 5pt] enc.west) to [out = 180, in = 90] node[above = 10pt] {$\begin{array}{c} \rho_1 \sigma_3 + \\  \rho_3 \sigma_1 + \\ \rho_{123} \sigma_{123} \end{array}$} ([xshift = -5pt] ang.north);

\draw[->] ([yshift = -5pt] amg.east) to [out = 0, in = -90] node[above left = 3pt] {$\rho_2 \sigma_2$} ([xshift = 5pt] emc.south);
\draw[->] ([xshift = -5pt] emc.south) to [out = -90, in = 0] node[below= 15pt] {$\begin{array}{c} \rho_1 \sigma_3 + \\  \rho_3 \sigma_1 + \\ \rho_{123} \sigma_{123} \end{array}$} ([yshift = 5pt] amg.east);

\end{tikzpicture}
\end{sideways}

\caption{$\CFDD(\Yselfgluer)$ in the middle spin$^c$-structure.}
\label{selfgluer_answer}
\end{figure}

\bigbreak

\noindent {\bf Extremal spin$^c$-structures:}
First consider the spin$^c$-structure in which both $\alpha^\rho$ arcs are occupied and neither $\alpha^\sigma$ arc is occupied. There are only two generators with those conditions: $aem$ and $aen$. There are two domains which have the right corner count to connect $aem$ and $aen$. The bigon $R_{12}$ contributes a differential from $aem$ to $aen$. The domain $R_7 R_8 R_9 R_{10} R_{11}$, as an element of $\pi_2(aem, aen)$, has index $-1$ and thus does not contribute. Canceling the differential and two generators, we find that $\CFDD(\Yselfgluer)$ in this extremal spin$^c$ structure has no generators.

The other extremal spin$^c$-structure has more generators (given the choice of Heegaard diagram $\Hselfgluer$), but the corresponding summand of $\CFDD(\Yselfgluer)$ is still trivial. Indeed, we could handleslide $\beta_0$ across the handles in Figure \ref{self_gluer} to produce a new Heegaard diagram with only two generators in this spin$^c$ structure. This diagram is a mirror image of $\Hselfgluer$, so the reasoning above applies and shows that the two generators are canceled by the single differential between them.

\subsection{Gradings}

As described in Section \ref{sec:gradings}, $\CFDD(\Yselfgluer)$ is graded by a set which is a quotient of the noncommutative group $G_{2,0}$. We will compute this (relative) grading using Figure \ref{selfgluer_answer}.

We choose $ebl$ to be the preferred generator and set $gr(ebl) = \vec0 = (0; 0, 0; 0, 0)$. The arrow labeled $\rho_3$ from $ebl$ to $afl$ determines the grading of $afl$.
\[ gr(afl) = \lambda^{-1} gr(\rho_3)^{-1} gr(ebl) = \lambda^{-1} \left( \frac{1}{2}; \frac{1}{2}, -\frac{1}{2}; 0, 0\right) \vec0 = \left( -\frac{1}{2}; \frac{1}{2}, -\frac{1}{2}; 0, 0\right)  \]
Similarly, the successive arrows labeled $\rho_2$ and $\rho_1$ (moving right from $afl$ in Figure \ref{selfgluer_answer}) determine the gradings of $edl$ and $ahl$.
\begin{eqnarray*}
 gr(edl) &=& \lambda^{-1} gr(\rho_2)^{-1} gr(afl) = \lambda^{-1} \left( \frac{1}{2}; -\frac{1}{2}, -\frac{1}{2}; 0, 0\right)  \left( -\frac{1}{2}; \frac{1}{2}, -\frac{1}{2}; 0, 0\right) \\
 &=&  \left( -\frac{1}{2}; 0, -1; 0, 0\right) \\
 gr(ahl) &=& \lambda^{-1} gr(\rho_1)^{-1} gr(edl) = \lambda^{-1} \left( \frac{1}{2}; -\frac{1}{2}, \frac{1}{2}; 0, 0\right) \left( -\frac{3}{2}; 0, -1; 0, 0\right) \\
&=&  \left( -\frac{1}{2}; -\frac{1}{2}, -\frac{1}{2}; 0, 0\right)
\end{eqnarray*}
Working down the chain of $\sigma$ labeled arrows starting from $ebl$ determines the gradings of $ebk$, $ebj$, and $ebi$.
\begin{eqnarray*}
gr(ebk) &=& \lambda^{-1} gr(\sigma_3)^{-1} gr(ebl) = \left( -\frac{1}{2}; 0, 0; \frac{1}{2}, -\frac{1}{2}\right) \\
gr(ebj) &=& \lambda^{-1} gr(\sigma_2)^{-1} gr(ebk) = \left( -\frac{1}{2}; 0, 0; 0, -1 \right) \\
gr(ebi) &=& \lambda^{-1} gr(\sigma_1)^{-1} gr(ebj) = \left( -\frac{1}{2}; 0, 0; -\frac{1}{2}, -\frac{1}{2}\right)
\end{eqnarray*}
The vertical chains of $\sigma$ labeled arrows from $afl$, $edl$, and $ahl$ determine the following gradings:
\begin{eqnarray*}
gr(afk) &=& \lambda^{-1} gr(\sigma_3)^{-1} gr(afl) = \left( -1; \frac{1}{2}, -\frac{1}{2}; \frac{1}{2}, -\frac{1}{2}\right) \\
gr(afj) &=& \lambda^{-1} gr(\sigma_2)^{-1} gr(afk) = \left( -1; \frac{1}{2}, -\frac{1}{2}; 0, -1 \right) \\
gr(afi) &=& \lambda^{-1} gr(\sigma_1)^{-1} gr(afj) = \left( -1; \frac{1}{2}, -\frac{1}{2}; -\frac{1}{2}, -\frac{1}{2}\right) \\
gr(edk) &=& \lambda^{-1} gr(\sigma_3)^{-1} gr(edl) = \left( -1; 0, -1; \frac{1}{2}, -\frac{1}{2}\right) \\
gr(edj) &=& \lambda^{-1} gr(\sigma_2)^{-1} gr(edk) = \left( -1; 0, -1; 0, -1 \right) \\
gr(edi) &=& \lambda^{-1} gr(\sigma_1)^{-1} gr(edj) = \left( -1; 0, -1; -\frac{1}{2}, -\frac{1}{2}\right) \\
gr(ahk) &=& \lambda^{-1} gr(\sigma_3)^{-1} gr(ahl) = \left( -1; -\frac{1}{2}, -\frac{1}{2}; \frac{1}{2}, -\frac{1}{2}\right) \\
gr(ahj) &=& \lambda^{-1} gr(\sigma_2)^{-1} gr(ahk) = \left( -1; -\frac{1}{2}, -\frac{1}{2}; 0, -1 \right) \\
gr(ahi) &=& \lambda^{-1} gr(\sigma_1)^{-1} gr(ahj) = \left( -1; -\frac{1}{2}, -\frac{1}{2}; -\frac{1}{2}, -\frac{1}{2}\right)
\end{eqnarray*}
The two unlabeled arrows in the diagram determine the gradings so $enc$ and $amg$.
\begin{eqnarray*}
gr(enc) &=& \lambda^{-1}  gr(ebl) = (-1; 0, 0; 0, 0) \\
gr(amg) &=& \lambda  gr(ahi) = \left( 0; -\frac{1}{2}, -\frac{1}{2}; -\frac{1}{2}, -\frac{1}{2}\right) \\
\end{eqnarray*}
Finally, the two arrows labeled $\rho_2 \sigma_2$ determine the gradings of $ang$ and $emc$.
\begin{eqnarray*}
gr(ang) &=& \lambda gr(\rho_2) gr(\sigma_2)  gr(enc) =  \left( -1; \frac{1}{2}, \frac{1}{2}; \frac{1}{2}, \frac{1}{2}\right) \\
gr(emc) &=& \lambda gr(\rho_2) gr(\sigma_2)  gr(amg) =  (0; 0, 0; 0, 0)
\end{eqnarray*}

It remains to compute the indeterminacy $\mathcal{P}(ebl)$. We compute equivalent values for the grading of $ebl$ by using the loop $ebl$ to $afl$ to $edl$ to $enc$ to $ebl$ and the loop $ebl$ to $ebk$ to $ebj$ to $enc$ to $ebl$. The first loop gives the element of $\mathcal{P}(ebl)$ corresponding to a periodic domain with boundary $\rho_{23} + \sigma_{12}$.
\begin{eqnarray*}
gr(afl) &=& \lambda^{-1} gr(\rho_3)^{-1} gr(ebl) = \left( -\frac{1}{2}; \frac{1}{2}, -\frac{1}{2}; 0, 0\right) \\
gr(edl) &=& \lambda^{-1} gr(\rho_2)^{-1} gr(afl) = \left( -\frac{1}{2}; 0, -1; 0, 0\right) \\
gr(enc) &=& \lambda^{-1} gr(\sigma_{12})^{-1} gr(afl) = \left( -1; 0, -1; -1, 0\right) \\
gr(edl) &=& \lambda gr(enc) = \left( 0; 0, -1; -1, 0\right)
\end{eqnarray*}
The second loop gives the element of $\mathcal{P}(ebl)$ corresponding to a periodic domain with boundary $\rho_{12} + \sigma_{23}$.
\begin{eqnarray*}
gr(ebk) &=& \lambda^{-1} gr(\sigma_3)^{-1} gr(ebl) = \left( -\frac{1}{2}; 0, 0; \frac{1}{2}, -\frac{1}{2} \right) \\
gr(ebj) &=& \lambda^{-1} gr(\sigma_2)^{-1} gr(ebk) = \left( -\frac{1}{2}; 0, 0; 0, -1 \right) \\
gr(enc) &=& \lambda^{-1} gr(\rho_{12})^{-1} gr(ebj) = \left( -1; -1, 0; 0, -1 \right) \\
gr(edl) &=& \lambda gr(enc) = \left( 0; -1, 0; 0, -1 \right)
\end{eqnarray*}
Thus $\mathcal{P}(ebl)$ is the subgroup of $G_{2,0}$ generated by $(0; 0, -1; -1, 0)$ and $(0; -1, 0; 0, -1)$.

\section{Computing $\HFhat$ of graph manifolds}
\label{sec:HFhat}

This section describes the procedure for computing $\HFhat$ of an arbitrary graph manifold given a connected plumbing graph $\Gamma$. For simplicity, we will assume that every vertex of $\Gamma$ has nonnegative genus. The manifold can be constructed from simpler bordered pieces using two types of gluing: \emph{extension} glues fibers to fibers and base surface to base surface, and \emph{plumbing} glues a fiber of one bundle to a curve in the base of the other bundle. Gluing two $S^1$-bundles by extension produces an $S^1$-bundle over the surface obtained by gluing the two bases.

Recall that in the Heegaard Diagram for $\Y$, $\alpha^\rho_1$, $\alpha^\sigma_2$, and $\alpha^\tau_1$ parametrize curves in the base surface $\pants$, while $\alpha^\rho_2$, $\alpha^\sigma_1$, and $\alpha^\tau_2$ parametrize fibers. If we glue two type $D$ boundaries together, $\alpha_1$ glues to $\alpha_2$ and vice versa (to combine the relevant modules we would first change one of the boundaries to type $A$, which switches $\alpha_1$ and $\alpha_2$). Thus gluing the $\rho$ boundary of one copy of $\Y$ to the $\sigma$ boundary of another is extension. Gluing the $\rho$ boundary to the $\tau$ boundary is plumbing. 

It will be convenient to introduce the bordered manifold $\Ybar$, the mirror image of $\Y$. The trimodule $\CFD{3}(\Ybar)$ can be obtained from $\CFD{3}(\Y)$ by interchanging 1's with 3's for all algebra elements and reversing the direction of the arrows. $\alpha_1$ and $\alpha_2$ are also interchanged on each boundary component.

\subsection{Trivial bundles over surfaces}

Recall that each vertex of $\Gamma$ represents a particular $S^1$-bundle over a surface $S_{g, b}$ with genus $g$ and $b$ boundary components. We first construct the trivial bundle over $S_{g, b}$.

If $g = 0$ and $b \ge 3$, then we simply glue copies of $\Y$ by extension until we have the right number of boundary components. The multimodule $\CFD{b}$ is obtained by taking box tensor products, inserting copies of $\CFAAid$ when two type $D$ boundaries are glued. For intance, $\CFD{4}(S^1 \times S_{0,4} )$ is given by
\[ \big( \CFAAid \boxtimes \CFD{3}(\Y) \big) \boxtimes \CFD{3}(\Y), \]
where the tensor products are with respect to the $\rho$ and $\sigma$ boundaries on the two copies of $\CFD{3}(\Y)$. The trivial bundle over $S_{0, 1}$ is just the solid torus, which has bordered invariant
\begin{center}
\begin{overpic}{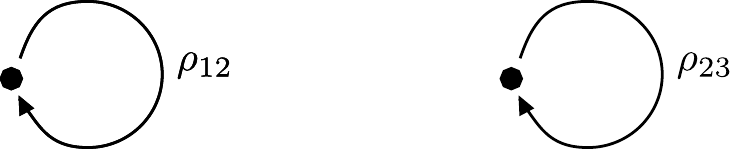}
\put(45, 10){ or }
\end{overpic}
\end{center}
depending on whether $\alpha_1$ parametrizes a curve in the base (left) or a fiber (right) \cite[Section 11.2]{LOT:Bordered}. The trivial bundle over $S_{0,2}$ is the same as the mapping cylinder of the identity map on the torus. The corresponding bimodule $\CFDD( \mathbb{I} )$ is computed in \cite[Proposition 10.1]{LOT:Bimodules}. Here either $\alpha$ arc can be the fiber, but $\alpha_1$ on one boundary is the same as $\alpha_2$ on the other boundary.

\begin{figure}[htbp]
\begin{center}
\begin{overpic}[scale = .7]{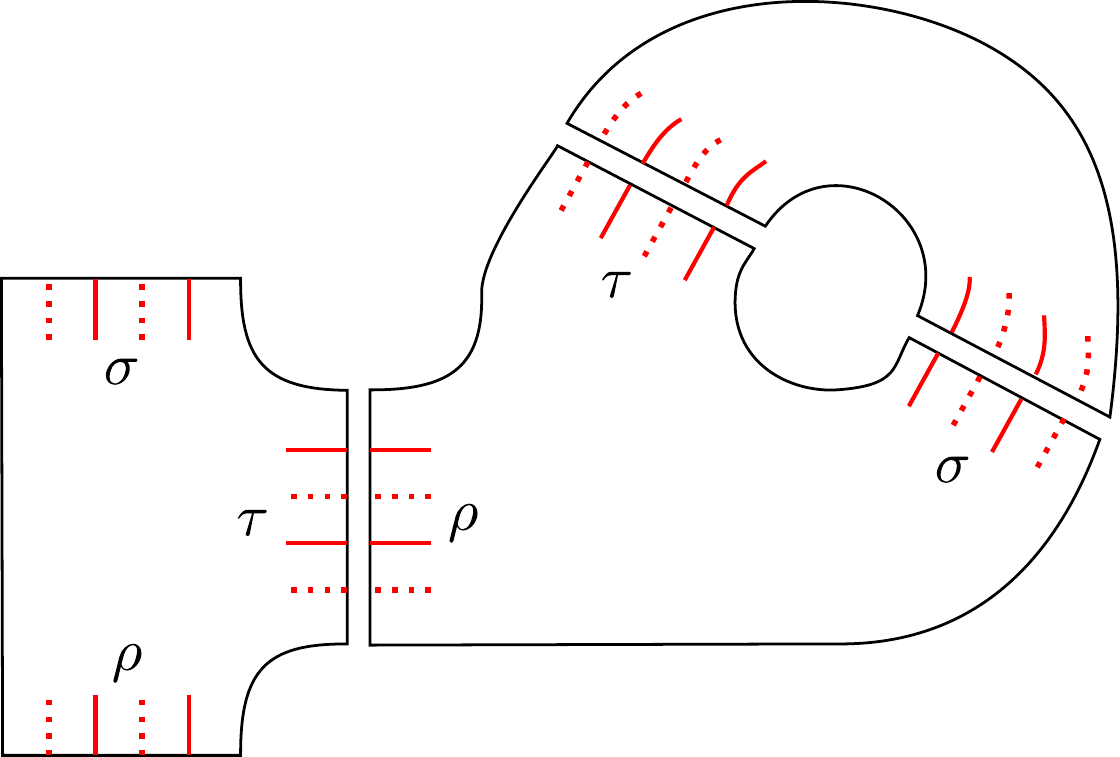}

\put(2,20){ $\HeegDiag(\Y)$ }
\put(55,20){ $\HeegDiag(\Ybar)$ }
\put(74,54){ $\Hselfgluer$ }

\end{overpic}
\caption{A Heegaard diagram for $S_{1,2} \times S^1$ can be obtained by gluing $\Hselfgluer$ and Heegaard diagrams for $\Y$ and $\Ybar$. Solid red lines indicate $\alpha$ arcs in the base of each bundle, while dotted red lines correspond to fibers.}
\label{plus_genus}
\end{center}
\end{figure}

We construct a Heegaard diagram for the trivial bundle over $S_{1,2}$ as indicated in Figure \ref{plus_genus}. Notice that we must insert a copy of $\Hselfgluer$ when we glue two components of $\partial \Ybar$ to each other. The bimodule $\CFDD(S^1 \times S_{1,2})$ can be computed as follows:
\begin{itemize}
\item
Change $\CFD{3}(\Ybar)$ to a type $DDA$ trimodule by tensoring with $\CFAAid$ along the $\sigma$ boundary;
\item
Tensor the type $A$ boundary of the resulting trimodule with $\CFDD(\Hselfgluer)$;
\item
Change the $\tau$ boundary to type $A$ by tensoring with $\CFAAid$, and then take the Hochschild homology with respect to the appropriate boundary components, resulting in a type $D$ module.
\item
Change this module to type $A$ by tensoring with $\CFAAid$ and tensor with the $\tau$ boundary of $\CFD{3}(\Y)$;
\item
The result is a type $DD$ bimodule with 16 generators. Note that it is still the case that $\alpha^\rho_2$ and $\alpha^\sigma_1$ represent fibers.
\end{itemize}
For $b>0$, the trivial bundle over $S_{g, b}$ can now be obtained easily by extending $S_{0, b}$ with $g$ copies of $S_{1, 2}$. For the case of $b = 0$, we simply extend the trivial bundle over $S_{g, 1}$ by capping off the boundary with the trivial bundle over $S_{0,1}$.

 \subsection{Nontrivial bundles}
In general the bundle associated to a vertex of $\Gamma$ is nontrivial, with a specified Euler number $e$. The Euler number of a circle bundle over a surface with boundary is well defined once a trivialization is chosen on the boundary. Choosing the trivialization over the boundary of an $S^1$-bundle is equivalent to choosing the $\alpha$ arcs to parametrize the boundary on a bordered Heegaard diagram of the total space. A trivialization over the boundary specifies two curves in the boundary $S^1 \times S^1$: a fiber $\gamma_f$, and a curve $\gamma_b$ meeting each fiber in one point. These in turn can specify a boundary parametrization by letting one be $\alpha_1$ and the other be $\alpha_2$.

\begin{figure}[htbp]
\begin{center}
\includegraphics[scale = .7]{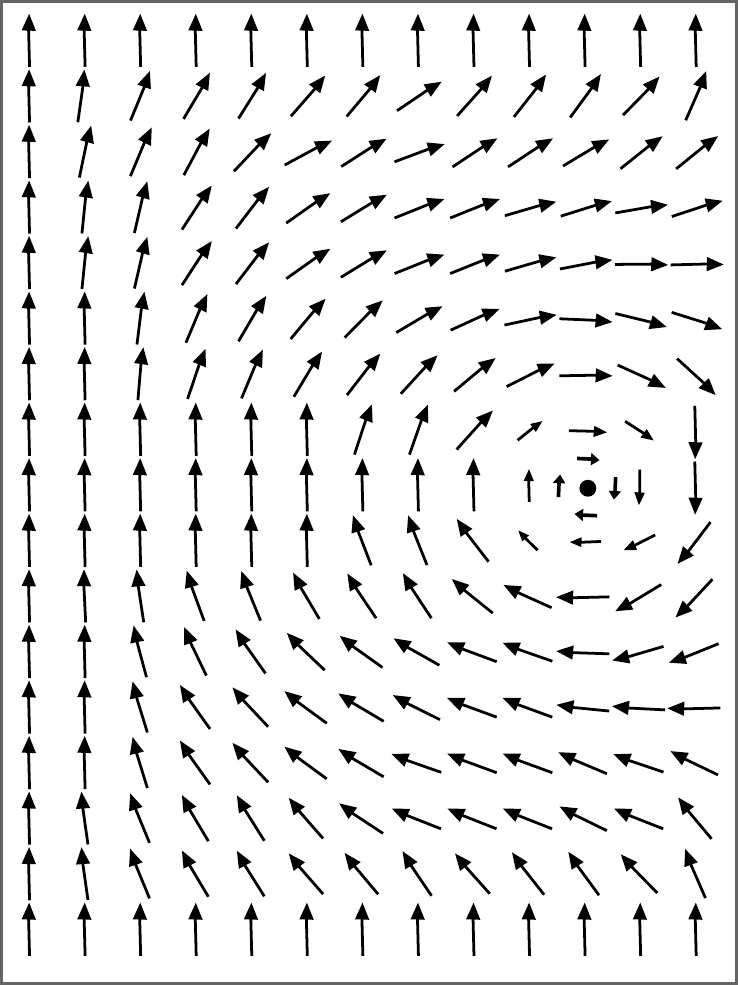}
\caption{A section of a $D^2$-bundle over the cylinder, where the top and bottom edges are identified. The section specifies a trivialization of the boundary on each end of the cylinder. The signed number of zeros indicates that the euler number is -1. The boundary $S^1$-bundle is equivalent to a the mapping cylinder of a Dehn twist.}
\label{euler_number}
\end{center}
\end{figure}

Consider changing the trivialization over one boundary component of an $S^1$-bundle from the trivialization on the left to the trivialization on the right in Figure \ref{euler_number}. On one hand, this change corresponds to gluing on the $S^1$-bundle over the cylinder which is the boundary of the $D^2$ bundle represented by Figure \ref{euler_number}. The figure shows a section of this $D_2$ bundle, which has a zero of sign $-1$. Since the Euler number can be defined as the signed number of zeros of such a section, it follows that attaching the corresponding $S^1$-bundle to a bundle over a surface decreases the Euler number by 1. On the other hand, this change of trivialization corresponds to the change in parametrization which fixes $\gamma_f$ and takes $\gamma_b$ to $\gamma_b \pm \gamma_f$ (depending on the orientation on $S^1 \times S^1$). This change is accomplished by attaching the mapping cylinder of a negative Dehn twist about $\gamma_f$.

In the same way, attaching a positive Dehn twist about the fiber $\gamma_f$ has the effect of increasing the Euler number of a circle bundle by 1. The bimodules for Dehn twists about $\alpha_1$ and $\alpha_2$ are known \cite[Section 10.2]{LOT:Bimodules}. By tensoring with enough of these bimodules we can obtain the bordered invariants for arbitrary $S^1$-bundles over arbitrary (oriented) surfaces with boundary.

\subsection{Combining vertices}

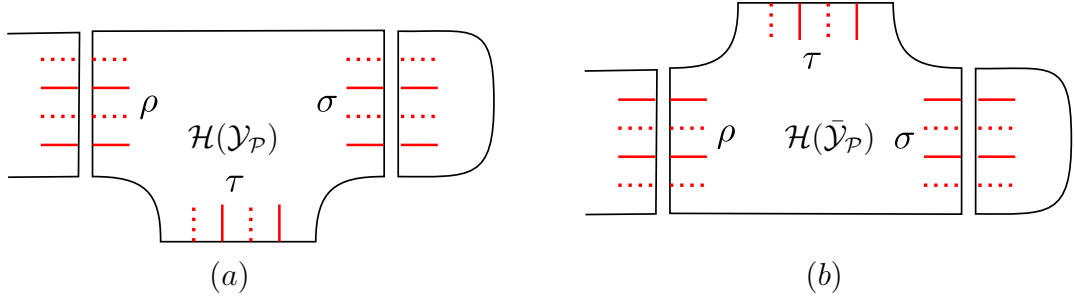
\begin{figure}[htbp]
\begin{center}
\begin{overpic}[scale = .8]{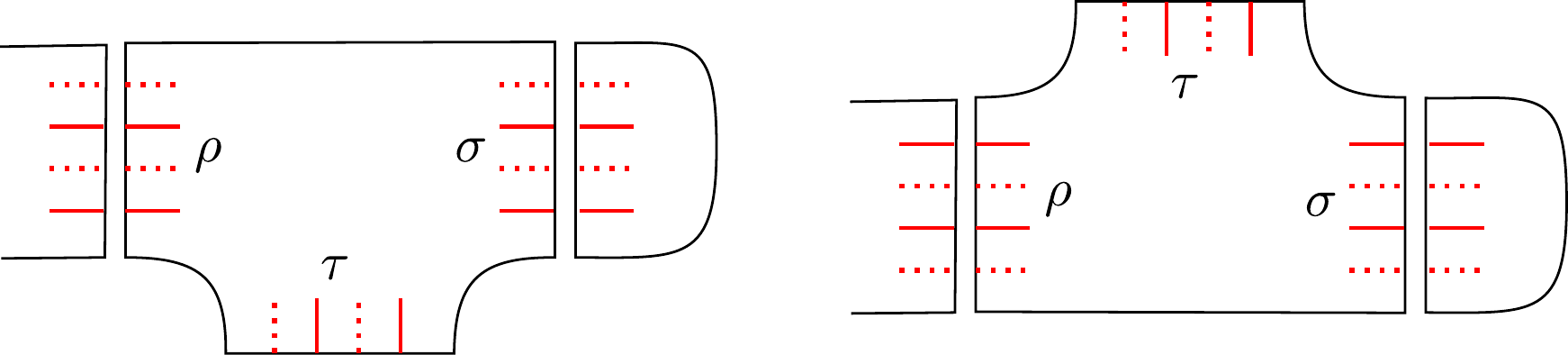}
\put(16, 9){ $\HeegDiag( \Y )$ }
\put(72, 9){ $\HeegDiag( \Ybar )$ }
\put(18, -4){ $(a)$ }
\put(74, -4){ $(b)$ }
\end{overpic}
\bigbreak

\caption{The fiber direction can be be changed by extending at a boundary with $\Y$ or $\Ybar$ and an appropriate solid torus. Dotted lines indicate $\alpha$ arcs which are fibers. The arrangement in $(a)$ changes the fiber from $\alpha_1$ to $\alpha_2$, while the arrangement in $(b)$ does the opposite.}
\label{change_fiber_direction}
\end{center}
\end{figure}

Once multimodules have been determined for each vertex of $\Gamma$, they can be combined according to the edges of $\Gamma$. If vertices $v_1$ and $v_2$ are connected by an edge, chose a boundary component of each circle bundle such that both boundaries have fiber $\alpha_1$ or both have fiber $\alpha_2$. Take the box tensor product (after changing one boundary component to type $A$) to compute the new multimodule. If there is no way to choose a boundary component with the desired $\alpha$ arc as fiber, the fiber direction can be changed as follows:
\begin{itemize}
\item
To change the fiber from $\alpha_1$ to $\alpha_2$, extend the bundle by $\Y$, attached along the $\rho$ boundary, with the $\sigma$ boundary capped off by a solid torus as in Figure \ref{change_fiber_direction}(a);
\item
To change the fiber from $\alpha_2$ to $\alpha_1$, extend the bundle by $\Ybar$, attached along the $\rho$ boundary, with the $\sigma$ boundary capped off by a solid torus as in Figure \ref{change_fiber_direction}(b).
\end{itemize}

For acyclic graphs any plumbing will work when combining vertices along an edge. In general, however, there is an additional consideration: edges are decorated by a sign, which distinguishes between two plumbing options. In terms of bordered Heegaard diagrams, the difference is between gluing two type $D$ boundaries with fiber $\alpha_1$ or gluing two boundaries with fiber $\alpha_2$. Suppose we orient each boundary component so that the positive fiber direction is to the left of the positive base direction at a fiber-base intersection. Then a type $D$ boundary with $\alpha_1$ a fiber has oriented fiber $-\alpha_1$ and oriented base $+\alpha_2$. Since gluing type $D$ boundaries glues $\alpha_1$ to $-\alpha_2$, this corresponds to the map $\scriptsize \left( \begin{array}{cc}  0 & 1 \\  1 & 0  \end{array} \right) $ in the standard $\{ \text{base, fiber} \}$ basis. That is, gluing two boundaries with $\alpha_1$ fibers corresponds to a $+$ edge. A type $D$ boundary with fiber $\alpha_2$ has oriented fiber $+\alpha_2$ and oriented base $+\alpha_1$, so gluing two of these boundaries corresponds to the map $\scriptsize \left( \begin{array}{cc}  0 & -1 \\  -1 & 0  \end{array} \right) $.

Once the bundles of two adjacent vertices have been plumbed, the result is no longer an $S^1$-bundle. However, continue to keep track of which $\alpha$ arc is the ``fiber" at each boundary component. Repeat the process above to add on successive vertices. If at any point an edge connects to a vertex that has already been incorporated, insert the bimodule $\CFDD(\Hselfgluer)$ and take the appropriate Hochshchild homology instead of a tensor product.

\subsection{Example computations}

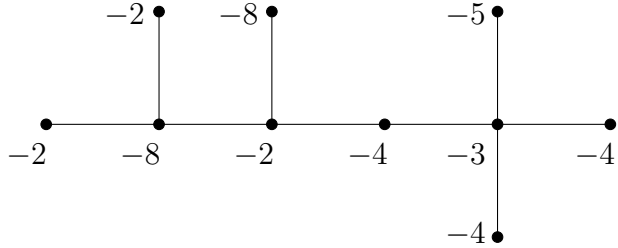
\begin{figure}[htbp]
\begin{center}
\begin{tikzpicture}

\tikzstyle{every node}=[draw,circle,fill=white,minimum size=4pt,
                            inner sep=0pt]

\node at (0,0) [circle, fill=black] {};
\node at (1.5,0) [circle, fill=black] {};
\node at (1.5,1.5) [circle, fill=black] {};
\node at (3,0) [circle, fill=black] {};
\node at (3,1.5) [circle, fill=black] {};
\node at (4.5,0) [circle, fill=black] {};
\node at (6,0) [circle, fill=black] {};
\node at (6,1.5) [circle, fill=black] {};
\node at (7.5,0) [circle, fill=black] {};
\node at (6,-1.5) [circle, fill=black] {};

\draw (0,0) -- (7.5,0); \draw (1.5,0) -- (1.5,1.5); \draw (3,0) -- (3,1.5); \draw (6,-1.5) -- (6,1.5);

\put(-15,-15){$-2$}
\put(28,-15){$-8$}
\put(22,38){$-2$}
\put(71,-15){$-2$}
\put(65,38){$-8$}
\put(114,-15){$-4$}
\put(151,-15){$-3$}
\put(151,38){$-5$}
\put(200,-15){$-4$}
\put(151,-45){$-4$}

\end{tikzpicture}

\caption{Plumbing graph for a graph manifold with $\text{rk}(\HFhat) = 213,312$. The weights on the vertices correspond to Euler numbers; the genus is zero for every vertex and we omit it from the notation.}
\label{huge_graph}
\end{center}
\end{figure}

The author has implemented a program\footnote{Available at \url{http://math.columbia.edu/~jhansel/graph_manifolds_program.html}} using the techniques described above to compute the total rank of $\HFhat$ of a closed graph manifold, or the bordered invariant of a graph manifold with boundary, from a plumbing graph. It can be used, for example, to see that the rank of $\HFhat$ of the manifold represented by the negative definite plumbing tree in Figure \ref{huge_graph} is 213,312. It is easy to compute $| H_1 |$ from the plumbing graph and see that this manifold is an $L$-space. This is as expected; the fact that this plumbing graph corresponds to an $L$-space follows from \cite[Theorem C]{Mauricio}.

With this algorithm, we can quickly run computations for large sets of graph manifolds and check, for instance, which are $L$-spaces. Consider as an example the plumbing graph $\Gamma$ below, with weights in the range $-5 \le a, b \le 5$ and $-5 \le c,d,e,f \le -2$ (the bound of $-2$ on the weights of the outer vertices is so that we only consider graphs in normal form, in the notation of \cite{Neumann}).
\begin{center}
\begin{tikzpicture}
\tikzstyle{every node}=[draw,circle,fill=white,minimum size=4pt,
                            inner sep=0pt]
\node at (0,0) [circle, fill=black] {};
\node at (1,0) [circle, fill=black] {};
\node at (-.5,.87) [circle, fill=black] {};
\node at (-.5,-.87) [circle, fill=black] {};
\node at (1.5,.87) [circle, fill=black] {};
\node at (1.5,-.87) [circle, fill=black] {};
\draw (0,0) -- (1,0); \draw (-.5,-.87) -- (0,0) -- (-.5,.87); \draw (1.5,-.87) -- (1,0) -- (1.5,.87);
\put(-11,-1){$a$}
\put(33,-1){$b$}
\put(-24,-26){$c$}
\put(-24,20){$d$}
\put(46,-26){$f$}
\put(46,20){$e$}
\end{tikzpicture}
\end{center}
There are 6106 distinct graphs of this form. Of the corresponding 3-manifolds, 5643 are $L$-spaces. Some of these trees are negative definite, but most are not. To the author's knowledge, there is currently no other way to compute $\HFhat$ for these non-definite examples. Results for a few examples are in Figure \ref{examples}.

\begin{figure}[htbp]
\begin{center}

\bigbreak
\bigbreak
\bigbreak

\begin{tikzpicture}
\tikzstyle{every node}=[draw,circle,fill=white,minimum size=4pt,
                            inner sep=0pt]
\node at (0,0) [circle, fill=black] {};
\node at (1,0) [circle, fill=black] {};
\node at (-.5,.87) [circle, fill=black] {};
\node at (-.5,-.87) [circle, fill=black] {};
\node at (1.5,.87) [circle, fill=black] {};
\node at (1.5,-.87) [circle, fill=black] {};
\draw (0,0) -- (1,0); \draw (-.5,-.87) -- (0,0) -- (-.5,.87); \draw (1.5,-.87) -- (1,0) -- (1.5,.87);
\put(-12,-1){$5$}
\put(35,-1){$5$}
\put(-34,-26){$-5$}
\put(-34,20){$-5$}
\put(46,-26){$-5$}
\put(46,20){$-5$}

\end{tikzpicture} \hspace{1 in}
\begin{tikzpicture}
\tikzstyle{every node}=[draw,circle,fill=white,minimum size=4pt,
                            inner sep=0pt]
\node at (0,0) [circle, fill=black] {};
\node at (1,0) [circle, fill=black] {};
\node at (-.5,.87) [circle, fill=black] {};
\node at (-.5,-.87) [circle, fill=black] {};
\node at (1.5,.87) [circle, fill=black] {};
\node at (1.5,-.87) [circle, fill=black] {};
\draw (0,0) -- (1,0); \draw (-.5,-.87) -- (0,0) -- (-.5,.87); \draw (1.5,-.87) -- (1,0) -- (1.5,.87);
\put(-12,-1){$0$}
\put(35,-1){$0$}
\put(-34,-26){$-5$}
\put(-34,20){$-4$}
\put(46,-26){$-5$}
\put(46,20){$-3$}
\end{tikzpicture} \hspace{1 in}
\begin{tikzpicture}
\tikzstyle{every node}=[draw,circle,fill=white,minimum size=4pt,
                            inner sep=0pt]
\node at (0,0) [circle, fill=black] {};
\node at (1,0) [circle, fill=black] {};
\node at (-.5,.87) [circle, fill=black] {};
\node at (-.5,-.87) [circle, fill=black] {};
\node at (1.5,.87) [circle, fill=black] {};
\node at (1.5,-.87) [circle, fill=black] {};
\draw (0,0) -- (1,0); \draw (-.5,-.87) -- (0,0) -- (-.5,.87); \draw (1.5,-.87) -- (1,0) -- (1.5,.87);
\put(-21,-1){$-4$}
\put(33,-1){$-1$}
\put(-34,-26){$-5$}
\put(-34,20){$-5$}
\put(46,-26){$-5$}
\put(46,20){$-2$}
\end{tikzpicture}

\bigbreak

$\begin{array}{c}
 \text{rk}( \HFhat ) = 17,\!600   \\
 | H_1 | = 17,\!600
\end{array}$ \hspace{.6 in}
$\begin{array}{c}
 \text{rk}( \HFhat ) = 230   \\
 | H_1 | = 228  
\end{array}$ \hspace{.6 in}
$\begin{array}{c}
 \text{rk}( \HFhat ) = 72   \\
 | H_1 | = 20
\end{array}$

\caption{The manifold corresponding to the graph on the left has the largest $\HFhat$ of the 6106 examples tested, and it is an $L$-space. The graph in the middle gives the smallest difference between $rk(\HFhat)$ and $ | H_1 |$ possible for a non $L$-space, and the third gives the largest difference among this set of examples.}
\label{examples}
\end{center}
\end{figure}
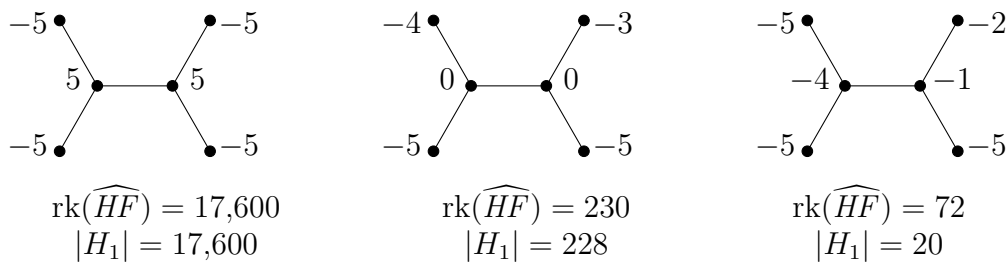

Our final example is the manifold $\Sigma \times S^1$, where $\Sigma$ is the surface of genus two. This manifold can be represented by a plumbing graph with just one vertex and no edges. The vertex carries the weights 2 and 0 for the genus and Euler number, respectively. Evaluating the rank of $\HFhat$ from this graph gives 24, which agrees with the result in \cite{JabukaMark}.

\bibliographystyle{alpha}
\bibliography{HFbibliography}

\end{document}